\documentclass[11pt]{article}
\usepackage[dvipsnames]{xcolor}

\usepackage{amsmath}
\usepackage{amssymb}
\usepackage{latexsym}
\usepackage{enumitem}
\usepackage{mathrsfs}
\usepackage{comment}
\usepackage{color}
\usepackage{bbm}

\usepackage{hyperref}
\hypersetup{
    colorlinks=true,
    linkcolor=MidnightBlue,   
    filecolor=MidnightBlue,
    urlcolor=MidnightBlue,     
    citecolor=MidnightBlue,
    pdftitle={Set Valued HJB Equations},
    pdfsubject={Mathematics},
    bookmarksnumbered=true,
    bookmarksopen=true,
    bookmarksopenlevel=1,
    pdfstartview=Fit,
    pdfpagemode=UseOutlines,
}

\newtheorem{assumption}{Assumption}

\def\qed{ \ \vrule width.2cm height.2cm depth0cm\smallskip}
\newenvironment{proof}{\noindent {\bf Proof.\/}}{$\qed$\vskip 0.1in}

\newcommand{\ol}{\overline}
\newcommand{\ul}{\underline}

\newcommand{\ba}{\begin{array}}
\newcommand{\ea}{\end{array}}
\newcommand{\be}{\begin{equation}}
\newcommand{\ee}{\end{equation}}
\newcommand{\bea}{\begin{eqnarray}}
\newcommand{\eea}{\end{eqnarray}}
\newcommand{\beaa}{\begin{eqnarray*}}
\newcommand{\eeaa}{\end{eqnarray*}}

\def\dbD{\mathbb{D}}
\def\dbE{\mathbb{E}}
\def\dbF{\mathbb{F}}
\def\dbG{\mathbb{G}}

\def\dbL{\mathbb{L}}

\def\dbP{\mathbb{P}}
\def\dbR{\mathbb{R}}

\def\dbT{\mathbb{T}}

\def\dbU{\mathbb{U}}
\def\dbV{\mathbb{V}}

\def\dbY{\mathbb{Y}}
\def\dbZ{\mathbb{Z}}

\def\a{\alpha}
\def\b{\beta}
\def\g{\gamma}
\def\d{\delta}
\def\e{\varepsilon}
\def\z{\zeta}
\def\k{\kappa}
\def\l{\lambda}

\def\si{\sigma}
\def\t{\tau}
\def\f{\varphi}
\def\th{\theta}
\def\o{\omega}

\def\G{\Gamma}
\def\D{\Delta}

\def\L{\Lambda}

\def\O{\Omega}
\def\cA{{\cal A}}

\def\cD{{\cal D}}

\def\cF{{\cal F}}

\def\cK{{\cal K}}
\def\cL{{\cal L}}

\def\cN{{\cal N}}

\def\cY{{\cal Y}}

\def\no{\noindent}

\def\ss{\smallskip}
\def\ms{\medskip}
\def\bs{\bigskip}
\def\q{\quad}
\def\qq{\qquad}

\def\pa{\partial}

\def\cd{\cdot}
\def\cds{\cdots}

\def\td{\nabla}

\def\tr{\hbox{\rm tr}}

\def\qed{ \hfill \vrule width.25cm height.25cm depth0cm\smallskip}

\newcommand{\basa}{\begin{assumption}}
\newcommand{\easa}{\end{assumption}}

\newcommand{\bas}{\begin{assum}}
\newcommand{\eas}{\end{assum}}

\def\limsup{\mathop{\overline{\rm lim}}}
\def\liminf{\mathop{\underline{\rm lim}}}

\def\esup{\mathop{\rm ess\;sup}}

\def\pa{\partial}

 \def\cd{\cdot}
\def\cds{\cdots}

\def\tr{\hbox{\rm tr$\,$}}

\def\dis{\displaystyle}

\def\wh{\widehat}

\def\bx{{\bf x}}

\def\1{{\bf 1}}
\def\by{{\bf y}}
\def\bn{{\bf n}}

\def\br{{\bf r}}

\def\:{\!:\!}

\def \proof{{\noindent \bf Proof\quad}}

\def\U{{\Upsilon}}

at 9pt

\newtheorem{thm}{Theorem}[section]
\newtheorem{lem}[thm]{Lemma}

\newtheorem{prop}[thm]{Proposition}
\newtheorem{rem}[thm]{Remark}
\newtheorem{eg}[thm]{Example}
\newtheorem{defn}[thm]{Definition}
\newtheorem{assum}[thm]{Assumption}

\numberwithin{equation}{section}

\newcommand{\ind}[1]{{\bf 1}_{\left\{ {#1} \right\}} }

\newcommand{\mb}[1]{\mathbb{#1}}

\begin{document}

\title{\bf Set-valued Hamilton-Jacobi-Bellman Equations} \author{Melih
  \.{I}\c{s}eri\footnote{Department of Mathematics, University of Michigan,
    United States, iseri@umich.edu.}  \quad Jianfeng
  Zhang\footnote{Department of Mathematics, University of Southern
    California, United States, jianfenz@usc.edu. This author is
    supported in part by NSF grant DMS-2205972.}}

\date{March 31, 2024}

\maketitle

  \begin{abstract}
    Building upon the dynamic programming principle for set-valued functions arising from many applications, in this paper we propose a new notion of set-valued PDEs. The key component in the theory is a set-valued It\^{o} formula, characterizing the flows on the surface of the dynamic  sets. In the contexts of multivariate control problems, we establish the wellposedness of the set-valued HJB equations, which extends the standard HJB equations in the scalar case to the multivariate case. As an application, a moving scalarization for certain time inconsistent problems is constructed by using the classical solution of the set-valued HJB equation.  
  \end{abstract}

  \vspace{3em}
  
\no{\bf MSC2020.} 49L12, 49J53,  60H30, 49J20, 47H04

\vspace{3mm}
\no{\bf Keywords.}  Set valued functions, dynamic programming principle, set valued It\^{o} formula, set valued HJB equation, multivariate control problems, time inconsistent problems, surface evolution equations

\vfill\eject

\section{Introduction}
\label{sect-Introduction}
\setcounter{equation}{0}
In this paper we consider set-valued functions taking the form:
\beaa
\label{setvalued}
\dbV: [0, T]\times \dbR^d \to 2^{\dbR^m}.
\eeaa
That is, for each $(t,x)\in [0, T]\times \dbR^d$, the value $\dbV(t,x)$ is a subset of $\dbR^m$ satisfying appropriate properties. 
Such set-valued functions, or their variants, have appeared in many applications, for example, stochastic
viability problems (cf. Aubin-Da Prato \cite{ADP}), multivariate
super-hedging problems (cf. Kabanov \cite{Kabanov} and Bouchard-Touzi \cite{BT2000}), multivariate
dynamic risk measures (cf. Feinstein-Rudloff \cite{FR2015}),  time inconsistent optimization problems (cf. Karman-Ma-Zhang \cite{KMZ}), stochastic
target problems (cf. Soner-Touzi \cite{ST2002, ST2002b, ST2003}), 
and recently, nonzero sum (mean field) games
with multiple equilibria (Feinstein-Rudloff-Zhang \cite{FRZ2022}
and \.I\c{s}eri-Zhang \cite{IZ2021}). Many of these problems were considered non-standard or even ill-posed in the literature, and overall we lacked convenient mathematical tools to treat them. When viewed as set values, however, their value functions (named set value functions\footnote{We use this for the value functions from the applications, to distinguish from general set-valued functions.} in this paper) enjoy many nice properties as the value function of standard control problems, in particular the crucial  {\it dynamic programming principle}, or say the  time consistency. Notice that, for a standard control problem, the combination of the dynamic programming principle and the It\^{o} formula leads to the popular PDE approach.  Then a natural question is:
\bea
\label{question}
\mbox{\it Can we characterize these set value functions via set-valued PDEs?} 
\eea

\no This is exactly the goal of the present paper: to introduce the PDE approach and hence recover the standard language for these challenging problems. To be precise, the main contributions of this paper are as follows:

\begin{itemize}
\item{} Introduce  derivatives for set-valued functions and establish the set-valued It\^{o} formula. 

\item{} Propose a notion of set-valued PDEs, and establish its wellposedness in the contexts of multivariate stochastic control problems. 

\item{} As an important application, construct a so called moving scalarization for a time inconsistent problem by using the classical solution of the corresponding set-valued PDE.
\end{itemize}
\no We hope this paper serves as the first step of our long term project on providing a convenient tool and systematic study for multivariate problems, including games.

Our first main result  is the {\it set-valued It\^{o} formula}, which roughly reads:  
\bea
\label{Ito0}
\qq d\dbV(t, X_t) = \Big[\pa_t \dbV + \pa_x\dbV \cd b + \frac{1}{2}\tr (\si^\top\pa_{xx}\dbV  \si) - \cK_\dbV \zeta + \xi\Big] dt + \Big[\pa_x\dbV \si + \zeta\Big]dB_t,
\eea
where  $d X_t = b_t dt + \si_t dB_t$ is an arbitrary diffusion. We refer to Theorem \ref{thm-Ito} below for the precise meaning of the above formula. In particular, $\pa_t \dbV, \pa_x\dbV, \pa_{xx}\dbV$ are derivatives of $\dbV$ defined on the graph $\dbG_\dbV$,  which is the graph of the boundary surface $\dbV_b$ and consists of all points $(t,x, y)$ where $y$ lies on $\dbV_b(t,x)$. The essence of the It\^o formula is to
  characterize flows on the boundary surface. Given the surface's
  invariance under tangential deformations, a key feature of the set
  valued It\^o formula is the inclusion of arbitrary driving forces
  $\xi$ and $\zeta$ on $\dbG_\dbV$, which take values in the tangent
  space. This, along with the appropriate correction term
  $\cK_\dbV \zeta$, ensures that the flows are not pushed away from
  the boundary surface.

Our set-valued HJB equation takes the following form: for some
  Hamiltonian function $h_\dbV$ and with appropriate terminal condition,
  \bea
\label{HJB0}
\hspace{1.5em}\sup_{a, \zeta} ~\bn_\dbV(t,x,y) \cd \Big[\pa_t \dbV(t,x,y) + h_\dbV(t,x, y, \pa_x\dbV, \pa_{xx}\dbV, a, \zeta)\Big] =0,\ \ (t,x, y) \in \dbG_\dbV,
\eea
where $a$ takes values in a control set,  $\zeta$ takes values in the tangent space, and $\bn_\dbV$ is the unit outward normal vector. This is derived by applying the above It\^o formula on the dynamic programming principle for the set value function of the multivariate control problem.
As we see, the introduction of $\zeta$ (and $\xi$) in \eqref{Ito0} is crucial. We note that $\xi$ disappears in the equation since $\bn_\dbV \cd \xi = 0$. However, $\cK_\dbV \zeta$ is nonlinear in $\zeta$ and thus $\bn_\dbV \cd \cK_\dbV \zeta$ is an important  component in the equation.  The equation \eqref{HJB0} can be rewritten equivalently in terms of the signed distance function $\br_\dbV$, see \eqref{HJBr} below.  We emphasize that $\bn_\dbV$ is part of the solution here and the equation is satisfied only on the graph $\dbG_\dbV$, so the wellposedness of \eqref{HJB0} has a completely different nature than that of standard PDEs.

In the scalar case: $m=1$, we can easily see that $\dbV(t,x) = [\ul v(t,x), \ol v(t,x)]$, where $\ul v$ and $\ol v$ are the value functions of the standard minimization and maximization problems, respectively. In this case,  $\bn_\dbV = 1$ or $-1$, and the tangent space is degenerate and thus $\zeta =0$. Then \eqref{HJB0} reduces exactly back to the standard HJB equations for $\ul v$ and $\ol v$. So our set-valued HJB equation is a natural extension of the standard HJB equation to the multivariate setting.  Indeed, \eqref{HJB0} means the optimization over each direction $\bn_\dbV(t,x,y)$. Moreover, note that $\ul v$ and $\ol v$ are the boundaries of $\dbV$ in this case, namely $\dbV_b(t,x) = \{\ul v(t,x), \ol v(t,x)\}$, which  inspires us to focus on the boundary surface $\dbV_b$ instead of on the whole set $\dbV$. We would like to point that, again since $\xi=0$, $\zeta =0$, in this case \eqref{Ito0} also reduces back to the standard It\^{o} formula for $\ul v(t, X_t)$ and $\ol v(t, X_t)$.

Our main result of the paper is that the dynamic set value function of
the multivariate stochastic control problem is the unique classical
solution of the set-valued HJB equation \eqref{HJB0}, provided $\dbV$
has sufficient regularity. We thus obtain the positive answer to our question \eqref{question} in this setting, which further opens the door to the PDE approach  for more general multivariate problems. Such PDE characterization helps to understand better the structure and the nice properties of the dynamic set value function. In particular, it helps to construct (approximate) optimal controls with certain Markovian structure. Indeed, when $\dbV$ is smooth, as in
standard verification theorem we may use the optimal arguments
$(a^*, \zeta^*)$ of the Hamiltonian in \eqref{HJB0} to construct an
optimal control for a scalarized optimization problem, as we will explain in a later paragraph. 

It should be noted that, under natural Lipschitz conditions our set values turn out to be convex and the set value function is uniformly continuous in $(t,x)$, both are basic properties important in theory and in applications. However, we emphasize that our set-valued It\^{o} formula does not require the convexity, and our arguments for the wellposedness of the set-valued HJB equations does not rely on this convexity either.

As an important application of our wellposedness result, we construct a moving scalarization for some time inconsistent problems, proposed by \cite{KMZ}\footnote{In \cite{KMZ} it's called dynamic utility function, instead of moving scalarization.} and Feinstein-Rudloff \cite{FR2022}. Note that we are in the multivariate setting and in general it is not feasible to optimize the multiple objects simultaneously. In practice quite often one considers a scalarized optimization problem: $\max_{y\in \dbV(0, x_0)} \f(y)$, where $\f: \dbR^m\to \dbR$. This scalarized problem, unfortunately, is typically time inconsistent. The idea of a moving scalarization is to find a dynamic scalarization function $\Phi(t, X_{[0, t]}; y)$, with $\Phi(0,x_0; y) = \f(y)$, such that the dynamic problem  $\max_{y\in \dbV(t, X_t)} \Phi(t, X_{[0, t]}; y)$ becomes time consistent.  In Section \ref{sect-moving} below we shall investigate this interesting application.  In particular, we shall construct a moving scalarization for the mean variance problem explicitly. 

At this point we would like to mention that the present paper considers classical solutions only. In particular, this requires that the set value $\dbV(t,x)$ is non degenerate and its boundary $\dbV_b(t,x)$ is a smooth $m-1$ dimensional manifold, namely the co-dimension is $1$. 
It is our strong interest to remove these constraints and study viscosity solutions of more general set-valued PDEs, thereby broadening the applicability of the theory. We shall leave this important question to future research.  
\ss

{\bf Some related literature.} There is a large literature on set-valued analysis, we refer to  the books  Aubin-Frankowska \cite{AF2009},   Kisielewicz \cite{Kisielewicz1, Kisielewicz2}, and the reference therein.  However, our approach is completely different from those in the set-valued analysis. First, the input $(t,x)$ of our  set-valued functions $\dbV$ is vector valued, while in  the standard literature it is typically  set-valued. Next, we focus on the dynamics of the boundary surface, rather than the dynamics of the whole set. Roughly speaking, we focus only on those special selectors whose flow remains on the boundary. These selectors have nice properties and are sufficient to characterize the whole sets. Moreover, the boundary is essentially the frontier which has intrinsic optimality and thus is also important from practical point of view.  We shall mention the recent paper Ararat-Ma-Wu \cite{AMW} on set-valued backward SDEs, which is highly relevant to our paper. Given our set value function $\dbV(t,x)$ and a state process $X$ (e.g. a Brownian motion), we may introduce a set-valued process $\dbY_t := \dbV(t, X_t)$. In spirit the process $\dbY$ should satisfy a set-valued backward SDE. However, besides that we employ completely different approaches, except in some simple cases the set-valued process $\dbY$  does not satisfy the equation in \cite{AMW}. That is, the objectives of the two works are different.  We should mention that the applications mentioned in the beginning of this introduction fall into our framework, although technically our current results do not cover many of them (which we intend to  study in our future research). 

Our approach is strongly motivated by the studies on surface evolution equations, see e.g. Sethian \cite{S1985}, Evans-Spruck \cite{ES1991}, Soner
\cite{S1993}, Barles-Soner-Souganidis \cite{BSS1993}, the monograph Giga \cite{G2006}, and the references therein. These equations arise in various
applications such as evolutions of phase boundaries, crystal growths, image processing, and mean-curvature flows, to mention a few. These works consider the dynamics of a set-valued function $\dbV(t)$,  more precisely the boundary $\dbV_b(t)$, without the state variable $x$. In our terms, roughly speaking these works study first order set-valued ODEs, while we extend to second order set-valued PDEs. In particular, the set-valued It\^{o} formula is not involved there. Another difference is, due to the nature of different applications, they study forward equations with initial conditions while we study backward problems with terminal conditions. This difference would be crucial when one concerns path dependent setting (not covered in this paper), where one cannot do the time change freely due to the intrinsic adaptedness requirement.

Furthermore, within the surface evolution literature,  our work is  closely related to Soner-Touzi  \cite{ST2003} which studies stochastic target problems by using mean curvature type geometric flows. In our contexts, their approach amounts to studying the following set value function via its signed distance function $\br_{\wh \dbV}$:
\begin{equation*}
\label{hatV0}
\wh \dbV(t): = \big\{(x, y): x\in \dbR^d, y\in \dbV(t,x)\big\},~\mbox{and thus}~ \wh \dbV_b(t): = \big\{(x, y): x\in \dbR^d, y\in \dbV_b(t,x)\big\}.  
\end{equation*}
Clearly $\wh \dbV$ and $\dbV$ are equivalent, with the same graph: $\dbG_{\wh \dbV} = \dbG_\dbV$.  The major difference here is that, while $\br_\dbV$ and $\br_{\wh\dbV}$ agree on the graph (both are $0$ by definition), their derivatives are different on the graph, and consequently, the equation derived in \cite{ST2003} is different from our set-valued HJB equation \eqref{HJB0}. In particular, in the scalar case: $m=1$, as mentioned \eqref{HJB0} reduces back to the standard HJB equations, but the equation for $\br_{\wh\dbV}$ appears in a quite different form.
Moreover, the normal vector $\bn_{\wh\dbV}$ of $\wh \dbV$ is also different from $\bn_\dbV$, and does not serve as a moving scalarization as we discussed. We shall provide more detailed discussions in Section \ref{sect-ST} below.

Finally, we remark that there are some very interesting studies on (possibly discontinuous) viscosity solutions along this line, see e.g. \cite{BSS1993}, Chen-Giga-Goto \cite{CGG1991}, Soner-Touzi \cite{ST2002, ST2002a, ST2002b}. It will be interesting to explore these ideas in our setting.

The rest of the paper is organized as follows. In Section \ref{sect-derivative} we introduce the setting and define the intrinsic derivatives of set-valued functions. In Section \ref{sect-Ito} we prove the crucial set-valued It\^{o} formula. In Section \ref{sect-control} we present the multivariate control problem. In Section \ref{sect-HJB} we introduce the set-valued HJB equation and show that the value function of the multivariate  control problem is a classical solution, and the uniqueness of classical solutions is established in Section \ref{sect-uniqueness}. Section \ref{sect-moving} is devoted to the application of moving scalarization. In particular we solve it explicitly for the mean variance problem.  In Section \ref{sect-further} we offer further discussions, including an extension to the case that the terminal condition is non-degenerate, and comparisons with \cite{AMW}  and \cite{ST2003}. Finally we complete some technical proofs in Appendix. 

\vskip 2mm

\no{\bf Notation.} For the convenience of the readers, we list some frequently used notation here.
\begin{itemize}
\item{}  $\dbD$: subsets of $\dbR^m$;

\item{} $\dbV: [0, T]\times \dbR^d\to 2^{\dbR^m}$: set-valued functions, where in particular $\dbV(t,x)\subset \dbR^m$ is a set;
 
\item{} $\dbV_b$ and $\dbV_o$: the boundary and  interior of  $\dbV$, respectively;

\item{} $ \dbG_\dbV$: the graph of $\dbV_b$, see \eqref{GV} and \eqref{GVtx};

\item{}  $\br$:  the signed distance function, see \eqref{rD};

\item{}  $\bn$: the outward unit normal vector, see \eqref{nD};

\item{} $\dbT$: the tangent space, see \eqref{TD};

\item{} $\pi$: the projection onto the boundary, see \eqref{piD};

\item{} $\mathbf{I}_m$: the $m\times m$ identity matrix; and $\mathbf{I}_m - \bn\bn^\top$ is the projection operator onto $\dbT$;

\item{} $\td \wh f$: the standard derivatives of a function $\wh f$;

\item{}  $\pa_\cdot f$: intrinsic derivatives of  $f:\dbG_\dbV\to\dbR$,\\ see \eqref{payfhat}, \eqref{paxfhat}, \eqref{patxV}, \eqref{paxfpaxV},  \eqref{paxxV}, and \eqref{pan};

\item{}  $X$ and $(Y,Z)$: solutions to  SDEs and BSDEs, respectively, see e.g. \eqref{FBSDE};

\item{} $\U$: forward dynamics typically on $\dbV_b$, see e.g. \eqref{Upsilon};

\item{} $(\xi,\zeta)$: vector fields taking values on the tangent space $\dbT$, see e.g. \eqref{Upsilon};

\item{} $\cK_\dbV^\sigma\zeta$: correction term of $\zeta$ in the It\^o formula, see \eqref{kappa};

\item{}    $\cL^{b,\sigma}$: differential operator in the  It\^o formula, see \eqref{kappa};

\item{} $(h_\dbV^0,h_\dbV,H_\dbV)$: the related Hamiltonians, see \eqref{H};

\item{} $\cL$: the differential operator of the HJB equation, see \eqref{HJB}.
\end{itemize}

\section{Intrinsic derivatives of set-valued functions}
\label{sect-derivative}
\setcounter{equation}{0}

Throughout the paper all vectors are viewed as column vectors, $\cd$ denotes the inner product, and $^\top, ^c$ denote the transpose and complement, respectively. We denote by $\td$ the gradient operator, and for a function $f: \dbR^d\times \dbR^m\to \dbR$, we take the convention that the second derivative $\td_{xy} f(x, y) := [\pa_{x_1 y} f, \cds, \pa_{x_d y} f] \in \dbR^{m\times d}$.

\subsection{Some basic materials}
In this subsection we present some basic materials in geometry, which will be the starting point of our set-valued functions in this paper.

Let $\cD^m_0\subset 2^{\dbR^m}$ denote the space of closed sets $\dbD$ in $\dbR^m$, 
and denote by $\dbD_o$, $\dbD_b$, $\dbD^c$, the interior, the boundary\footnote{It may seem more natural to use $\pa \dbD$ to denote the boundary, however, in this paper we reserve the notation $\pa$ for the intrinsic (partial) derivatives.}, and the complement of $\dbD$, respectively.
Define the Hausdorff distance:
\bea
\label{dDD}
~~~~~~~~d(\dbD, \tilde \dbD) := \big(\sup_{y\in \dbD} d(y, \tilde \dbD)\big) \vee \big(\sup_{\tilde y\in \tilde \dbD} d(\tilde y, \dbD)\big),~\mbox{where}~  d(y, \dbD) := \inf_{y'\in \dbD} |y-y'|,~\forall y\in \dbR^m.
\eea
The following signed distance function of $\dbD$ plays a central role in our analysis: 
\bea
\label{rD}
\br_{\dbD}(y) := \left\{\ba{lll} d(y, \dbD_b),\q~ y\in \dbD^c;\\ -d(y, \dbD_b),\q y\in \dbD.\ea\right. 
\eea
It is obvious that
\beaa
\dbD_o = \{y\in \dbR^m: \br_\dbD(y) <0\},~ \dbD_b = \{y\in \dbR^m: \br_\dbD(y) = 0\},~ \dbD^c = \{y\in \dbR^m: \br_\dbD(y) > 0\}.
\eeaa

We next let $\cD^m_2$ denote the space of $\dbD\in \cD^m_0$ such that $\br_{\dbD}$ is twice continuously differentiable with bounded derivatives on $O_\e(\dbD_b) := \{y\in \dbR^m: |\br_{\dbD}(y)|<\e\}$ for some $\e>0$. We remark that the boundary $\dbD_b$ is a manifold without boundary, as regular as $\br_\dbD$.
For each $y\in \dbD_b$, let $\bn_\dbD(y)\in \dbR^m$ denote the outward unit normal vector at $y$.
Moreover, for any $y\in O_\e(\dbD_b)$, let $\pi_{\dbD}(y)$ denote the unique projection of $y$ on $\dbD_b$, i.e. $\pi_{\dbD}(y)\in \dbD_b$ satisfies:
\bea
\label{piD}
y = \pi_{\dbD}(y) + \br_\dbD(y) \bn_\dbD(\pi_{\dbD}(y)),\q y\in O_\e(\dbD_b).
\eea
We refer to \cite{KP1981} and \cite[Corollary 3.4.5]{CS2004} for
  discussions on the connections between the regularity of the signed
  distance function $\br_\dbD$, the surface regularity of $\dbD_b$,
  and the projection $\pi_\dbD$.
It is clear that: 
  \bea
\label{nD}
\bn_\dbD(\pi_\dbD(y)) = \td_y \br_{\dbD}(y),\q\mbox{and}\q |\td_y \br_{\dbD}(y)|=1,\q y\in O_\e(\dbD_b).
\eea

For any $y\in \dbD_b$, let $\dbT_\dbD(y)$ denote the tangent space:
\bea
\label{TD}
 \dbT_\dbD(y):= \big\{ \xi\in \dbR^m: \xi\cd \bn_\dbD(y) =0\big\},\q y\in \dbD_b.
\eea
For a differentiable function $f: \dbD_b\to \dbR$, we define its intrinsic derivative $\pa_y f(y)\in \dbT_\dbD(y)$:
\bea
\label{payf}
\ \ \lim_{\e\to 0} {f(\th(\e)) - f(y)\over \e} = \pa_y f(y) \cd \th'(0),~\mbox{for any smooth curve $\th: \dbR \to \dbD_b$ with $\th(0)=y$}.\hspace{0em}
\eea
Alternatively, for any smooth extension $\wh f: \dbR^m\to \dbR$, i.e. $\wh f = f$ on $\dbD_b$, we have
\bea
\label{payfhat}
\ \ \pa_y f(y) = \td_y \wh f(y) - [\td_y \wh f(y) \cd \bn_\dbD(y)] \bn_\dbD(y)= (\mathbf{I}_m - \bn_\dbD(y) \bn_\dbD(y)^\top) \td_y \wh f(y).
\eea
where $\mathbf{I}_m\in\dbR^{m\times m}$ is the identity matrix, and $\mathbf{I}_m - \bn_\dbD(y)\bn_\dbD(y)^\top$ is the projection operator onto $\dbT_\dbD(y)$.
We emphasize that $\pa_y f(y)$ does not depend on the choice of the extension $\wh f$.\footnote{We note that here we require $f$ to be differentiable (as a function on a manifold). Then \eqref{payf} and \eqref{payfhat}  are the usual definition of a derivative of a differentiable function on a manifold, however, 
we can identify them as vectors since $\dbD_b$ is a submanifold of $\dbR^m$. We also note  that tangent space is defined as an equivalence class of curves on the manifold, but we similarly identify it by vectors.} 

We also recall the shape operator $\pa_y \bn_\dbD(y) = [\pa_y \bn^1_\dbD(y),\cds, \pa_y \bn^m_\dbD(y)]\in \dbR^{m\times m}$, which captures the curvatures of $\dbD_b$ at $y$.

\subsection{Set-valued functions}
Consider a function $\dbV: \dbR \to \cD^m_2$. Denote
\bea
\label{Vxy}
\left.\ba{c}
\dis \dbV_b(x):= (\dbV(x))_b,\q \br_\dbV(x,y) := \br_{\dbV(x)}(y),\q \bn_\dbV(x, y) := \bn_{\dbV(x)}(y),\\
\dis \pi_{\dbV}(x, y) := \pi_{\dbV(x)}(y),\q \dbT_\dbV(x, y) := \dbT_{\dbV(x)}(y),
\ea\right.
\eea 
and introduce the graph of $\dbV_b$:
\bea
\label{GV}
\dbG_\dbV := \big\{(x, y): x\in \dbR, y\in \dbV_b(x)\big\}.
\eea
When there is no notational confusion, we may drop the subscript $_\dbV$ in $\br_\dbV$, $\bn_\dbV$, $\pi_{\dbV}$  and denote them as $\br, \bn, \pi$. 
We say $\dbV\in C^2(\dbR; \cD^m_2)$ if $\br_{\dbV}$ is twice continuously differentiable with bounded derivatives on $O_\e(\dbG_\dbV)$ for some $\e>0$.  Note that in this case $\dbV(x)\in \cD^m_2, \forall x\in \dbR$.

\begin{rem}
\label{rem-rextension}
(i) We note that our results in the paper will only involve $\br_{\dbV}$ and its derivatives near $\dbG_\dbV$. For the convenience of our arguments, throughout the paper, we shall modify $\br_\dbV$ outside of $O_\e(\dbG_\dbV)$, so that the modified function $\wh \br_\dbV$ satisfies: 
\begin{itemize}
\item{} $\wh \br_\dbV = \br_\dbV$ on $O_\e(\dbG_\dbV)$;

\item{} $\wh \br_\dbV \in C^2(\dbR\times \dbR^m;\dbR)$ with bounded derivatives;

\item{} $|\wh \br_{\dbV}(x, y)| \ge {\e\over 2}$ for all $(x, y) \notin O_\e(\dbG_\dbV)$.
\end{itemize}
We emphasize that all our results will not rely on the choice of such a modification. For notational simplicity, we may identify the notation $\wh \br_\dbV$ with $\br_\dbV$.

\vskip 3mm

(ii) Similarly we may extend $\pi_{\dbV}$ outside of $O_\e(\dbG_\dbV)$, still denoted as $\pi_{\dbV}$, such that
\begin{itemize}
\item{} On $O_\e(\dbG_\dbV)$, $\pi_{\dbV}(x, y)$ is the original unique projection of $y$ on $\dbV_b(x)$ such that the counterpart of \eqref{piD} holds true;

\item{} $\pi_{\dbV}(x, y)$ is jointly measurable and $\pi_{\dbV}(x, y)\in \dbV_b(x)$ for all $(x, y)\in \dbR\times \dbR^m$;

\item{} There exists a constant $C=C_\dbV$ such that 
\bea
\label{pir}
 |y-\pi_{\dbV}(x, y)|\le C |\br_{\dbV}(x, y)|.  
\eea
\end{itemize}

 We remark that \eqref{pir} follows from \eqref{piD} when $(x, y)\in O_\e(\dbG_\dbV)$, and the existence of $C$ for arbitrary $(x, y)\notin O_\e(\dbG_\dbV)$ is due to the facts that $|\br_\dbV(x, y)|\ge {\e\over 2}$ and that, if needed, we may modify the extension of $\br_\dbV$ further when $|y-\pi_{\dbV}(x, y)|$ is large.

\vskip 2mm

(iii) We can also extend $\bn_\dbV$  to the whole space $\dbR\times \dbR^m$, still denoted as $\bn_\dbV$, such that $\bn_\dbV$ is continuously differentiable with bounded derivatives. For convenience, we fix $\bn_\dbV(x, y) := \td_y \br_\dbV(x, y)$ for all $(x,y)\in\dbR\times \dbR^m$.
\qed
\end{rem}

Fix $x_0\in \dbR$. For each $y_0\in \dbV_b(x_0)$, in light of the modifications in Remark \ref{rem-rextension}, we consider the following ODE which clearly has a unique solution $\U^{x_0,y_0}:\dbR\to \dbR^m$: 
\bea
\label{ODE}
\U^{x_0,y_0}(x) = y_0 - \int_{x_0}^x \td_x \br_\dbV~ \bn_\dbV(\tilde x, \U^{x_0,y_0}(\tilde x)) ~d\tilde x.
\eea
The operator $\U$ is crucial for characterizing the flows on the boundary surface, see Proposition \ref{prop-chain}. We note that the drift  $-\td_x \br_\dbV \bn_\dbV$ here is intrinsic in certain sense, and $\U$ leads to a local geodesic, see Remark \ref{rem-intrinsic-derivative}, Proposition~\ref{prop-geodesic}, and Remark \ref{rem-intrinsic2}  below.

\begin{prop}
\label{prop-chain}
Assume $\dbV\in C^2(\dbR; \cD^m_2)$ and $x_0\in \dbR$. Then, for any $x\in \dbR$,
 \bea
 \label{fundamental}
 \dbV_b(x) = \big\{\U^{x_0,y_0}(x): y_0\in \dbV_b(x_0)\big\}.
 \eea
  Consequently, \eqref{ODE} involves $\br_\dbV$ and $\bn_\dbV$ only on $\dbG_\dbV$ and thus does not depend on the modifications of $\br_\dbV$ and $\bn_\dbV$ in Remark \ref{rem-rextension}.
\end{prop}
\proof For notational simplicity, we drop the subscripts and denote $\br, \bn, \pi$. 

We first show that, for any $y_0\in \dbV_b(x_0)$ and $x>x_0$, $\U(x):= \U^{x_0,y_0}(x)\in \dbV_b(x)$. Let $\e>0$ be such that the original $\br$ in \eqref{rD} is twice continuously differentiable on $O_\e(\dbG_\dbV)$. Note that $(x_0, y_0) \in\dbG_\dbV\subset O_\e(\dbG_\dbV)$. Denote
\beaa
\t:= \inf\big\{x>x_0: (x, \U(x)) \notin  O_\e(\dbG_\dbV)\big\}.
\eeaa
Then, for $x\in [x_0, \t)$, apply the chain rule we have
\beaa
{d\over dx} \br(x, \U(x)) &=& \td_x \br (x, \U(x)) - \td_y \br (x, \U(x)) \cd [\td_x \br~\bn( x, \U(x))]\\
&=& \td_x \br (x, \U(x)) \Big[1- \td_y \br (x, \U(x)) \cd \bn( x, \U(x))\Big].
\eeaa
 Here we used the fact that $\td_x \br$ is scalar. Recall the extension $\bn(x,\U(x)) = \td_y\br(x,\U(x))$ in Remark \ref{rem-rextension} (iii) and by the second equality in \eqref{nD} we have
\begin{equation*}
  {d\over dx} \br(x, \U(x)) = \td_x \br (x, \U(x))\Big[1 - |\td_y \br (x, \U(x))|^2 \Big] = 0,\q x\in [x_0, \t).
\end{equation*}
Note that $\br (x_0, \U(x_0)) = \br (x_0, y_0) =0$. Then the above equation implies $ \br (x, \U(x)) =0$ for all $x\in [x_0, \t)$, which in turn implies $\t=\infty$. Thus  $ \br (x, \U(x)) =0$  and hence $\U(x)\in \dbV_b(x)$ for all  $x\ge x_0$. This implies that $ \{\U^{y_0}(x): y_0\in \dbV_b(x_0)\} \subset \dbV_b(x)$ for all  $x\ge x_0$. Similarly we can show that  $ \{\U^{y_0}(x): y_0\in \dbV_b(x_0)\} \subset \dbV_b(x)$ for all $x\le x_0$, and hence for all $x\in \dbR$.

On the other hand, for any $y\in \dbV_b(x)$, set $y_0 := \U^{x,y}(x_0)\in \dbV_b(x_0)$. By the flow property of the ODE \eqref{ODE}, we have $\U^{x,y}(x') =  \U^{x_0,y_0}(x')$ for all $x'\in \dbR$, and thus $y=  \U^{x,y}(x) = \U^{x_0,y_0}(x)$. This proves the opposite inclusion in \eqref{fundamental}.
\qed

\begin{rem}
\label{rem-fundamental}
(i) Later on we will define $\pa_x\dbV(x, y) :=  -  \td_x \br_\dbV \bn_\dbV(x, y)$ for $(x, y)\in \dbG_\dbV$, see \eqref{patxV} below. Then \eqref{ODE} can be rewritten as
\bea
\label{ODE2}
\U^{x_0,y_0}(x) = y_0 + \int_{x_0}^x \pa_x\dbV(\tilde x, \U^{x_0,y_0}(\tilde x)) d\tilde x.
\eea
Thus \eqref{fundamental} can be viewed as the fundamental theorem of calculus for set-valued functions: 
\bea
\label{fundamental2}
\dbV_b(x) = \dbV_b(x_0) + \int_{x_0}^x \pa_x\dbV(\tilde x, \dbV_b(\tilde x)) d\tilde x.
\eea

(ii) However, \eqref{fundamental2} should be interpreted as \eqref{ODE2} and \eqref{fundamental}, rather than the meaning in the standard set-valued analysis, which roughly speaking considers 
\beaa
\label{setUpsilon}
\tilde \U(x) :=  y_0 +  \int_{x_0}^x \pa_x\dbV(\tilde x, \tilde \g(\tilde x)) d\tilde x,\q\forall y_0\in \dbV_b(x_0), \tilde \g(\tilde x)\in \dbV_b(\tilde x).
\eeaa
The above $\tilde \U(x)$ is in general not in $\dbV_b(x)$. See also related  discussions in Section \ref{sect-AMW} below.

\vspace{1mm}

(iii) The set-valued It\^{o} formula in the next section, which is one of the main results of this paper, can be viewed as the stochastic version of  Proposition \ref{prop-chain}.
\qed
\end{rem}

\begin{rem}
\label{rem-intrinsic-derivative}
For any $\U\in C^1(\dbR; \dbR^m)$ taking values on the boundary surface, namely $\U(x)\in \dbV_b(x)$ for all $x\in \dbR$, by \eqref{nD} we have
\bea
\label{paxVU}
\left.\ba{c}
\dis 0 = {d\over dx}\br_\dbV(x, \U(x)) = \td_x \br_\dbV(x, \U(x)) + \td_y \br_\dbV(x, \U(x)) \cd \td_x \U(x) \\
\dis =  \td_x \br_\dbV(x, \U(x)) + \bn_\dbV(x, \U(x)) \cd \td_x \U(x).
\ea\right.
\eea
That is, for any such $\U$ with $\U(x) = y$, the projection of $\td_x \U(x)$ on $\bn_\dbV(x,y)$ 
is independent of the choice of $\U$: $(\td_x \U \cd \bn_\dbV) \bn_\dbV(x, y) = -  \td_x \br_\dbV \bn_\dbV(x, y) = \pa_x\dbV(x,y)$,  so our definition of $\pa_x\dbV(x,y)$ is intrinsic in the sense that it is independent of the choice of $\U$.
 
 However, if we do not require $\pa_x \dbV$ to go along the normal direction, then there are alternative ways to define $\pa_x \dbV$. This relaxation will not affect our main results though, see Remark \ref{rem-intrinsic2} below.
  \qed
\end{rem}

The next result, although technically not used  in this paper, is interesting in its own right. We postpone its proof to Appendix. 
\begin{prop}
\label{prop-geodesic}
Assume $\dbV\in C^2(\dbR; \cD^m_2)$ and $(x_0, y_0)\in \dbG_\dbV$. Then the curve $\U(x):= \U^{x_0,y_0}(x)$ determined by \eqref{ODE} is a local geodesic of the flow $\dbV$ in the following sense. For any continuous curve $\th(x)\in \dbV_b(x)$ with $\th(x_0)=y_0$, we have
\beaa
\label{geodesic}
\limsup_{x\to x_0} {1\over |x-x_0|}\big[L_\U(x_0, x) - L_\th(x_0, x)\big] \le 0,
\eeaa
where $L_\U(x_0, x)$ (resp. $L_\th(x_0, x)$) denotes the length of $\U$ (resp. $\th$) between $x_0$ and $x$. 
\end{prop}

We now turn to functions $f: \dbG_\dbV \to \dbR$. For fixed $x$, the intrinsic derivative $\pa_y f(x, y)$ for $y\in \dbV_b(x)$ is defined by \eqref{payf} or equivalently by \eqref{payfhat}. We next define the intrinsic derivative of $f$ with respect to $x$ following the local geodesic $\U$ defined by \eqref{ODE}: 
\bea
\label{paxf}
\pa_x f(x_0, y_0) := \lim_{x\to x_0} {f(x, \U^{x_0,y_0}(x))-f(x_0, y_0)\over x-x_0},\q (x_0, y_0)\in \dbG_\dbV.
\eea
Equivalently, for any smooth extension $\wh f$ of $f$, we have
 \bea
\label{paxfhat}
\left.\ba{c}
\dis \pa_x f(x_0, y_0) = \lim_{x\to x_0} {\wh f(x, \U^{x_0,y_0}(x))-\wh f(x_0, y_0)\over x-x_0} \ms\\
\dis = \td_x \wh f(x_0, y_0) - \td_x \br_\dbV \td_y \wh f \cd \bn_\dbV (x_0, y_0).
\ea\right.
\eea
Again, the right side above does not depend on the choice of the extension $\wh f$.

We say $f\in C^1(\dbG_\dbV;\dbR)$ if $f$ has continuous intrinsic derivatives $\pa_y f$ and $\pa_x f$. By \eqref{paxf}, it is obvious that $\pa_x f$ is linear on $f$, and the  product rule and the chain rule remain true:
\bea
\label{product}
\left.\ba{c}
\dis \pa_x (fg) = g\pa_x f  + f\pa_x g,~ \mbox{for all}~ f, g\in C^1(\dbG_\dbV;\dbR);\ms\\
\dis \pa_x [g(f)] = g'(f) \pa_x f, ~\mbox{for all}~  f\in C^1(\dbG_\dbV;\dbR), g\in C^1(\dbR; \dbR).
\ea\right.
\eea

\subsection{Intrinsic derivatives of set-valued functions}

We now extend all the above analyses to functions $\dbV: [0, T]\times \dbR^d \to \cD^m_2$. In this and the next section we may allow infinite time horizon $[0, \infty)$. However, in later sections we require $T$ to be finite, so for simplicity we consider finite $T$ here as well.
Introduce  $\dbV_b(t,x)$, $\br_\dbV(t,x,y)$, $\bn_\dbV(t,x,y)$, $\pi_{\dbV}(t,x,y)$, $\dbT_\dbV(t,x,y)$ in an obvious manner as in \eqref{Vxy} and denote
\bea
\label{GVtx}
\dbG_\dbV := \big\{(t, x, y): (t,x)\in [0, T]\times\dbR^d, y\in \dbV_b(t, x)\big\}.
\eea
As before we may use the simplified notations $\br, \bn, \pi$ when there is no confusion, and we will always use their modified version or extension as in Remark \ref{rem-rextension}.

Recall \eqref{paxf} and \eqref{paxfhat} when $\dbV$ is defined on $\dbR$. Now for our more general $\dbV$ and for any function $f: \dbG_\dbV \to \dbR$, we define its intrinsic partial derivatives $\pa_t f\in \dbR, \pa_x f\in \dbR^d, \pa_y f\in \dbR^m$, and the higher order intrinsic derivatives  in an obvious manner, for example, the second order derivatives are defined as:
\beaa%
\label{2ndorder}
\pa_{x_i x_j} f := \pa_{x_i} (\pa_{x_j} f).  
\eeaa%
Moreover, for $f: \dbG_\dbV \to \dbR^n$, we define its intrinsic derivatives component wise.  

Finally, by considering the special function $f_0(t, x, y) := y$ and its natural extension $\wh f_0(t, x, y) = y$, applying \eqref{paxfhat}  component wise  we define the intrinsic derivatives of $\dbV$.
\begin{defn}
\label{defn-paV}
 For any $(t,x,y)\in \dbG_\dbV$ and by denoting  $f_0(t,x, y) := y$, define
\bea
\label{patxV}
\left.\ba{c}
\dis \pa_t \dbV(t,x, y) := \pa_t f_0(t,x, y) = -\td_t \br ~\! \bn(t,x,y)\in \dbR^m;\ms\\
\dis  \pa_{x_i} \dbV(t,x, y) := \pa_{x_i} f_0(t,x, y) = -\td_{x_i} \br ~\! \bn(t,x,y)\in \dbR^m, \q i=1,\cds, d.
\ea\right.
\eea
\end{defn}
We recall Remark \ref{rem-fundamental} and note that  \eqref{paxfhat} becomes: for any $f\in C^1(\dbG_\dbV;\dbR)$,
 \bea
\label{paxfpaxV}
\pa_t f =  \td_t \wh f +  \td_y \wh f \cd \pa_t\dbV,\q \pa_{x_i} f =  \td_{x_i} \wh f +  \td_y \wh f \cd \pa_{x_i}\dbV,\q \mbox{on}~ \dbG_\dbV.
\eea
Note that $\pa_t\dbV$ and $\pa_x\dbV$ are functions on $\dbG_\dbV$, then we may continue to define higher order derivatives of $\dbV$ by applying \eqref{paxfhat} or \eqref{paxfpaxV} repeatedly. 

\begin{lem}
\label{lem-paxxV}
Assume $\br_\dbV\in C^2([0, T]\times \dbR^d;\cD^m_2)$. Then
\bea
\label{paxxV}
&\dis \pa_{x_ix_j} \dbV(t,x,y) = -  \td_{x_ix_j} \br ~\bn(t,x, y) - \td_{x_j} \br~ \pa_{x_i} \bn(t,x, y);\\
\label{pan}
&\dis  \pa_x \bn^i = \td_{xy_i} \br ;\qq  \pa_y \bn^i = \td_{y_iy} \br .
\eea
\end{lem}
The proof is quite straightforward, we thus postpone it to Appendix.
Throughout the paper, we view
  $\pa_{xx}\dbV\in (\dbR^{d\times d})^m$ as a tensor with components
  $\pa_{xx}\dbV^i\in\dbR^{d\times d}$ for $i=1,\cdots,m$, and we shall
take the notational convention: \bea
\label{pavector}
\left.\ba{c}
\pa_x\dbV := [\pa_{x_1}\dbV, \cds, \pa_{x_d}\dbV]\in \dbR^{m\times d},\q \pa_{xx}\dbV^i := [\pa_{x_1x}\dbV^i, \cds, \pa_{x_dx}\dbV^i]\in \dbR^{d\times d};\ss\\
\pa_x \bn = [\pa_{x_1}\bn,\cds, \pa_{x_d}\bn]\in \dbR^{m\times d},\q \pa_y\bn =[\pa_{y_1}\bn,\cds, \pa_{y_m}\bn]\in \dbR^{m\times m}.
\ea\right.
\eea

\vskip 2mm
\begin{rem}
\label{rem-symmetric}
(i) At $(t,x,y)\in \dbG_\dbV$, since $|\bn|^2=1$, by \eqref{product} we have $\pa_{x_j}\bn \cd \bn=0$. That is, $\pa_{x_j}\bn(t,x,y) \in \dbT_\dbV(t,x,y)$. So \eqref{paxxV} provides an orthogonal decomposition of $\pa_{xx}V$. In particular, unlike the first order derivatives in \eqref{patxV}, $\pa_{x_ix_j} \dbV$ is in general not parallel to $\bn$.

\vspace{2mm}
(ii) It is clear that
\beaa
\pa_{xx} \dbV \cd \bn := \big[\pa_{x_i x_j} \dbV \cd \bn\big]_{1\le i, j\le d} = -  \td_{xx} \br \in \dbR^{d\times d}
\eeaa 
is symmetric. However, in general $\pa_{xx} \dbV$ is  not symmetric: $\td_{x_j} \br \pa_{x_i} \bn \neq \td_{x_i} \br \pa_{x_j} \bn$, $ i\neq j$.

\vspace{2mm}
(iii) $\pa_\by \bn = \td_{yy} \br \in \dbR^{m\times m}$ is symmetric. Moreover, since $\pa_{y_i}\bn \cd \bn =0$, we see that $0$ is an eigenvalue of $\pa_\by \bn$  with eigenvector $\bn$.
\qed
\end{rem}

\begin{eg}
\label{eg-paV}
(i) Let $w: [0, T]\times \dbR^2\to \dbR^2$ and $u: [0, T]\times \dbR^2 \to (0, \infty)$ be continuously differentiable. Set, with $d=m=2$, 
\beaa
&\dis \dbV(t,x) := \big\{y\in \dbR^2: |y-w(t,x)|\le u(t,x)\big\}, \\
&\dis \mbox{and thus}\q \dbV_b(t,x)=\big\{y\in \dbR^2: |y-w(t,x)| = u(t,x)\big\}.
\eeaa
It is clear that
\beaa
 \br (t,x, y) = |y-w(t,x)| - u(t,x),~ (t,x,y) \in [0, T]\times \dbR^2\times \dbR^2.
 \eeaa
 Then, for $(t,x,y)$ near $\dbG_\dbV$ (so that $|y-w(t,x)|>0$),  by straightforward calculation, 
{\small \beaa
 \left.\ba{c}
 \dis \td_t \br = -{(y-w) \cd\td_t w\over |y-w|} - \td_t u;\q \td_{x_i} \br = -{(y-w) \cd\td_{x_i} w\over |y-w|} - \td_{x_i} u;\q \td_{y_i} \br = {y_i-w^i\over |y-w|};\ms\\
 \dis \td_{x_ix_j} \br = {\td_{x_i}w \cd\td_{x_j} w - (y-w)\cd \td_{x_ix_j}w\over |y-w|} - {[(y-w)\cd\td_{x_i} w][(y-w)\cd\td_{x_j}w]\over |y-w|^3}- \td_{x_ix_j} u;\ms\\
 \dis \td_{x_i y_j} \br = - {\td_{x_i}w^j\over |y-w|} + {[(y-w)\cd\td_{x_i}w] (y_j-w^j)\over |y-w|^3},\ \ \td_{y_iy_j} \br = {\1_{\{i=j\}}\over |y-w|} - {[y_i-w^i][y_j-w^j]\over |y-w|^3},
 \ea\right.
 \eeaa}
 Then, by \eqref{nD}, Definition \ref{defn-paV}, and Lemma \ref{lem-paxxV},  at $(t,x,y)\in \dbG_\dbV$ we have
 \bea
 \label{circle}
 \left.\ba{c}
\dis \bn = {y-w\over u};\q \pa_t \dbV = \big[ \td_t w \cd \bn+ \td_t u\big] \bn;\q \pa_{x_i} \dbV = \big[ \td_{x_i} w \cd \bn+ \td_{x_i} u\big] \bn;\ss\\
\dis \pa_{x_i}\bn = {1\over u}\big[- \td_{x_i}w + [\bn\cd\td_{x_i}w] \bn\big],\q \pa_{y_i} \bn^j = {1\over u}\big[\1_{\{i=j\}} - \bn^i\bn^j\big];\ss\\
\dis \pa_{x_ix_j} \dbV = -\Big[{1\over u}\big[\td_{x_i }w\cd\td_{x_j}w - (\td_{x_i} w\cd\bn)(\td_{x_j}w\cd \bn)\big] - \td_{x_ix_j}w \cd \bn - \td_{x_ix_j} u\Big]\bn \\
\dis - {1\over u}\big[\td_{x_j}w\cd \bn + \td_{x_j}u\big]\big[\td_{x_i}w - (\td_{x_i}w\cd\bn)\bn\big].
\ea\right.
\eea
In particular, we see that in general $\pa_{x_ix_j}\dbV \neq \pa_{x_jx_i}\dbV$ for $i\neq j$.

\vspace{2mm}

(ii) Consider a special case that $w=0$ and $u$ satisfies the heat equation:
\beaa%
\label{heat0}
\td_t u + {1\over 2} \tr(\td_{xx} u) =0.
\eeaa%
Then by \eqref{circle} we have $\pa_t \dbV = \td_t u ~\!\bn$, $\pa_{x_ix_j} \dbV = \td_{x_ix_j} u ~\! \bn$, on $\dbG_\dbV$.
 Thus $\dbV$ satisfies the following equation:
$\pa_t \dbV + {1\over 2} \tr(\pa_{xx} \dbV) =0$, on $\dbG_\dbV$.
This clearly implies the following set-valued heat equation:
 \beaa%
\label{heat}
\bn \cd\big[\pa_t \dbV + {1\over 2} \tr(\pa_{xx} \dbV)\big] =0,~\mbox{on}~\dbG_\dbV.
\eeaa

\vspace{-8mm}
\qed
\end{eg}

We remark that we assumed $\br$ had bounded derivatives on $\dbG_\dbV$ in all above analyses. For our applications later, however, $\dbV(T,x)$ could be degenerate, in the sense that $\dbV(T,x) = \{g(x)\}$ is a singleton and hence a degenerate manifold in $\dbR^m$. Note that in \eqref{circle}, $\pa_x\bn$, $\pa_y\bn$, and $\pa_{xx}\dbV$ explode when $u=0$. This motivates us to define the following space.
 
\begin{defn}
\label{defn-C12V}
(i) We say $\dbV \in C^{1,2}([0, T]\times \dbR^d; \cD^m_2)$ if $\br_\dbV \in C^{1,2}(O_\e(\dbG_\dbV); \dbR)$ for some $\e>0$ such that all the related derivatives are bounded and uniformly Lipschitz continuous in $y$.  Consequently, $\pa_t \dbV, \pa_x\dbV, \pa_{xx}\dbV$, $\pa_x\bn$, $\pa_y\bn$ are bounded and uniformly Lipschitz continuous in $y$ on $\dbG_\dbV$.

(ii) We say $\dbV \in C^{1,2}([0, T)\times \dbR^d; \cD^m_2)$ if $\dbV \in C^0([0, T]\times \dbR^d; \cD^m_0)$, and $\dbV \in C^{1,2}([0, T-\d]\times \dbR^d; \cD^m_2)$ for all $0<\d<T$. Note that we do not require $\dbV(T,x)\in \cD^m_2$ here.
\end{defn}

\section{The set-valued It\^{o} formula}
\label{sect-Ito} 
\setcounter{equation}{0}
We start with the probabilistic setting. Let  $\O := \{\o\in C([0, T], \dbR^d): \o_0=0\}$ be the canonical space, $B$ the canonical process, i.e. $B(\o)=\o$,  $\dbP$ the Wiener measure, i.e. $B$ is an $\dbP$-Brownian motion, and $\dbF = \dbF^B$ the augmented filtration generated by $B$.  For a generic Euclidean space $E$ and $p\ge 1$, let $\dbL^p_{loc}(E)$ denote the space of $\dbF$-progressively measurable $E$-valued processes $\th$ with $\int_0^T|\th_t|^pdt <\infty$, a.s., and $\dbL^p_{loc}(\dbR^m; E)$ the space of $\dbF$-progressively measurable random fields $\xi: (t,\o, y)\in [0, T]\times \O\times \dbR^m\to E$ such that $\xi(\cd,\cd, 0)\in \dbL^p_{loc}(E)$ and $\xi$ is uniformly Lipschitz continuous in $y$.

Fix $\dbV\in C^{1,2}([0, T]\times \dbR^d; \cD^m_2)$ as in Definition \ref{defn-C12V} (i) with corresponding $\e>0$, and $x_0\in \dbR^d$, $b\in \dbL^1_{loc}(\dbR^d)$, $\si\in \dbL^2_{loc}(\dbR^{d\times d})$,  $\xi \in \dbL^1_{loc}(\dbR^m; \dbR^{m})$, $\zeta \in \dbL^2_{loc}(\dbR^m; \dbR^{m\times d})$. Denote,
\beaa
\label{X0}
X_t := x_0 + \int_0^t b_s ds + \int_0^t \si_s dB_s,
\eeaa
and introduce the (random) differential operators: recalling \eqref{pavector},
\bea
 \label{kappa}
\left.\ba{c}
\dis \cL^{b,\sigma} \dbV (t, \o, x,y) :=  \big[\pa_t\dbV + \pa_x\dbV b + \tfrac{1}{2}\tr(\sigma^\top\pa_{xx}\dbV\sigma)\big](t, \o, x, y),\\[1.2ex]
\dis  \cK^{\si}_\dbV \zeta(t,\o,x,y) := \big[\tr \big(\zeta^\top\pa_x\bn\sigma + \tfrac{1}{2}\zeta^\top \pa_y\bn\zeta\big) ~\! \bn\big](t, \o, x, y),
\ea\right.
 \eea
 where $\tr(\sigma^\top\pa_{xx}\dbV\sigma)\in \dbR^m$ with $i$-th component $\tr(\sigma^\top\pa_{xx}\dbV^i\sigma)$, and recalling Remark \ref{rem-rextension}, we may extend the derivatives of $\dbV$ and $\bn$ to $[0, T]\times \dbR^d\times \dbR^m$. Here $\pa_x\dbV b(t,\o,x,y) = \pa_x\dbV(t,x,y)b(t,\o)$, and we take this convention for all combinations of functions involving different variables. Moreover, as usual in the literature, when there is no confusion we omit the variable $\omega$. 

 We are now ready to establish the set-valued It\^{o} formula.
\begin{thm}
  \label{thm-Ito} Let $\dbV, x_0, b, \si, X, \xi, \zeta$ be as above. Assume, for each $i$, $\zeta_i(t,\o, y) \in \dbT_{\dbV}(t,X_t(\o), y)$ holds for all $y\in \dbV_b(t,X_t(\o))$, for $dt\times d\dbP$-a.e. $(t,\o)$.   For each $y\in\dbR^m$, recall Remark \ref{rem-rextension} and let $\U^y$ denote the unique strong solution of SDE: 
  \bea
  \label{Upsilon} 
  \U_t^y = y + \int_0^t \big[\cL^{b,\sigma} \dbV - \cK^{\si}_\dbV \zeta+ \xi\big](s,X_s,\U_s^y)ds + \int_0^t \big[\pa_x\mb{V}\sigma + \zeta\big] (s,X_s,\U_s^y)dB_s. \hspace{-1em}
  \eea
  
  (i) Assume $\xi_t(y) \in \dbT_{\dbV}(t,X_t, y)$, for all $y\in \dbV_b(t,X_t)$, for $dt\times d\dbP$-a.e. $(t,\o)$. Then 
  \beaa
  \{\U_t^y: y\in \dbV_b(0,x_0)\} \subset \dbV_b(t,X_t), ~\mbox{for all $0\leq t\leq T$, a.s.}
  \eeaa
   In particular, in this case no extension is needed in \eqref{Upsilon}.
  
  Moreover, if $\dbV_b$ takes values in connected compact sets, then the equality holds: 
  \beaa
  \big\{\U_t^y: y\in \dbV_b(0,x_0)\big\} = \dbV_b(t,X_t),~\mbox{ for all $0\leq t\leq T$, a.s.}
  \eeaa

(ii) Assume $\xi_t(y) \cd \bn_\dbV (t,X_t, y)\le 0$  for all $y\in \dbV_b(t,X_t(\o))$, for $dt\times d\dbP$-a.e. $(t,\o)$.  Then 
\beaa
\{\U_t^y: y\in \dbV_o(0,x_0)\} \subset \dbV_o(t,X_t), ~\mbox{for all $0\leq t\leq T$, a.s.}
\eeaa 

(iii) Assume $\xi_t(y) \cd \bn_\dbV (t,X_t, y)\ge 0$  for all $y\in \dbV_b(t,X_t(\o))$, for $dt\times d\dbP$-a.e. $(t,\o)$.  Then 
\beaa
\{\U_t^y: y\in \dbV^c(0,x_0)\} \subset \dbV^c(t,X_t),~\mbox{for all $0\leq t\leq T$, a.s.}
\eeaa
\end{thm}

\begin{rem}
\label{rem-tangent}
Note that Theorem \ref{thm-Ito} (i) is mainly to characterize the follows on the boundary surface. Since the forces by $\xi, \zeta$ in the tangent space will not push the flow away from the surface, so the It\^{o} formula holds true for arbitrary $\xi, \zeta$, where $\cK^{\si}_\dbV \zeta$ is the correction term due to the quadratic covariation. 

We remark that, while all the results in this theorem remain true if we choose the trivial one: $\xi \equiv 0$, $\zeta\equiv 0$, this special choice is not sufficient for our main results later. For example, in the dynamic programming principle \eqref{DPP} below for the set value function $\dbV$, besides the standard control $\a$, one may choose another free term $\phi\in \dbV(\t,X_\t^{t,x,\a})$ in the right side, which stands for an appropriate selector of $\dbV$. Roughly speaking, the presence of $\phi$ will induce appropriate $\xi, \zeta$. In the standard literature where $\dbV$ is a singleton, $\phi$ is fixed and correspondingly $\xi \equiv 0$, $\zeta \equiv 0$.    

Moreover, since the HJB equation is derived by applying the It\^{o} formula on the dynamic programming equation, our set-valued HJB equation \eqref{HJB} below will also involve $\zeta$.\footnote{$\xi$ is not involved though as we will explain later.}
\qed
\end{rem}

\begin{rem}
\label{rem-measurability}
In this remark we specify the measurability involved in the theorem. Note that the Hausdorff distance $d$ in \eqref{dDD} induces naturally  a topology on $\cD^m_0$ and we may consider the measurability of $\dbV$ with respect to the corresponding Borel field. However, since our analysis goes through $\br_\dbV$, such a topology or measurability is not needed. Indeed, the measurability of the events involved in the above theorem and in the rest of the paper can be derived from $\br_\dbV$ directly. For example, in Theorem \ref{thm-Ito} (i), the event
\beaa
\Big\{ \{\U_t^y: y\in \dbV_b(0,x_0)\} \subset \dbV_b(t,X_t)\Big\} = \Big\{\br_\dbV(t, X_t, \U^y_t) = 0, ~ \forall y\in \dbV_b(0,x_0)\Big\}. 
\eeaa
Since $\br_\dbV$ is continuous on $(t,x,y)$ and $\U^y_t$ is continuous in $y$ in $\dbL^2$-sense, one can easily show that the above event is $\cF_t$-measurable (augmenting $\cF_t$ as usual). 
\qed
\end{rem}

\no{\bf Proof of Theorem \ref{thm-Ito}.} 
  (i) Inspired by the arguments for Proposition
  \ref{prop-chain}, we shall show that 
  \bea
  \label{r=0}
  \bar \br_t := \br(t,X_t,\U_t) = 0.
  \eea
  Fix $y_0\in \dbV_b(0, x_0)$ and denote $\U := \U^{x_0,y_0}$. Introduce 
\bea
\label{Ito-tau}
\t := \inf\big\{t\ge 0: (t, X_t, \U_t) \notin O_\e(\dbG_\dbV)\big\}\wedge T.
\eea
Since $(0, x_0, y_0)\in \dbG_\dbV$, then $\t>0$, and $\br$ is smooth on $[0, \t]$. By the standard It\^o formula, 
\bea
\label{dr} 
&&\dis\q d\bar \br_t = d\br (t,X_t,\U_t) = \L(t,X_t,\U_t)dt  + M(t,X_t,\U_t)dB_t,\q\mbox{where}\\
&&\dis\q \L:=\td_t \br  +  \td_x \br \cd b +   \td_y \br \cd (\cL^{b,\sigma}\dbV- \cK^{\si}_\dbV \zeta + \xi)  \nonumber\\
&&\dis\qq+ \frac{1}{2}\tr\Big(\sigma^\top \td_{xx}\br\sigma + (\pa_x\dbV\sigma + \zeta)^\top \td_{yy}\br (\pa_x\dbV\sigma + \zeta) + 2(\pa_x\dbV\sigma + \zeta)^\top \td_{xy}\br\sigma\Big);\nonumber\\
&&\dis\q M:= \td_x\br^\top \sigma + \td_y\br^\top (\pa_x\dbV\sigma + \zeta).\nonumber
\eea
We proceed the rest of the proof in three steps. 

\ss
{\bf Step 1.} In this step we show that, for $0\le t\le \t$,
  \begin{equation}
      \label{LM}
      \L(t,X_t,\U_t) = \bn(t,X_t,\U_t) \cd \xi_t(\U_t),\q
      M(t,X_t, \U_t) =
      \big[\bn \cd \zeta_1, \cds, \bn \cd \zeta_d\big](t,X_t,\U_t).
  \end{equation}
  Indeed, recall Remark \ref{rem-rextension} (iii) that we have fixed the extension
  $\bn(t,x,y) := \td_y \br(t,x,y)$, and similarly we consider this extension in \eqref{patxV}, \eqref{paxxV},
  and \eqref{pan}.  Then,
  \beaa
  &\!\!\!\! \bn \cd (\cL^{b,\sigma}\dbV- \cK^{\si}_\dbV \zeta) = -\td_t \br -  \td_x \br \cd b - \frac{1}{2}\tr\big(\sigma^\top\td_{xx} \br \sigma + 2\zeta^\top \td_{xy} \br \si + \zeta^\top \td_{yy}\br \zeta\big);\\
  &\big(\td_{yy} \br \pa_x \dbV\big)_{ij} = \td_{y_iy}\br \cd \pa_{x_j} \dbV = -\td_{y_iy}\br \cd \td_{x_j} \br\ \bn= 0,\q 1\le i\le m, 1\le j\le d;\\
  &\Big((\pa_x \dbV)^\top \td_{xy} \br \Big)_{ij}  = \pa_{x_i} \dbV \cd \td_{x_jy}\br   = -\td_{x_i} \br\ \bn \cd \td_{x_jy}\br = 0,\q  1\le i, j\le d.
  \eeaa
  Here we used the facts $\pa_{x_i}\bn \cd \bn=0$, $\td_{y_iy}\br \cd  \bn=0$, $\td_{x_iy}\br \cd \bn=0$. Plug these into the expression of $\L$ in \eqref{dr} we obtain $\L = \bn\cd\xi$ straightforwardly. Similarly,
  \beaa
  M:= \td_x \br^\top \sigma - \bn^\top[ \td_{x_1} \br\ \bn, \cds, \td_{x_d}\br\ \bn]\sigma + \bn^\top \zeta = \bn^\top \zeta.
  \eeaa
  Thus \eqref{LM} holds true.
  
  \ss
  {\bf Step 2.} Recall \eqref{nD} that $\bn(t,x,y) = \bn(t,x,\pi(t,x,y))$ for all $(t,x,y)\in
  O_\e(\dbG_\dbV)$.  Since $\xi(t, y), \zeta_i(t,y)\in \dbT_\dbV(t, X_t, y)$ for any
  $y\in \dbV_b(t,X_t)$, then,   for $1\le i\le d$,
  \beaa
   \bn(t,X_t,\U_t) \cd \xi(t, \pi(t, X_t, \U_t)) = \bn(t,X_t,\U_t) \cd \zeta_i(t, \pi(t, X_t, \U_t))=0, ~ 0\le t\le \t.
   \eeaa
Note further that, by \eqref{nD} again,  $|\U_t - \pi(t,X_t,\U_t)| = |\bar\br_t|$, $0\le t\le \t$. Then
  \begin{equation}
    \begin{aligned}
      \label{LMest}
      \L(t, X_t, \U_t)
      &= \bn(t, X_t, \U_t) \cd \Big(\xi(t,\U_t)- \xi(t, \pi(t, X_t, \U_t))\Big)
      = \tilde b_t \bar\br_t;\\
      M^i(t, X_t, \U_t) &=
       \bn(t, X_t, \U_t) \cd \Big(\zeta_i(t,\U_t)- \zeta_i(t, \pi(t, X_t, \U_t))\Big) =
      \tilde \si^i_t  \bar\br_t;
    \end{aligned}
  \end{equation}
    where $\tilde b_t, \tilde \si^i_t$ are appropriate processes bounded by the Lipschitz
    constants of $\xi, \zeta_i$ with respect to $y$, and $M^i$ denotes the $i$-th column of $M$.  
Thus \eqref{dr} becomes:
 \bea
 \label{drlinear}
 d\bar\br_t  =  \tilde b_t \bar\br_tdt + \tilde \si_t  \bar\br_t dB_t,\q 0\le t\le \t.
 \eea
 Introduce
 \bea
 \label{tildeG}
\dis \tilde \G_t:= \exp\Big(-\int_0^t \tilde \si_s \cd dB_s - \int_0^t [\tilde b_s + {1\over 2}|\tilde \si_s|^2] ds\Big).
 \eea
 Then $\bar \br_t = \bar\br_0 \tilde \G_t^{-1}$.  Since $\bar\br_0  = \br (0,x_0,y_0)=0$, then $\bar\br_t=0$, $0\le t\le \t$, and thus, by \eqref{Ito-tau} we must have 
 $\t=T$, a.s. This implies that $\U_t \in \dbV_b(t, X_t)$, $0\le t\le T$, a.s.
 
 \ss
{\bf Step 3.} Moreover,  assume further that $\dbV_b$ takes values in connected, compact sets. Note that $y\mapsto \U^y_t$ is a homeomorphism almost surely (See Kunita \cite{K1984}). In particular, it is continuous and locally one-to-one. Since $\dbV_b(0,x_0)$ is compact, it is mapped to a closed set in $\dbV_b(0,x_0)$. By invariance of domains for manifolds without boundaries, $y\mapsto \U^y_t$ is an open mapping in relative topologies of $\dbV_b(0,x_0)$,$\dbV_b(t,X_t)$. Therefore, $\dbV_b(0,x_0)$ maps to an open set in $\dbV_b(t,X_t)$. This concludes the equality as we assumed connectedness.

\ss
(ii)  Fix $y_0\in \dbV_o(0,x_0)$ and denote  $\U := \U^{x_0,y_0}$.   It suffices to prove
\bea
\label{r<0}
\bar \br_t := \br(t, X_t, \U_t) <0, \q 0\le t\le T,\q\mbox{a.s.}
\eea
For this purpose, let $\d<{\e\over 2}$ be small enough so that $\bar\br_0 <-\d$. Introduce recursively a sequence of stopping times:  $\t_0:=0$, and for $n=0,1,\cds$,
\beaa 
&\dis \t_{2n+1} := \inf\big\{t>\t_{2n}: \bar\br_t  = -\d\big\}\wedge T;\q
&\dis \t_{2n+2} := \inf\big\{t>\t_{2n+1}: |\bar\br_t| =2\d\big\}\wedge T.
\eeaa
Since $\bar\br_0 <-\d$, by the desired continuity it is clear that 
\beaa
\t_1>\t_0,\q \bar\br_t \le -\d, ~\t_0\le t\le \t_1,\q\mbox{and}\q \bar\br_{\t_1}= -\d~\mbox{when} ~\t_1 <T.
\eeaa

We now assume $\t_1<T$. Then $\t_2>\t_1$, due to the desired continuity, and we consider  $\t_1 \le t\le \t_2$. Since $|\bar\br_{\t_1}|=\d$, we have $|\bar\br_t|\le 2\d<\e$, namely $(t, X_t, \U_t)\in O_\e(\dbG_\dbV)$. Then \eqref{dr} and \eqref{LM} remain true for $\t_1 \le t\le \t_2$, and similarly to \eqref{LMest} and \eqref{drlinear} we have
\beaa
d\bar\br_t  = \big[ \tilde b_t \bar\br_t + \eta_t\big]dt + \tilde \si_t  \bar\br_t  dB_t,~ \t_1 \le t\le \t_2,~\mbox{where}~ \eta_t := \bn(t,X_t,\U_t) \cd \xi(t, \pi(t, X_t, \U_t)).
\eeaa
for some bounded processes $\tilde b, \tilde \si$ constructed as in Step 2. This implies that $\tilde \G_t\bar \br_t =  \bar \br_0 + \int_0^t \tilde \G_s \eta_s ds$ for the $\tilde \G$ in \eqref{tildeG}. Since $\bar\br_{\t_1}= -\d<0$, and by our assumption $\eta_t \le 0$,  then 
\beaa
\bar\br_t <0,~ \t_1 \le t\le \t_2,\q\mbox{and in particular},~ \bar\br_{\t_2}= -2\d~\mbox{when} ~\t_2 < T.
\eeaa

Repeat the arguments, we see that, on $\{\t_n<T\}$,
\bea
\label{taun}
&\dis \t_{n+1}>\t_n,\q\mbox{and}\q \bar\br_t <0,~\mbox{for all}~t \le \t_{n+1};\\
\label{rtaun}
&\dis \bar\br_{\t_k} = -\d,~~ \mbox{for $k\le n$ odd, \q and}\q  \bar\br_{\t_k} = -2\d,~~ \mbox{for $2\le k\le n$ even}.
\eea
Denote $E:= \cap_{n\ge 1} \{\t_n <T\}$. For each $\o\in E$, by \eqref{rtaun} it is clear that $\bar \br$ is not left continuous at $\t_*(\o):= \lim_{n\to\infty} \t_n(\o)$, where the limit exists since $\t_n$ is increasing. However, $\bar \br$ is continuous in $t$, a.s., then we must have $\dbP(E)=0$. That is, for a.e. $\o\in \O$, $\t_n(\o) = T$ for $n$ large enough. Then by \eqref{taun} we see that $\bar \br_t<0$ for $0\le t\le T$, a.s.

  \ss
(iii) follows from similar  arguments as in (ii). 
\qed

\section{A multivariate control problem}
\label{sect-control}
\setcounter{equation}{0}

Recall the canonical setting introduced in the beginning of Section \ref{sect-Ito}.  Given $0\le t<T$, we shall also consider the shifted Brownian motion $B^t_s:= B_s - B_t$, and the shifted filtration $\dbF^t := \dbF^{B^t}$ on $[t, T]$. For a generic Euclidean space $E$, let $\dbL^2(\cF_t, E)$ denote the set of  $\cF_t$-measurable square integrable $E$-valued random variables, and $\dbL^2(\dbF^t, E)$ the set of $\dbF^t$-progressively measurable square integrable $E$-valued processes on $[t, T]$.

 Let $A$ be a domain in some Euclidean space. For each $t\in [0, T]$, our set of admissible controls $\cA_t$ consists of all $\dbF^t$-progressively measurable $A$-valued processes $\a$. We remark that in this paper we consider open loop controls, which is more convenient to study the regularities and to construct desired approximations for our value functions. However, as in standard stochastic control problems, one can easily see that the set values in this section will remain the same if we consider appropriate closed loop controls. 
 
Given $(t, x)\in [0, T]\times \dbR^d$, consider the following controlled dynamics: for each $\a\in \cA_t$,
\bea
\label{FBSDE}
\left.\ba{c}
\dis X^{t,x,\a}_s = x + \int_t^s b(r, X^{t,x,\a}_r, \a_r) dr + \int_t^s \si(r, X^{t,x,\a}_r, \a_r) dB_r,\\
\dis Y^{t,x,\a}_s = g(X^{t,x,\a}_T) + \int_s^T f(r, X^{t,x,\a}_r, Y^{t,x,\a}_r, Z^{t,x,\a}_r, \a_r) dr - \int_s^T Z^{t,x,\a}_r dB_r.
\ea\right.
\eea
Here $X, Y, Z$ take values in $\dbR^d, \dbR^m, \dbR^{m\times d}$, respectively, and $b, \si, f, g$ are in appropriate dimensions and satisfy certain technical conditions which will be specified later. We emphasize that  $Y$ is typically multiple dimensional: $m>1$. Our set value is defined as:
\bea
\label{Vtx}
\dbV(t,x) := \mbox{cl}\big\{Y^{t,x,\a}_t: \a\in \cA_t\big\} \subset \dbR^m.
\eea
Here cl denotes the closure. Thus $\dbV$ is a set-valued mapping $[0, T]\times \dbR^d \to 2^{\dbR^m}$. We now motivate this set value function in the following remarks.

\begin{rem}
\label{rem-1d}
In the scalar case: $m=1$, consider the standard control problems:
\beaa
\label{1dv}
\ul v(t,x) := \inf_{\a\in \cA_t} Y^{t,x,\a}_t,\q  \ol v(t,x) := \sup_{\a\in \cA_t} Y^{t,x,\a}_t.
\eeaa
Then it is obvious that
\bea
\label{1dV}
\dbV(t,x) = [\ul v(t,x), \ol v(t,x)].
\eea
That is, the standard optimization problems are characterizing the boundary of our set value function. In this paper, we will characterize the boundary of $\dbV$ through a set-valued HJB equation, and thus we extend the scalar optimization problem to the multivariate setting. 
\qed
\end{rem}

\begin{rem}
\label{rem-MV}
The set-valued functions can be used to analyze some time inconsistent optimization problems. Consider the well known mean variance optimization problem:
\bea
\label{MVV0}
\left.\ba{c}
 \dis V_0 := \sup_{\a\in \cA} \Big[\dbE[X^{0,x_0,\a}_T] - {\l\over 2} \mbox{Var}(X^{0, x_0,\a}_T)\Big],\\
 \dis\mbox{where}~ X^{t,x,\a}_s = x + \int_t^s \a_r dr + \int_t^s \a_r dB_r.
 \ea\right.
\eea
Here $X, B, \a$ are all scalar processes. Note that $\text{Var}(X_T) = \dbE[ |X_T|^2] - |\dbE[X_T]|^2$. Introduce
\bea
\label{MVdbV}
\left.\ba{c}
\dis \dbV(t,x) := {\rm cl}\big\{Y^{t,x,\a}_t: \a\in \cA_t\big\},\q\mbox{where}\\[1.5ex]
\dis  Y^{t,x,\a, 1}_s = X^{t,x,\a}_T - \int_s^T Z^{t,x,\a, 1}_r dB_r,\q Y^{t,x,\a, 2}_s = |X^{t,x,\a}_T|^2 - \int_s^T Z^{t,x,\a, 2}_r dB_r.
\ea\right.
\eea
Then one can easily verify that 
\bea
\label{MVdbV=V0}
V_0 := \sup_{y\in \dbV(0, x_0)} \f(y),\q\mbox{where} \q \f(y) := y_1 + {\l\over 2} |y_1|^2 - {\l\over 2} y_2.
\eea 
Our goal of this paper is to characterize the dynamic set value function $\dbV$. In fact, in this special case we can solve $\dbV$ explicitly, following the calculation in Pedersen-Peskir \cite[Theorem 3, Part 2]{PP} \footnote{In the  mean variance portfolio selection literature, including \cite{PP}, typically one uses geometric Brownian motion setting and the controlled dynamics (the wealth process) becomes: with $u$ denoting the control,
\bea
\label{GBM}
\tilde X^{t,x,u}_s = x + \int_t^s u_r \tilde X^{t,x,u}_r dr + \int_t^s u_r \tilde X^{t,x,u}_r dB_r.
\eea
Clearly, this is equivalent to our formulation by setting $\a = u \tilde X$ and then $\tilde X^{t,x,u} = X^{t,x,\a}$. Moreover, as already observed in \cite{PP}, the optimal $u^*_s$ explodes when $ \tilde X^{t,x,u}_s=0$. Our $\a^*$ always exists however, as implied by \cite{PP}.  }:
\bea
\label{MVV}
\dbV(t,x) := \Big\{(y_1, y_2): y_1\in \dbR, y_2 \ge e^{-(T-t)} x^2 + \frac{\big(y_1 - xe^{-(T-t)}\big)^2}{1-e^{-(T-t)}}\Big\}\subset \dbR^2.
\eea
Then, given the set $\dbV(0, x_0)$, it is trivial to solve the deterministic optimization problem \eqref{MVdbV=V0}:
\bea
\label{MVoptimal}
\left.\ba{c}
\dis V_0  = x_0 + {1\over 2\l}[e^T-1],\q\mbox{with optimal arguments in \eqref{MVdbV=V0}:}\\
\dis y^\l_1 := x_0+{1\over \l}[e^T-1],\q y^\l_2 :=|y^\l_1|^2 + {1\over \l^2}[e^T-1].
\ea\right.
\eea
 We also refer to Section \ref{sect-moving} below for the related time inconsistency issue.
\qed
\end{rem}

\begin{rem}
\label{rem-target}
(i) The problem \eqref{Vtx} can also be viewed as a stochastic target problem:
\bea
\label{targetV}
\left.\ba{c}
\dis \dbV(t,x) = {\rm cl}\big\{y\in \dbR^m: \exists \a\in \cA_t, Z\in \dbL^2(\dbF^t, \dbR^{m\times d}) ~\mbox{s.t.}~ Y^{t,x, y, \a, Z}_T = g(X^{t,x,\a}_T),~\mbox{a.s.}\big\},\ss\\
\dis \mbox{where}\q Y^{t,x,y, \a, Z}_s = y - \int_t^s f(r, X^{t,x,\a}_r, Y^{t,x,y, \a, Z}_r, Z_r, \a_r) dr + \int_t^s Z_r dB_r.
\ea\right.
\eea

(ii) Note that at above $\dbV(T,x) = \{g(x)\}$ is a singleton. In this respect we may easily extend our setting to non-degenerate terminal $\dbG: \dbR^d \to \cD^m_0$. That is,
\bea
\label{VtxG}
 \dbV(t,x) := {\rm cl} \big\{y\in \dbR^m: \exists \a, Z ~\mbox{such that}~ Y^{t,x, y, \a, Z}_T \in \dbG(X^{t,x,\a}_T),~\mbox{a.s.}\big\}.
  \eea
All our results in this paper can be extended to this case as well, see Section \ref{sect-nondegenerate} below.
\qed
\end{rem}

Throughout the paper, we shall impose the following technical conditions.

\begin{assum}
\label{assum-coefficients}
(i) $(b,\sigma): (t,x,a)\in [0,T]\times \dbR^d\times A \to (\dbR^d,\dbR^{d\times d})$ are bounded, uniformly continuous in $(t,a)$, and uniformly Lipschitz continuous in $x$. 

(ii) $f: (t,x,y,z, a)\in [0,T]\times \dbR^d \times \dbR^m \times \dbR^{m\times d} \times A \to \dbR^m$ is uniformly continuous in $(t, x, a)$ and $f(t,x,0,0,a)$ is bounded. Moreover, $f$  is continuously differentiable in $(y,z)$ with $\td_y f, \td_z f$ bounded and uniformly Lipschitz continuous in  $(y, z)$. 

(iii) $g: x\in \dbR^d\to\dbR^m$ is bounded and uniformly continuous in  $x$.
\end{assum}

It is clear that \eqref{FBSDE} is wellposed for any $\a\in \cA_t$, and thus $\dbV$ is well defined by \eqref{Vtx}.

\subsection{Some basic properties} In this subsection we establish two basic properties of $\dbV$: 
\begin{itemize}
\item For each $(t,x)$, the set value $\dbV(t,x)$ is convex;

\item The set value function $\dbV$ is uniformly continuous in $(t,x)$ under the Hausdorff distance $d$. Consequently, the signed distance function $\br_\dbV$ is  also uniformly continuous in $(t,x)$.
\end{itemize}
These properties are interesting both in theory and in applications. However, we note that, by assuming the smoothness of $\br_\dbV$, our main results in this paper do not rely on these properties. All the proofs in this subsection are postponed to Appendix.

\begin{prop}
\label{prop-convex}
 Under Assumption \ref{assum-coefficients}, the set value $\dbV(t,x)$ is compact and convex.
\end{prop}

\begin{rem}
\label{rem-convex1} The Lipschitz property of $f$ in $(y,z)$ is crucial for the convexity. When $f$ has quadratic growth, as we see in Example \ref{eg-nonconvex} below, this convexity is not guaranteed. We shall remark though such quadratic growth violates Assumption \ref{assum-coefficients}, which is assumed for technical reasons and can be weakened. We refer to Remark \ref{rem-convex} which provides an intuition to this convexity from the PDE perspective.
\qed
\end{rem}

\begin{eg}
\label{eg-nonconvex} Consider the example where $f, g$ and hence $\dbV$ are independent of $x\in\dbR$. Set $m=2$, $A=\{a\in \dbR^2: |a| \le 1\}$, $g = 0$, and $f = (f_1,f_2)^\top$ is specified at below:
\beaa 
 f_1(a, y, z) = a_1,\q f_2(a, y, z)  = \frac{2y_1y_2a_1}{1+y_1^2} + (1+y_1^2)a_2  - \frac{y_2 z_1^2 + 2y_1z_1z_2}{1+y_1^2} + \frac{4y_1^2y_2z_1^2}{(1+y_1^2)^2}.
\eeaa 
Then $\dbV(t)$ can be solved explicitly and is  nonconvex when $T-t>{1\over \sqrt{2}}$:
 \bea
 \label{V-quadratic}
 \left.\ba{c}
 \dis \dbV_b(t) = \Big\{\cY(t,\th): \forall \theta\in [0,2\pi]\Big\},\q\mbox{where}\ms\\
  \dis \cY_1(t,\th):=(T-t)\cos\theta,\q  \cY_2(t,\th):=  (T-t)\big[1 + (T-t)^2\cos^2\theta\big]\sin\theta.
  \ea\right.
  \eea  
  
  \vspace{-6mm}
  \qed
  \end{eg}

\begin{prop}
\label{prop-reg}
Let Assumption \ref{assum-coefficients} hold.

(i) The set value function $\dbV$ is uniformly continuous in $(t,x)$ under the Hausdorff distance $d$ in \eqref{dDD}. That is, for some modulus of continuity function $\rho$,
\bea
\label{Vcont}
 d\big(\dbV(t, x), \dbV(\tilde t, \tilde x)\big) \le \rho(|t-\tilde t|+|x-\tilde x|).
  \eea

(ii) The signed distance function $\br_\dbV$  is also uniformly continuous in $(t,x)$: for the above $\rho$,
\bea
\label{rVcont}
\sup_{y\in \dbR^m}|\br_{\dbV}(t, x, y) -  \br_{\dbV}(\tilde t, \tilde x, y)| \le 3\rho(|t-\tilde t|+|x-\tilde x|).
 \eea
\end{prop}
\no We note that  $\br_\dbV$ is by definition Lipschitz continuous in $y$ with Lipschitz constant $1$.

\begin{rem}
\label{rem-HolderLipschitz}
If we assume further that $f, g$ are uniformly Lipschitz continuous in $x$ and $b,\si, f$ are uniformly H\"{o}lder-${1\over 2}$ continuous in $t$, then following the same arguments one can easily improve the regularity: both $\dbV$ and $\br_{\dbV}$ are uniformly Lipschitz continuous in $x$ and uniformly H\"{o}lder-${1\over 2}$ continuous in $t$.
\end{rem}

\subsection{The dynamic programming principle (DPP)}
We now establish the DPP for $\dbV$.  For $0\le t<T$, $x\in \dbR^d$, $\dbF^t$-stopping time $\t\ge t$, $\phi\in \dbL^2(\cF^t_\t, \dbR^m)$, and $\a\in \cA_t$,  introduce:
\bea
\label{YT0}
\qq Y^{\t, \phi; t,x,\a}_s= \phi + \int_s^\t f(r, X^{t,x,\a}_r, Y^{\t, \phi; t,x,\a}_r, Z^{\t,\phi; t,x,\a}_r, \a_r) dr - \int_s^\t Z^{\t,\phi; t,x,\a}_r dB_r.
\eea

\begin{thm}
  \label{thm:DPP} Let Assumption \ref{assum-coefficients} hold and $\dbV$ be defined by \eqref{Vtx}. For any $0\le t< T$, $x\in \dbR^d$, and any $\dbF^t$-stopping time $\t\ge t$, it holds
  \bea
 \label{DPP} 
\qq \dbV(t,x) = {\rm cl} \big\{Y^{\t, \phi; t,x,\a}_t:  \forall \a\in \cA_t,\ \phi\in \dbL^2(\cF^t_\t,\dbR^m)~\mbox{s.t.}~ \phi\in \dbV(\t,X_\t^{t,x,\a})~\mbox{a.s.} \big\}.
  \eea
\end{thm}
\proof Without loss of generality we prove \eqref{DPP} only at $(t,x) = (0,x_0)$, and for notational simplicity we omit the superscripts $^{0,x_0}$.  Denote the right side of \eqref{DPP} as $\tilde\dbV(0,x_0)$. 

{\bf Step 1.} In this step we show that $\dbV(0,x_0) \subset \tilde \dbV(0,x_0)$. Fix arbitrary $y_0\in \dbV(0,x_0)$ and $\e>0$. By definition of $\dbV(0,x_0)$ there exists $\a = \a^\e\in \cA_0$ such that $|y_0 - Y^\a_0|\le \e$. Denote $\phi := Y^\a_\t$. It is clear that $Y^\a_0 = Y^{\t, \phi; \a}_0$ and thus $|y_0 - Y^{\t, \phi; \a}_0|\le \e$. We claim that
\bea
\label{DPP-claim}
\phi\in \dbV(\t,X_\t^{\a})~\mbox{a.s.} 
\eea
Then $Y^{\t, \phi; \a}_0\in \tilde \dbV(0,x_0)$, and by the arbitrariness of $\e>0$ we obtain $y_0 \in  \tilde \dbV(0,x_0)$.

We first prove \eqref{DPP-claim} in the case that $\t\equiv t$ is deterministic. Consider the shifted canonical space:  $\O_t := \{\o\in C([t, T], \dbR^d): \o_t=0\}$. For any $\o\in \O, \tilde \o\in \O_t$, and $\xi\in \dbL^2(\cF_T)$, introduce 
\beaa
\label{concatenation}
(\o\oplus_t\tilde \o)_s := \o_s \1_{[0, t)} + (\o_t + \tilde \o_s) \1_{[t, T]}(s),\q \xi^{t,\o}(\tilde \o) := \xi(\o\oplus_t\tilde \o).
\eeaa
Then it is clear that $\o\oplus_t\tilde \o\in \O$, and $\xi^{t,\o}\in \dbL^2(\cF^t_T)$ for a.e. $\o\in \O$. In particular, for a.e. $\o\in \O$, we have $\a^{t,\o}\in \cA_t$, and by \eqref{FBSDE} and denoting $\psi^{\a,t,\o} := (\psi^{0, x_0, \a})^{t,\o}$ for $\psi=X, Y, Z$, 
\beaa
\left.\ba{c}
\dis X^{\a, t,\o}_s = X^\a_t(\o) + \int_t^s b(r, X^{\a,t,\o}_r, \a^{t,\o}_r) dr + \int_t^s \si(r, X^{\a,t,\o}_r, \a^{t,\o}_r) dB^t_r,\ms\\
\dis Y^{\a, t, \o}_s = g(X^{\a, t, \o}_T) + \int_s^T  f(r, X^{\a, t, \o}_r, Y^{\a,t, \o}_r, Z^{\a, t, \o}_r, \a^{t, \o}_r) dr - \int_s^T  Z^{\a, t, \o}_r dB^t_r.
\ea\right.
\eeaa
This implies \eqref{DPP-claim} immediately in the case $\t\equiv t$. 

Next, assume $\t$ is discrete, namely it takes finitely many values $t_1<\cds<t_n$. Note that 
\beaa
\phi = Y^\a_\t = \sum_{i=1}^n Y^\a_{t_i}\1_{\{\t=t_i\}},\qq \dbV(\t,X_\t^{\a}) = \sum_{i=1}^n \dbV(t_i,X_{t_i}^{\a}) \1_{\{\t=t_i\}}.
\eeaa
Since $Y^\a_{t_i}\in \dbV(t_i,X_{t_i}^{\a})$, a.s., then it is clear that $\phi\in \dbV(\t,X_\t^{\a})$, a.s.

Finally, for general $\t$, by standard arguments one can construct discrete $\t_n$ such that $\t_n\downarrow \t$. By above $Y^\a_{\t_n}\in \dbV(\t_n,X_{\t_n}^{\a})$, a.s. Note that $Y^\a, X^\a$ are continuous in $t$, a.s., and $\dbV$ is continuous in $(t,x)$, then $Y^\a_{\t_n} \to Y^\a_{\t}$, $\dbV(\t_n,X_{\t_n}^{\a})\to \dbV(\t,X_{\t}^{\a})$, a.s. This, together with the fact that $\dbV(\t,X_{\t}^{\a})$ is closed, implies \eqref{DPP-claim} in the general case.

\ss

{\bf Step 2.} We next prove the opposite inclusion: $\tilde \dbV(0,x_0) \subset \dbV(0,x_0)$.  Fix arbitrary $y_0\in \tilde \dbV(0,x_0)$ and $\e>0$. By definition of $\tilde\dbV(0,x_0)$ there exist $\a = \a^\e\in \cA_0$ and $\phi =\phi^\e \in \dbL^2(\cF_\t, \dbR^m)$ such that  $|y_0 - Y^{\t,\phi;\a}_0|\le \e$ and $\dbP(E) =1$, where
\beaa
E:= \big\{\o\in \O: \phi(\o)\in \dbV(\t(\o), X^\a_\t(\o))\big\}.
\eeaa
Our goal is to construct an $\hat\a\in \cA_0$ such that
\bea
\label{hata-a}
\big|Y^{\hat\a}_0 - Y^{\t,\phi;\a}_0\big| \le C\e.
\eea
Then $|y_0 - Y^{\hat\a}_0|\le \e+C\e$. Since $Y^{\hat\a}_0 \in \dbV(0, x_0)$ by definition, then $y_0\in \dbV(0, x_0)$. 

We construct $\hat\a$ by utilizing the desired regularities as in the standard literature. First let $0=t_0<\cds<t_n=T$ be a partition such that $t_i-t_{i-1}\le \e^2$, $i=1,\cds, n$, and let $\{O^m_j\}_{j\ge 1}$ be a partition of $\dbR^m$ and  $\{O^d_k\}_{k\ge 1}$ a partition of $\dbR^d$ such that the diameter of each  $O^m_j$ and $O^d_k$ is less than $\e$. We now denote
\bea
\label{Eth}
\left.\ba{c}
\dis E^\t_i := \{t_{i-1} < \t \le t_i\},\q E^\phi_j := \{\phi \in O^m_j\},\q E^\a_k := \{X^\a_\t \in O^d_k\},\\
\dis \t_\e:= \sum_{i=1}^n t_i \1_{E^\t_i},\q E_\th := E\cap E^\t_i \cap E^\phi_j \cap E^\a_k,~ \mbox{where}~ \th = (i,j,k).
\ea\right.
\eea
For any $\th=(i, j, k)$ such that $E_\th\neq \emptyset$, choose $\o^\th\in E_\th$ such that
\bea
\label{oth}
\dbP\big(\{\t >t_\th\}\cap E_\th\big) \le \e^2 \dbP(E_\th),\q \mbox{where}\q t_\th := \t(\o^\th),\q x_\th := X^\a_{t_\th}(\o^\th).
\eea
Moreover, since $\phi(\o^\th) \in \dbV(t_\th, x_\th)$, choose $\a^\th\in \cA_{t_\th}$ such that 
\bea
\label{phith}
\big|\phi(\o^\th) - Y^{t_\th, x_\th, \a^\th}_{t_\th}\big| \le \e.
\eea
We then construct $\hat\a\in \cA_0$  by: denoting $\o^t_s := \o_s - \o_t$, $0\le t\le s\le T$, 
\beaa
\label{hata}
\left.\ba{c}
\dis \hat\a_t(\o) := \a_t(\o) \1_{[0, \t_\e(\o))}(t) + \1_{[\t_\e(\o), T]}(t)\Big[\sum_{\th} \1_{E_\th}(\o) \a^\th_t(\o^{t_\th}) + a_0 \1_{E^c}\Big],
\ea\right.
\eeaa
where the summations are over all $\th=(i,j,k)$ with $i=1,\cds, n$ and $j, k\ge 1$, and $a_0\in A$ is an arbitrary value.

\ss
{\bf Step 3.} We now verify \eqref{hata-a}. First, for any $\th=(i, j, k)$ such that $E_\th\neq \emptyset$,  a.s. on $E_\th$ we have $\t_\e = t_i \ge t_\th$ and, denoting $(X^\th, Y^\th, Z^\th) := (X^{t_\th, x_\th, \a^\th}, Y^{t_\th, x_\th, \a^\th}, Z^{t_\th, x_\th, \a^\th})$, 
\beaa
\left.\ba{c}
\dis X^{\hat\a}_t = X^\a_{t_\th} + \int_{t_\th}^t b(s, X^{\hat\a}_s, \a_s) ds + \int_{t_\th}^t \si(s, X^{\hat\a}_s, \a_s) dB_s,\q t\in [t_\th, t_i],\ms\\
\dis X^{\hat\a}_t = X^{\hat\a}_{t_i} + \int_{t_i}^t b(s, X^{\hat\a}_s, \a^\th_s(B^{t_\th})) ds + \int_{t_i}^t \si(s, X^{\hat\a}_s, \a^\th_s(B^{t_\th})) dB_s,\q t\in [t_i, T],\ms\\
\dis X^{\th}_t = x_\th + \int_{t_\th}^t b(s, X^{\th}_s, \a^\th_s(B^{t_\th})) ds + \int_{t_\th}^t \si(s, X^{\th}_s, \a^\th_s(B^{t_\th})) dB_s,\q t\in [t_\th, T].
\ea\right.
\eeaa
Since $b, \si$ are bounded, and $|X^\a_{t_\th} -x_\th|\le \e$, $t_i-t_\th \le \e^2$, by standard SDE estimates we get
\bea
\label{SDEest}
\q \dbE_{\cF_{t_\th}}\Big[\sup_{t_\th \le t\le t_i}|X^{\hat\a}_t - X^\th_t|^2\Big] \le C\e^2,\ \mbox{and then}\  \dbE_{\cF_{t_\th}}\Big[\sup_{t_i \le t\le T}|X^{\hat\a}_t - X^\th_t|^2\Big] \le C\e^2.
\eea
Similarly, note that
\beaa
\left.\ba{c}
\dis Y^{\hat\a}_t = g(X^{\hat\a}_T) + \int_t^T f(s, X^{\hat\a}_s, Y^{\hat\a}_s, Z^{\hat\a}_s, \a^\th_s(B^{t_\th})) ds - \int_t^T Z^{\hat\a}_s dB_s,\q t\in [t_i, T],\ms\\
\dis Y^{\hat\a}_t = Y^{\hat\a}_{t_i} + \int_t^{t_i} f(s, X^{\hat\a}_s, Y^{\hat\a}_s, Z^{\hat\a}_s, \a_s) ds - \int_t^{t_i} Z^{\hat\a}_s dB_s,\q t\in [t_\th, t_i],\ms\\
\dis Y^{\th}_t = g(X^\th_T) + \int_t^T f(s, X^\th_s, Y^\th_s, Z^\th_s, \a^\th_s(B^{t_\th})) ds - \int_t^T Z^\th_s dB_s,\q t\in [t_\th, T].
\ea\right.
\eeaa
Then, by \eqref{SDEest} and standard BSDE estimates we have
\beaa
 \dbE_{\cF_{t_\th}}\Big[\sup_{t_i \le t\le T}|Y^{\hat\a}_t - Y^\th_t|^2\Big] \le C\e^2,\q\mbox{and then}\q   \dbE_{\cF_{t_\th}}\Big[\sup_{t_\th \le t\le t_i}|Y^{\hat\a}_t - Y^\th_t|^2\Big] \le C\e^2.
\eeaa
In particular, this implies that
\bea
\label{BSDEest}
 |Y^{\hat\a}_{t_\th} - Y^\th_{t_\th}| \le C\e,\q \mbox{a.s. on}~ E_\th.
\eea

By Assumption \ref{assum-coefficients}, one can easily see that $Y^\th$, $Y^{\hat\a}$ are bounded. Consider the BSDE:
\bea
\label{tildeYth}
\tilde Y^\th_t = Y^\th_{t_\th}+ \int_t^{t_\th} f(s, x_\th, \tilde Y^\th_s, \tilde Z^\th_s, a_0) ds - \int_t^{t_\th}\tilde Z^\th_s dB_s,\q t\in [t_{i-1}, t_\th].
\eea
Note that $Y^\th_{t_\th}$ is deterministic, then so is $\tilde Y^\th_t$ and thus $\tilde Z^\th_t =0$. Therefore,
\beaa
\sup_{t_{i-1}\le t\le t_\th} |\tilde Y^\th_t - Y^\th_{t_\th}| \le \int_{t_{i-1}}^{t_\th}|f(s, x_\th, \tilde Y^\th_s, 0, a_0)| ds \le C(t_i-t_{i-1}) \le C\e^2.
\eeaa
Moreover, note that $E_\th \in \cF_\t$ and 
\beaa
Y^{\hat\a}_t = Y^{\hat\a}_{t_\th} + \int_t^{t_\th} f(s, X^{\hat\a}_s, Y^{\hat\a}_s, Z^{\hat\a}_s, \a_s) ds - \int_t^{t_\th} Z^{\hat\a}_s dB_s,\q t\in [t_{i-1}, t_\th].
\eeaa
Compare this with \eqref{tildeYth}, by \eqref{BSDEest} and standard BSDE estimate we have
\beaa
\dbE_{\cF_\t}\Big[\sup_{\t\le t\le t_\th} |Y^{\hat\a}_t - \tilde Y^\th_t|^2\Big] \le C\e^2,\q\mbox{a.s. on}~ \{\t\le t_\th\}\cap E_\th.
\eeaa
Then
\beaa
\label{BSDEest1}
|Y^{\hat\a}_\t - Y^\th_{t_\th}| \le |Y^{\hat\a}_\t - \tilde Y^\th_\t| + |\tilde Y^\th_\t - Y^\th_{t_\th}| \le C\e, \q\mbox{a.s. on}~ \{\t\le t_\th\}\cap E_\th.
\eeaa
This, together with \eqref{Eth} and \eqref{phith}, implies that, for a.e. $\o\in \{\t\le t_\th\}\cap E_\th$,
\beaa
|Y^{\hat\a}_\t (\o)- \phi(\o)| \le |Y^{\hat\a}_\t(\o) - Y^\th_{t_\th}|  + |Y^\th_{t_\th} - \phi(\o^\th)| + |\phi(\o^\th) - \phi(\o)|\le C\e.
\eeaa
Note that $\{E_\th\}$ form a partition of $E$ and $\dbP(E)=1$, then by \eqref{oth} we have
\beaa
\dbP\big(|Y^{\hat\a}_\t- \phi| > C\e\big) \le \sum_\th \dbP\big(\{\t> t_\th\}\cap E_\th\big) \le \sum_\th \e^2 \dbP( E_\th) = \e^2.
\eeaa
Note again that $Y^{\hat\a}$ and $\phi$ are bounded. Then
\bea
\label{BSDEest2}
\dbE\big[|Y^{\hat\a}_\t- \phi|^2\big] \le C\e^2 + C\dbP\big(|Y^{\hat\a}_\t- \phi| > C\e\big) \le C\e^2.
\eea 

Finally, note that $Y^{\hat\a}_t = Y^{\t, Y^{\hat\a}_\t; \a}_t$, $0\le t\le \t$. Then, by \eqref{BSDEest2} and standard BSDE estimates we have
\beaa
\dbE\Big[\sup_{0\le t\le \t} |Y^{\hat\a}_t - Y^{\t, \phi; \a}_t|^2\Big] \le C\e^2.
\eeaa
This clearly implies \eqref{hata-a} and hence completes the proof.
\qed

\section{Set-valued HJB equations}
\label{sect-HJB}
\setcounter{equation}{0}

We now derive the set-valued HJB equation for $\dbV$, as usual by combining the DPP \eqref{DPP} and the It\^{o} formula. 
Given $(t,x, y)$, $\a$, and tangent fields $\xi, \zeta$, consider the system starting from $(t,x,y)$. First, by It\^{o}'s formula we have $\U^{\xi,\zeta}_{t+\d}\in \dbV_b(t+\d, X^\a_{t+\d})$; Next, setting $\t:=t+\d$ and $\phi := \U^{\xi,\zeta}_{t+\d}$ in the DPP \eqref{DPP}, we have $Y^{\phi,\a}_t \in \dbV(t,x)$. Here we are using simplified notations. Then it follows from  straightforward calculation that
\beaa
Y^{\phi, \a}_t - y = \dbE\Big[\int_t^T \big[\pa_t \dbV + h_\dbV+ \xi\big] ds\Big],
\eeaa
 for some appropriate $h_\dbV$ which will be defined later. Unlike standard optimization problems, here $Y^{\phi,\a}_t-y$ is a vector, thus one cannot obtain a desired $0$  by optimizing over all $\a$ (and/or all $\phi$). Our remedy relies on the following crucial observation\footnote{\label{boundarymax2}When $\dbV(t,x)$ is convex, we have the stronger inequality: 
$\dis
 \sup_{y' \in \dbV(t,x)} (y'-y)\cd n_\dbV(t,x, y) \le 0$.}: since $y\in \dbV_b(t,x)$, 
\bea
\label{boundarymax}
\lim_{\d\to 0} {1\over \d}\sup_{y' \in \dbV(t,x), |y'-y|\le \d}  \bn_\dbV(t,x, y) \cd (y'-y)\le 0.
\eea
Moreover, by optimizing over all $\a$ and all $\phi$, where the latter effectively means to run over all $\xi$ and $\zeta$, we may obtain $0$ at above. Thus formally we have: noting further that $\bn_\dbV \cd \xi = 0$,
\beaa
0 = \sup_{\a, \xi, \zeta} \lim_{\d\to 0}\ \frac{1}{\d}\ \bn_\dbV \cd  [Y^{\phi, \a}_t - y] = \sup_{\a, \xi, \zeta}\bn_\dbV \cd  \big[\pa_t \dbV + h_\dbV+ \xi\big] = \sup_{\a, \zeta}\bn_\dbV \cd  \big[\pa_t \dbV + h_\dbV\big].
\eeaa
That is, by focusing on the boundary points, we are actually considering scalarized optimization problems along the direction $\bn_\dbV$.

We are now ready to present the set-valued HJB equation precisely. For $(t,x,y) \in \dbG_\dbV$, $z\in \dbR^{m\times d}$, $\g\in (\dbR^{d\times d})^m$, $a\in A$, $\zeta\in (\dbT_\dbV(t,x,y))^d$, inspired by \eqref{Upsilon}  in the It\^{o} formula and \eqref{DPP} we introduce the Hamiltonian function $H_\dbV$:
\bea
\label{H}
\left.\ba{lll}
 \dis \q\cK_\dbV(t,x,y,a, \zeta):=\tr \Big(\zeta^\top\pa_x\bn_\dbV(t,x,y)\sigma(t,x,a) + \frac{1}{2}\zeta^\top \pa_y\bn_\dbV(t,x,y)\zeta\Big)\bn_\dbV(t,x,y);\ss\\
  \dis\q h^0_\dbV(t, x,y, z, \g, a,\zeta) :=  z b(t,x,a)  + \frac{1}{2}\tr(\sigma^\top\g \sigma(t,x,a))- \cK_\dbV(t,x,y,a, \zeta); \ss\\
  \dis\q  h_\dbV(t, x,y, z, \g, a,\zeta) := h^0_\dbV(t, x,y, z, \g, a,\zeta) + f\big(t,x, y, z\si(t,x,a)+\zeta, a\big);\ss\\
  \dis\q H_\dbV(t,x,y,z,\g) := \sup_{ a\in A, \zeta\in (\dbT_\dbV(t,x,y))^d} \bn_\dbV(t,x,y)\cd h_\dbV(t, x,y, z, \g, a, \zeta);
\ea\right.
\eea
We note that, in the scalar case: $m=1$, then $\bn = 1$ or $-1$, and $H_\dbV$ reduces to standard Hamiltonian functions. See Remark \ref{rem-HJB} (i) below. 

Our set-valued HJB equation takes the form:
\bea
\label{HJB}
\left.\ba{c}
\cL \dbV(t,x,y)  =0,~\forall (t,x, y)\in \dbG_\dbV,\q \mbox{where}\ss\\
\cL \dbV(t,x,y) :=\pa_t \dbV(t,x,y) \cd \bn_\dbV(t,x,y)+ H_\dbV\big(t,x, y, \pa_x\dbV(t,x,y), \pa_{xx}\dbV(t,x,y)\big).
\ea\right.
\eea
Equivalently, by \eqref{nD}, \eqref{patxV}, \eqref{paxxV}, \eqref{pan}, we may rewrite the above equation:
\bea
\label{HJBr}
\left.\ba{c}
\dis \td_t \br_\dbV + \inf_{a\in A, \zeta\in (\dbT_\dbV(t,x,y))^d}\Big[ \td_x \br_\dbV \cd b + {1\over 2} \tr\big(\si^\top \td_{xx} \br_\dbV \si + 2 \zeta^\top \td_{xy}\br_\dbV \si + \zeta^\top \td_{yy}\br_\dbV \zeta\big) \ss\\
\dis - \td_y \br_\dbV \cd f(t ,x, y, - \td_y\br_\dbV(\td_x \br_\dbV)^\top\si+ \zeta, a)\Big]=0,\q (t,x,y)\in \dbG_\dbV.
\ea\right.
\eea

\begin{rem}
\label{rem-HJB}
 (i)  For the scalar case as in Remark \ref{rem-1d}, by \eqref{1dV} we have 
 \beaa
 &\dbV_b(t,x) = \{\ul v(t,x), \ol v(t,x)\}, \q \bn(t,x, \ul v(t,x))=-1,\q  \bn(t,x, \ol v(t,x))=1,\\
 & \dbT_\dbV(t,x, \ul v(t,x)) =\dbT_\dbV(t,x, \ol v(t,x))=\{0\}.
 \eeaa
  In the neighborhood of  $y=\ul v(t,x)$, we have $\br_\dbV(t,x, y) = \ul v(t,x) - y$ and $\bn(t,x, \ul v(t,x))=-1$. Then  \eqref{HJBr} reduces to the standard HJB equation for $\ul v$:
  \beaa
  \td_t \ul v + \inf_{a\in A}\Big[ \td_x \ul v \cd b + {1\over 2} \tr\big(\si^\top \td_{xx} \ul v \si\big) +  f\big(t ,x, \ul v, (\td_x\ul v)^\top\si, a\big)\Big]=0.
  \eeaa
  Similarly, in the neighborhood of  $y=\ol v(t,x)$, we have $\br_\dbV(t,x, y) = y-\ol v(t,x)$ and $\bn(t,x, \ol v(t,x))=1$. Then  \eqref{HJBr} reduces to the standard HJB equation for $\ol v$:
  \beaa
  \td_t \ol v + \sup_{a\in A}\Big[ \td_x \ol v \cd b + {1\over 2} \tr\big(\si^\top \td_{xx} \ol v \si\big) +  f\big(t ,x, \ol v, (\td_x\ol v)^\top\si, a\big)\Big]=0.
  \eeaa  
  
  (ii)  Although $\br_\dbV$ is scalar, we emphasize that \eqref{HJBr} holds true only on $\dbG_\dbV$, and the set $\dbG_\dbV$ is in turn determined by the solution $\br_\dbV$. So \eqref{HJBr} is  actually quite involved, and we can not apply the standard PDE theory on it.
  
  (iii) It is clear that $\dbV(T,x)=\{g(x)\}$ is degenerate, so we do not require the smoothness of $\dbV$ at $T$. See Definition \ref{defn-C12V} and the paragraph above it.
  \qed
\end{rem}

\begin{rem}
  \label{rem-convex} Note that $\cK_\dbV$ relies on $\zeta$ quadratically, and the space of $\zeta$ is typically unbounded.  This is associated with the fact that the PDE \eqref{WHJB} in the proof of Proposition \ref{prop-convex} is quadratic in $z$, and is related to the convexity of $\dbV$ in Proposition \ref{prop-convex}.
   
 (i) In the scalar case: $m=1$, we have $\dbT_\dbV(t,x, y) = \{0\}$ for all $(t,x,y)\in \dbG_\dbV$. Then this issue is trivial.
   Indeed, in this case the set-valued PDE reduces back to the standard HJB equations, as we saw in Remark \ref{rem-HJB} (i).
 
 (ii) For $m\ge 2$, recall Remark \ref{rem-symmetric} (iii) that $\pa_y \bn_\dbV = \td_{yy} \br_\dbV$ is symmetric and $0$ is an eigenvalue with eigenvector $\bn$. At any fixed $(t,x,y)\in \dbG_\dbV$, let $\l_1\le \cds\le \l_{m-1}$ be the other eigenvalues. When $f$ has linear growth in $z$, which is implied by the Lipschitz continuity, and $\l_1>0$, then clearly $H_\dbV<\infty$. In this case $\dbV(t,x)$ is strictly convex. 
 
 (iii) When $f$ has linear growth in $z$, one may easily derive from
 $H_\dbV<\infty$ that $\l_1 \ge 0$. Thus, our classical solution
 $\dbV$ has to be convex.

(iv) When $f$ has quadratic growth, this convexity is not required, as we saw in Example \ref{eg-nonconvex}. We shall remark though such quadratic growth violates Assumption \ref{assum-coefficients}, which is assumed for technical reasons and can be weakened.   
\qed
\end{rem}

\begin{rem}
\label{rem-intrinsic2}
Consider Remarks \ref{rem-fundamental} and \ref{rem-intrinsic-derivative} again with $x\in \dbR$. One can show that, roughly speaking, $\dbV$ satisfies the fundamental theorem \eqref{fundamental2} if and only if $\pa_x\dbV :=  -  \td_x \br_\dbV \bn_\dbV + \zeta_0$ for some appropriate field $\zeta_0$ taking values in the tangent space $\dbT_\dbV$. So, by restricting $\pa_x \dbV(x,y)$ to be parallel to $\bn_\dbV(t,x,y)$, we see that $\pa_x\dbV :=  -  \td_x \br_\dbV \bn_\dbV$ is the unique choice. If we relax this requirement, then we may define $\pa_x \dbV$ alternatively by choosing certain $\zeta_0$, denoted as $\pa_x^{\zeta_0}\dbV$. This choice of $\zeta_0$ will not affect our main results though. Indeed, one can show that $\bn^\top\pa_{xx}^{\zeta_0}\dbV = \bn^\top\pa_{xx}\dbV - \zeta_0^\top \pa_x\bn - \pa_x\bn^\top\zeta_0 - \zeta_0^\top \pa_y\bn \zeta_0$. Thus, in the set-valued It\^{o} formula, the presence of $\zeta_0$ will only lead to a different choice of $(\xi,\z)$ in \eqref{Upsilon}. Then, by taking the supremum over $\zeta$, the Hamiltonian $H_\dbV$ in \eqref{H} will remain the same. Consequently,  the set-valued PDE \eqref{HJB} also does not depend on $\zeta_0$.

Finally, due to the local geodesic property in Proposition \ref{prop-geodesic}, throughout the paper we fix $\zeta_0\equiv 0$ so that $\pa_x \dbV$ is align with $\bn_\dbV$.
  \qed
\end{rem}

We now turn to the wellposedness of \eqref{HJB}.  We first define  classical solutions rigorously. 
\begin{defn}
\label{defn-classical}
(i) Let $C^{1,2}_0([0, T)\times \dbR^d;\cD^m_2)$ denote the set of $\dbV \in C^{1,2}([0, T)\times \dbR^d;\cD^m_2)$ such that, for any $T_0<T$, the eigenvalue $\l_1$ of $\pa_y \bn (t,x,y)$ in Remark \ref{rem-convex} (ii) has a lower bound  $c_{T_0}>0$ for all $(t,x)\in [0, T_0]\times \dbR^d$, $y\in \dbV_b(t,x)$. That is, 
\bea
\label{HVbound}
\tr \Big(\zeta^\top \pa_y\bn_\dbV(t,x,y)\zeta\Big) \ge c_{T_0} |\zeta|^2\q \forall (t,x)\in [0, T_0]\times \dbR^d, y\in \dbV_b(t,x), \zeta\in \dbT_\dbV(t,x,y).
\hspace{-2em}
\eea
This implies that $H_\dbV(\cd, \pa_x\dbV, \pa_{xx}\dbV)$ is finite and uniformly continuous whenever $t\le T_0$.

(ii) We say $\dbV \in C^{1,2}_0([0, T)\times \dbR^d;\cD^m_2)$ is a classical solution to \eqref{HJB} if 
it satisfies \eqref{HJB} for all $(t,x)\in [0, T)\times \dbR^d$ and $y\in \dbV_b(t,x)$.
\end{defn}

We shall provide an example in Example \ref{eg-classical} below. We next establish a crucial estimate, whose proof is postponed to Appendix.
\begin{lem}
\label{lem-boundary}
Let Assumption \ref{assum-coefficients} hold and $\dbV$ be defined by \eqref{Vtx}. Assume $\dbV \in C^{1,2}_0([0, T)\times \dbR^d;\cD^m_2)$. Fix $T_0<T$ and $x_0\in \dbR^d$. Let $\e, \d>0$ and $\a \in \cA_0$ be such that $|\br_\dbV(0, x_0, Y^\a_0)|\le \e$, where $(X^\a, Y^\a, Z^\a)=(X^{0,x_0, \a}, Y^{0,x_0, \a}, Z^{0,x_0, \a})$ are defined by \eqref{FBSDE}. Then there exists a constant $C_{T_0}$, which may depend on $T_0$ but not on $\e, \d, \a$, such that
\bea
\label{boundaryest}
\dbP\Big(\sup_{0\le t\le T_0} |\br_\dbV(t, X^\a_t, Y^\a_t)|  > \d\Big) \le C_{T_0}\sqrt{\e\over \d}.
\eea
In particular, if $Y^\a_0 \in \dbV_b(0, x_0)$, then $Y^\a_t \in \dbV_b(t, X^\a_t)$, $0\le t\le T$, a.s.
\end{lem}

The main result of this section is the following theorem.
\begin{thm}
\label{thm-existence}
Let Assumption \ref{assum-coefficients} hold and $\dbV$ be defined by \eqref{Vtx}. Assume $\dbV \in C^{1,2}_0([0, T)\times \dbR^d;\cD^m_2)$. Then $\dbV$ is a classical solution of \eqref{HJB} with terminal condition $\dbV(T,x) = \{g(x)\}$.
\end{thm}

\proof It is clear that $\dbV(T,x) = \{g(x)\}$. Without loss of generality, we shall verify \eqref{HJB} only at a fixed $(0,x_0, y_0)\in \dbG_\dbV$, and for notational simplicity, in this proof we omit the supscripts $^{0,x_0, y_0}$ and the subscript $_\dbV$ in $\br, \bn, \pi$.  We proceed in two steps, following the ideas in the standard stochastic control literature, but by adopting the scalarization with $\bn$ as explained in the beginning of this section.

{\bf Step 1.} We first show that $\cL\dbV(0,x_0, y_0) \le 0$. For this purpose, we fix an arbitrary $a\in A$ and let $X^a:= X^{0,x_0,\a}$ be defined by \eqref{FBSDE} for constant control process $\a\equiv a$. Moreover, we fix arbitrary $\xi,\zeta_i:[0, T]\times\Omega\times\dbR^m\to\dbR^m$, $i=1,\cds, d$, which are $\dbF$-progressively measurable,  bounded,  continuous in $t$,  uniformly Lipschitz continuous in $y$, and $\xi(t,\o, y), \zeta_i(t,\o,y) \in \dbT_\dbV(t, X^a_t(\o), y)$ for all $y\in \dbV_b(t, X^a_t(\o))$, for $dt\times d\dbP$-a.e. $(t,\o)$. Denote $\zeta=(\zeta_1,\cds, \zeta_d)$ and consider the SDE:
\bea
\label{Upsilona}
\left.\ba{c}
\dis \U^{a,\xi,\zeta}_t = y_0 + \int_0^t \big[\pa_t \dbV+ h^0_\dbV(\cd, \pa_x \dbV, \pa_{xx}\dbV, a, \zeta_s) + \xi\big](s, X^a_s, \U^{a,\xi,\zeta}_s) ds \\
\dis + \int_0^t [\pa_x\dbV \si(\cd,a) +\zeta\big](s, X^a_s, \U^{a,\xi,\zeta}_s) dB_s.
\ea\right.
\eea
Applying the It\^{o} formula Theorem \ref{thm-Ito} we have $\U^{a,\xi,\zeta}_t \in \dbV_b(t, X^a_t)$, for all $0\le t<T$.

Now for any $\d>0$ small, consider the BSDE \eqref{YT0} with terminal condition $(\d, \U^{a,\xi,\zeta}_\d)$:
\bea
\label{Ya}
Y^{\d,a,\xi,\zeta}_t= \U^{a,\xi,\zeta}_\d + \int_t^\d f(s, X^{a}_s, Y^{\d,a,\xi,\zeta}_s, Z^{\d,a,\xi,\zeta}_s, a) ds - \int_t^\d Z^{\d,a,\xi,\zeta}_s dB_s, 0\le t\le \d. \hspace{-1em}
\eea
Since $\U^{a,\xi,\zeta}_\d \in \dbV_b(\d, X^a_\d)$, by the DPP \eqref{DPP} we see that $Y^{\d,a,\xi,\zeta}_0\in \dbV(0,x_0)$. Denote 
\beaa
\D Y^\d_t:= Y^{\d,a,\xi,\zeta}_t-\U^{a,\xi,\zeta}_t,\q \D Z^\d_t:= Z^{\d,a,\xi,\zeta}_t- [\pa_x\dbV \si(\cd,a) +\zeta\big](t, X^a_t, \U^{a,\xi,\zeta}_t).
\eeaa
Then, by \eqref{Upsilona} and \eqref{Ya} we have
\beaa
&&\dis\!\!\!\!\!\! \D Y^\d_t = \int_t^\d \Big[\pa_t \dbV+h_\dbV(\cd, \pa_x\dbV, \pa_{xx}\dbV, a, \zeta) + \xi\Big](s, X^a_s, \U^{a,\xi,\zeta}_s) ds - \int_t^\d \D Z^\d_s dB_s\\
&&\dis + \int_t^\d\Big[f(s, X^{a}_s, Y^{\d, a,\xi,\zeta}_s, Z^{\d,a,\xi,\zeta}_s, a) - f(s, X^{a}_s, Y^{\d,a,\xi,\zeta}_s-\D Y^\d_s,Z^{\d,a,\xi,\zeta}_s-\D Z^\d_s, a)\Big]ds\\
&&\!\!\!\!\!\! =\int_t^\d \Big[\pa_t \dbV+h_\dbV(\cd, \pa_x\dbV, \pa_{xx}\dbV, a, \zeta) + \xi\Big](s, X^a_s, \U^{a,\xi,\zeta}_s) ds - \int_t^\d \D Z^\d_s dB_s\\
&&\dis + \int_t^\d\Big[\tilde b_s \D Y^\d_s + \tilde\si_s \D Z^\d_s \Big]ds,
\eeaa
where $\tilde b, \tilde \si$ are appropriate $\dbF$-progressively measurable bounded processes. Then, for the $\tilde \G$ defined by \eqref{tildeG}, we have
\beaa
\tilde \G_t\D Y^\d_t &=&\int_t^\d \tilde \G_s\Big[\pa_t\dbV+h_\dbV(\cd, \pa_x\dbV, \pa_{xx}\dbV, a, \zeta) + \xi\Big](s, X^a_s, \U^{a,\xi,\zeta}_s) ds\\
&& - \int_t^\d \tilde \G_s\big[\D Z^\d_s +  \D Y^\d_t \tilde\si \big] \cd dB_s.
\eeaa
In particular,
\bea
\label{DYd0}
\D Y^\d_0 = \dbE\Big[\int_0^\d \tilde \G_s\big[\pa_t\dbV+h_\dbV(\cd, \pa_x\dbV, \pa_{xx}\dbV, a, \zeta) + \xi\big](s, X^a_s, \U^{a,\xi,\zeta}_s) ds\Big].
\eea
 
Given our conditions, it is clear that $|\D Y^\d_0|\leq C\delta$, which implies that $\dis\lim_{\d\to 0}Y^{\d,\a,\xi,\zeta}_0=y_0$. Since $Y^{\d,a,\xi,\zeta}_0 \in \dbV(0,x_0)$ and $y_0\in \dbV_b(0,x_0)$, then it follows from \eqref{boundarymax}\footnote{ We emphasize that here we are using \eqref{boundarymax}, rather than the stronger inequality in Footnote \ref{boundarymax2} for convex sets. That is, although our set values turn out to be convex, the arguments here do not rely on this convexity.}   that:
\bea
\label{DYd0<0}
\limsup_{\d\to 0} {1\over \d} \big[\bn(0,x_0, y_0) \cd \D Y^\d_0\big] \le 0.
\eea
 Now by \eqref{DYd0} and the desired continuity of $\tilde \G_s, X^a_s, \U^{a,\xi,\zeta}_s, \xi_s$ in $s$ as well as the desired regularity of all the involved functions in $(x, y)$, we have
 \beaa
 0&\ge& \limsup_{\d\to 0} {1\over \d}  \dbE\Big[\bn(0,x_0, y_0) \cd\int_0^\d \tilde \G_s\Big[\pa_t\dbV+h_\dbV(\cd, \pa_x\dbV, \pa_{xx}\dbV, a, \zeta) + \xi\Big](s, X^a_s, \U^{a,\xi,\zeta}_s) ds\Big]\\
 &=&\limsup_{\d\to 0} {1\over \d}  \Big[\bn(0,x_0, y_0) \cd\int_0^\d \Big[\pa_t\dbV+h_\dbV(\cd, \pa_x\dbV, \pa_{xx}\dbV, a, \zeta) + \xi\Big](0, x_0, y_0) ds\Big]\\
 &=&\bn(0,x_0, y_0) \cd \Big[\pa_t\dbV+h_\dbV(\cd, \pa_x\dbV, \pa_{xx}\dbV, a, \zeta)\Big](0, x_0, y_0),
 \eeaa
 where the last equality is due to the assumption $\xi\in \dbT_\dbV$. Now by the arbitrariness of $a, \zeta$ we obtain $\cL\dbV(0,x_0, y_0) \le 0$.

{\bf Step 2.} We next show that  $\cL\dbV(0,x_0, y_0) \ge 0$.  This step is more involved, because one needs to construct approximately "optimal" $\a, \zeta, \xi$ to make the inequality \eqref{boundarymax} essentially reaching $0$.  For this purpose, fix $T_0<T$, and throughout this proof, the generic constant $C$ may depend on $T_0$. Since $\dbV\in C^{1,2}([0, T_0]\times \dbR^d; \cD^m_2)$, there exists $\e_0>0$ such that $\br\in C^{1,2}(O^{T_0}_{\e_0}(\dbG_\dbV);\dbR)$, where $O^{T_0}_{\e_0}(\dbG_\dbV):=\{(t,x,y)\in [0, T_0]\times \dbR^d\times \dbR^m: |\br (t,x,y)|\le \e_0\}$. Fix a sufficiently small constant $\e>0$.  Since $y_0\in \dbV_b(0,x_0)$, there exists $\a = \a^\e\in \cA_0$ such that
  \beaa
\pi(0, x_0, Y^\a_0) = y_0\q\mbox{and}\q |y_0 - Y^\a_0|\le \e^4,
 \eeaa
 where $(X^\a, Y^\a, Z^\a) = (X^{0,x_0, \a}, Y^{0,x_0, \a}, Z^{0,x_0,\a})$ are defined by \eqref{FBSDE}.   Define\footnote{Since we are using the augmented filtration, the $\t$ here is an $\dbF$-stopping time. }
 \beaa
 \label{tau}
 \t:=\t_{\e, \a} :=\inf\{t>0: |\br (t, X^\a_t, Y^\a_t)|  > \e^2\} \wedge T_0.
 \eeaa
 By Lemma \ref{lem-boundary} we have
 \bea
 \label{tauest}
 \dbP(\t < T_0) \le \dbP\big(\sup_{0\le t\le T_0}  |\br (t, X^\a_t, Y^\a_t)|  >\e^2\big) \le  C\sqrt{\e^4\over \e^2} = C\e.
 \eea
 
 {\bf Step 2.1.} Introduce two random fields:
 \bea
 \label{zeta}
 \left.\ba{lll}
\dis \zeta_t(y) := Z^\a_t - \pa_x\dbV\sigma(t,X_t^\a, y, \a_t), \ms\\
\dis \xi_t(y) := -\big[\pa_t\dbV+h_\dbV(\cd, \pa_x\dbV, \pa_{xx}\dbV, \a_t, \zeta_t)\big](t,X_t^\a,y). 
\ea\right.
  \eea
  Then we may rewrite the BSDE for $(Y^\a, Z^\a)$ forwardly: 
\beaa 
Y_t^{\a} &=& Y_0^{\a} + \int_0^t \Big[\pa_t\dbV+h^0_\dbV(\cd, \pa_x\dbV, \pa_{xx}\dbV, \a_s, \zeta)+ \xi\Big](s,X_s^\a,Y_s^\a)ds \\
&&+ \int_0^t \big[\pa_x\dbV\sigma(\cd,\a_s)+ \zeta\big](s,X_s^\a,Y_s^\a)dB_s. 
\eeaa 
We note that, if $Y^\a_t \in \dbV_b(t, X^\a_t)$ and $\zeta_t(Y^\a_t), \xi_t(Y^\a_t)$ are in the tangent space $\dbT_\dbV(t, X^\a_t, Y^\a_t)$, then by the optimality of $H_\dbV$ we have $\cL\dbV(t, X^\a_t, Y^\a_t) \ge -\bn(t, X^\a_t, Y^\a_t) \cd \xi_t(Y^\a_t)=0$, which is the desired inequality. In this and the next substeps, we shall prove these properties in approximate sense.  
  
 Denote $\pi^\a_t := \pi(t, X^\a_t, Y^\a_t)$.
By \eqref{dr} and \eqref{LM}, similarly to \eqref{LMest} and \eqref{drlinear} we have, 
\bea
\label{dr2}
\left.\ba{c}
\dis d \br (t, X^\a_t, Y^\a_t) = \Big[\br (t,X_t^\a,Y^\a_t) \tilde b_t -\bn\cd \xi (t, X^\a_t, \pi^\a_t) \Big]dt \ms\\
\dis\qq\qq\qq + \Big[\br (t,X_t^\a,Y^\a_t) \tilde \si_t - \bn^\top\zeta (t, X^\a_t, \pi^\a_t) \Big]dB_t, \q0\le t\le \t,
\ea\right.
\eea
where $\tilde b, \tilde \si$ are $\dbF$-progressively measurable and  satisfy: for some constant $C=C_{T_0}$,
\bea
\label{tildebbound}
|\tilde \si_t|\le C,\q |\tilde b_t|\le C\big[1+|Z_t^\a|\big].
\eea
Recall the process $\tilde\G$ defined in \eqref{tildeG}, we have
\bea
\label{rtau}
\br (0, x_0, Y^\a_0) - \tilde \G_\t \br (\t, X^\a_\t, Y^\a_\t) =\! \int_0^\t\!\! \tilde \G_s \bn\cd \xi (s, X^\a_s, \pi^\a_s) ds + \! \int_0^\t \!\! \tilde \G_s \bn^\top \zeta(s, X^\a_s, \pi^\a_s) dB_s. \hspace{0em}
\eea
Moreover, by Assumption \ref{assum-coefficients} one can easily see that $Y^\a$ is bounded, then $\int_0^{\t\wedge \cd} Z^\a_s dB_s$ is a BMO martingale, and thus there exist $c_0, C>0$,  such that (cf. \cite[Section 7.2]{Zhang})
\beaa
\label{BMOest}
\dbE\Big[\exp\big(c_0 \int_0^\t |Z^\a_s|^2 ds\big)\Big] \le C < \infty.
\eeaa
In particular, this implies that, for any $p\ge 1$, there exists a constant $C_p>0$ such that
\bea
\label{tildeGest}
\dbE\Big[\sup_{0\le t\le \t} \big[|\tilde \G_t|^p+|\tilde \G_t|^{-p}\big]\Big] \le C_{p}.
\eea
Applying the standard It\^o formula on $|\tilde \G_t \br(t,X_t^\a,Y_t^\a)|^2$, by \eqref{dr2} we have
\bea
\label{zetaest}
&&\qq \dbE\Big[\int_0^\tau |\tilde\G_s\bn^\top\zeta(s,X_s^\a,\pi_s^\a)|^2ds\Big] \nonumber\\
&&\qq = \dbE\Big[|\tilde\G_\tau \br(\tau,X_\tau^\a,Y_\tau^\a)|^2- |\br(0,x_0,Y_0^\a)|^2 + 2\int_0^\tau \tilde\G_s^2 \br(s,X_s^\a,Y_s^\a)\bn\cdot\xi(s,X_s^\a,\pi_s^\a)ds\Big]\nonumber\\
&&\qq \leq C\e^4 + C\e^2 \dbE\bigg[\sup_{0\le t\le \t} |\tilde \G_t|^2\int_0^\tau [1+|Z^\a_s|^2]ds\bigg] \leq C\e^2.
\eea

{\bf Step 2.2.} Introduce two processes
\bea
\label{hatzeta}
\left.\ba{c}
\dis \hat\zeta_s := (\mathbf{I}_m-\bn\bn^\top)\zeta_s(\pi^\a_s) \in \big(\dbT_\dbV(s, X^\a_s, \pi^\a_s)\big)^d;\ms\\
\dis\hat\xi_s:= - \big[\pa_t\dbV+h_\dbV(\cd, \pa_x\dbV, \pa_{xx}\dbV, \a_s, \hat\zeta_s)\big](s,X^\a_s, \pi^\a_s).
\ea\right.
\eea
By \eqref{HJB}, \eqref{HJBr}, and by Step 1, we have
\bea
\label{cLVhatxi}
0\le  -  \cL \dbV(s, X^\a_s, \pi^\a_s)  \le  \bn(s, X^\a_s, \pi^\a_s) \cd \hat\xi_s.
\eea
Then, by taking expectation on both sides of \eqref{rtau} we have
\beaa
&& - \dbE\Big[\int_0^\t \tilde \G_s \cL \dbV(s, X^\a_s, \pi^\a_s) ds\Big]\le \dbE\Big[\int_0^\t \tilde \G_s \bn(s, X^\a_s, \pi^\a_s) \cd \hat\xi_s ds\Big] \\
&&\le  \dbE\Big[\int_0^\t \tilde \G_s \bn(s, X^\a_s, \pi^\a_s) \cd \xi_s ds\Big] + \dbE\Big[\int_0^\t \tilde \G_s  |\xi_s-\hat \xi_s| ds\Big]\\
&&=\dbE\Big[ \br (0, x_0, Y^\a_0) - \tilde \G_\t \br (\t, X^\a_\t, Y^\a_\t)\Big] + \dbE\Big[\int_0^\t \tilde \G_s  |\xi_s-\hat \xi_s| ds\Big].
\eeaa
Recall Remark \ref{rem-symmetric} (iii), we see that  $\dis \tr \big(( \bn\bn^\top \zeta)^\top \pa_y\bn( \bn\bn^\top \zeta)\big) =0$. Then, by \eqref{H}, 
\beaa
 |\xi_s-\hat \xi_s| \le C|\bn^\top \zeta(s, X^\a_s, \pi^\a_s)|,
 \eeaa
and thus, by  \eqref{zetaest},
\bea
\label{cLVest1}
&& - \dbE\Big[\int_0^\t \tilde \G_s \cL \dbV(s, X^\a_s, \pi^\a_s) ds\Big]\le \dbE\Big[\int_0^\t \tilde \G_s \bn(s, X^\a_s, \pi^\a_s) \cd \hat\xi_s ds\Big]\nonumber \\
&&\le\dbE\Big[ \br (0, x_0, Y^\a_0) - \tilde \G_\t \br (\t, X^\a_\t, Y^\a_\t)\Big] + C \dbE\Big[\int_0^\t \tilde \G_s |\bn^\top \zeta(s, X^\a_s, \pi^\a_s)|ds \Big]\nonumber\\
&&\le C\e^2 + C\bigg(\dbE\Big[\int_0^\t |\tilde \G_s \bn^\top \zeta(s, X^\a_s, \pi^\a_s)|^2ds \Big]\bigg)^{1/2}\leq C\e.
\eea

{\bf Step 2.3.} Fix another small constant $\d>0$. Since $\cL \dbV \le 0$, by \eqref{tauest} we have
\beaa
&& - \cL \dbV(0, x_0, y_0) = \dbE\Big[-{1\over \d}\int_0^{\t\wedge \d} \cL \dbV(0, x_0, y_0) ds  - {(\d-\t)^+\over \d} \cL \dbV(0, x_0, y_0) \Big]\\
&&= - {1\over \d}\dbE\Big[\int_0^{\t\wedge \d} \tilde \G_s \cL \dbV(s, X^\a_s, \pi^\a_s) ds\Big] -  \cL \dbV(0, x_0, y_0) \dbE\Big[(1-{\t\over \d})\1_{\{\t<\d\}}\Big]\\
&&\q + {1\over \d} \dbE\Big[\int_0^{\t\wedge \d} \big[\tilde \G_s \cL \dbV(s, X^\a_s, \pi^\a_s)-  \cL \dbV(0, x_0, y_0)\big]ds\Big]\\
&&\le {C\e\over \d}   +  \dbE\Big[\sup_{0\le t\le \t\wedge \d} \big|\tilde \G_t \cL \dbV(t, X^\a_t, \pi^\a_t)-  \cL \dbV(0, x_0, y_0)\big|\Big].
\eeaa
Since $\dbV\in C^{1,2}([0, T_0]\times \dbR^d;\cD^m_2)$, $H_\dbV$ is bounded  and uniform continuous. Then, for some  modulus of continuity function $\rho$ we have
\bea
\label{cLVest2}
\qq - \cL \dbV(0, x_0, y_0) \le {C\e\over \d}   + C \dbE\Big[\sup_{0\le t\le \t\wedge \d} \big[|\tilde \G_t - 1| + \rho\big(\d + |X^\a_t -x_0| + |\pi^\a_t-y_0|\big)\big]\Big].\hspace{1em}
\eea
Recall \eqref{tildeG}, \eqref{tildebbound}, and \eqref{tildeGest}, we have
\beaa
&&\dis \dbE\Big[\sup_{0\le t\le \t\wedge \d} |\tilde \G_t - 1| \Big] \le \dbE\Big[\sup_{0\le t\le \t\wedge \d} \big[(\tilde \G_t + 1)(|\int_0^t \tilde \si_s \cd dB_s| + \int_0^t [|\tilde b_s| + {1\over 2}|\tilde \si_s|^2] ds \big)\Big]\\
&&\dis\q \le C\sqrt{\d} + C\Big(\dbE\Big[\big(\int_0^{\t\wedge \d} |Z^\a_t|dt\big)^2\Big]\Big)^{1\over 2}\le  C\sqrt{\d} + C\sqrt{\d}\Big(\dbE\Big[\int_0^T |Z^\a_t|^2dt\Big]\Big)^{1\over 2}\le C\sqrt{\d};\\
&&\dis \dbE\Big[\sup_{0\le t\le \t\wedge \d}  |X^\a_t -x_0|\Big] \le C\sqrt{\d};\\
&&\dis \dbE\Big[\sup_{0\le t\le \t\wedge \d}  |\pi^\a_t -y_0|\Big] \le C\e^2 + \dbE\Big[\sup_{0\le t\le \t\wedge \d}  |Y^\a_t -y_0|\Big] \\
&&\dis\q \le C[\e^2 + \sqrt{\d}] + C\dbE\Big[\big(\int_0^{\t\wedge \d} |Z^\a_t|^2dt\big)^{1\over 2}\Big].
\eeaa
Then \eqref{cLVest2} implies
\bea
\label{cLVest3}
\dis - \cL \dbV(0, x_0, y_0) &\le& C\Big[{\e\over \d} +\sqrt{\d}\Big]  + C\rho(\d+\d^{1\over 3}) +  {C\over \d^{1\over 3}} \dbE\Big[\sup_{0\le t\le \t\wedge \d} [ |X^\a_t -x_0| +  |\pi^\a_t -y_0|] \Big] \nonumber\\
& \le& C\Big[{\e\over \d} +\sqrt{\d}   + \rho(\d+\d^{1\over 3}) +  {\e^2+\sqrt{\d}\over \d^{1\over 3}}\Big] +  {C\over \d^{1\over 3}}  \dbE\Big[\big(\int_0^{\t\wedge \d} |Z^\a_t|^2dt\big)^{1\over 2}\Big].
\eea

\ms
{\bf Step 2.4.} Recall \eqref{HVbound} and \eqref{cLVhatxi}. Then by \eqref{cLVest1} we have
\beaa
C\e  \ge \dbE\Big[\int_0^{\t\wedge \d} \tilde \G_s \bn(s, X^\a_s, \pi^\a_s) \cd \hat\xi_s ds\Big] \ge \dbE\Big[\int_0^{\t\wedge \d} \tilde \G_s \big[{c_{T_0}\over 2} |\hat \zeta_s|^2 - C |\hat \zeta_s| - C\big] ds\Big].
\eeaa
This, together with \eqref{tildeGest},  implies that 
\bea
\label{hatzetaest}
\dbE\Big[\int_0^{\t\wedge \d} \tilde \G_s|\hat \zeta_s|^2ds\Big] \le C[\e +\d].
\eea
By \eqref{zeta} and \eqref{hatzeta} we have 
\beaa
|Z^\a_t| \le |\zeta_t(\pi^\a_t)| + C \le |\hat \zeta_t| + |\bn^\top \zeta(t, X^\a_t,\pi^\a_t)| + C.
\eeaa
Then, by \eqref{zetaest}, \eqref{hatzetaest},  and \eqref{tildeGest}, we have
\beaa
&&\dis \dbE\Big[\big(\int_0^{\t\wedge \d} |Z^\a_t|^2dt\big)^{1\over 2}\Big] \le C\sqrt{\d}+C\dbE\Big[\big(\int_0^{\t\wedge \d} [|\hat \zeta_t|^2 + |\bn \cd \zeta(t, X^\a_t,\pi^\a_t)|^2]dt\big)^{1\over 2}\Big] \\
&&\dis\le C\sqrt{\d}+C\dbE\Big[\big(\sup_{0\le t\le \t} \tilde \G_t^{-1} \int_0^{\t\wedge \d} \tilde \G_t |\hat \zeta_t|^2 dt \big)^{1\over 2}\\
&&\dis \qq\qq\qq +\big(\sup_{0\le t\le \t} \tilde \G_t^{-2} \int_0^{\t} \tilde \G_t^2 |\bn \cd \zeta(t, X^\a_t,\pi^\a_t)|^2]dt\big)^{1\over 2}\Big] \\
&&\dis\le C\sqrt{\d}+C\Big(\dbE\Big[\int_0^{\t\wedge\d} \tilde \G_t |\hat \zeta_t|^2 dt + \int_0^{\t} \tilde \G_t^2 |\bn \cd \zeta(t, X^\a_t,\pi^\a_t)|^2dt\Big]\Big)^{1\over 2}\\
&&\dis\le C\sqrt{\d}+C\big(\e+\d +\e^2 \big)^{1\over 2} \le C[\sqrt{\e} + \sqrt{\d}].
\eeaa
Plug this into \eqref{cLVest3}, we get
\beaa
&&\dis - \cL \dbV(0, x_0, y_0)  \le C\Big[{\e\over \d} +\sqrt{\d}   + \rho(\d+\d^{1\over 3}) +  {\e^2+\sqrt{\d}\over \d^{1\over 3}}\Big] +  {C\over \d^{1\over 3}} [\sqrt{\e} + \sqrt{\d}].
\eeaa
By first send $\e\to 0$ and then $\d\to 0$, we obtain $- \cL \dbV(0, x_0, y_0)  \le 0$.
\qed

\section{The uniqueness of the classical solution}
\label{sect-uniqueness}
\setcounter{equation}{0}
We now turn to the uniqueness of the classical solution, including the verification result.
\begin{thm}
  \label{thm-unique} Let Assumption \ref{assum-coefficients} hold and $\dbV$ be defined by \eqref{Vtx}. 
  
  (i) Assume $\dbU \in C^{1,2}_0([0, T)\times \dbR^d;\cD^m_2)$ is a classical solution of \eqref{HJB} with terminal condition $\dbV(T,x) = \{g(x)\}$. Then $\dbU = \dbV$, and consequently \eqref{HJB} has a unique classical solution with terminal condition $\{g(x)\}$.

(ii) Assume further that the Hamiltonian $H_\dbU(\cd, \pa_x\dbU, \pa_{xx}\dbU)$ has an optimal argument: 
\beaa
\label{optimal1}
a^*=I^\dbU_1(t,x,y)\in A, \q \zeta^* = I^\dbU_2(t,x,y)\in (\dbT_\dbU(t,x,y))^d.
\eeaa
Moreover, recall Remark \ref{rem-rextension} and denote
\bea
\label{xi*}
\left.\ba{c}
\dis \tilde I^\dbU_3(t,x, y) := -\Big[\pa_t \dbU+h_\dbU(\cd, \pa_x \dbU, \pa_{xx}\dbU, I^\dbU_1, I^\dbU_2)\Big](t,x,y);\ms\\
\dis I^\dbU_3(t,x, y) := (\mathbf{I}_m-\bn_\dbU\bn_\dbU^\top)\ \tilde I^\dbU_3(t,x, y); 
 \ea\right.
 \eea
and assume, for given $(0, x_0, y_0)\in \dbG_\dbU$,  the following SDE has a strong solution: 
\bea
\label{X*}
\left.\ba{lll}
 \dis X_t^* = x_0 + \int_{0}^t b(\cd,I^\dbU_1)(s,X_s^*,\U_s^*)ds + \int_0^t \sigma(\cd,I^\dbU_1)(s,X_s^*,\U_s^*)dB_s;\ms \\
 \dis \U_t^* = y_0 + \int_0^t \Big[\pa_t \dbU+h^0_\dbU(\cd, \pa_x \dbU, \pa_{xx}\dbU, I^\dbU_1, I^\dbU_2) + I^\dbU_3\Big](s,X_s^*,\U_s^*)ds\ms \\
 \dis\qq + \int_{0}^t [\pa_x\dbU\sigma(\cd, I^\dbU_1)+I^\dbU_2](s,X_s^*,\U_s^*)dB_s.
 \ea\right.
 \eea
 Then, for $\a^*_t := I^\dbU_1(t,X_t^*,\U_t^*)$, we have $Y^{0,x_0, \a^*}_t=\U_t^* \in \dbV_b(t,X_t^*)$, $0\le t\le T$, a.s.  In particular, $Y^{0,x_0, \a^*}_0=y_0$.
\end{thm}

\begin{rem}
\label{rem-X*}
From Step 2 in the proof, especially \eqref{IU3=0} below, we see that \eqref{X*} actually becomes the following simpler and more natural SDE:
\bea
\label{X*2}
\left.\ba{lll}
 \dis X_t^* = x_0 + \int_{0}^t b(\cd,I^\dbU_1)(s,X_s^*,\U_s^*)ds + \int_0^t \sigma(\cd,I^\dbU_1)(s,X_s^*,\U_s^*)dB_s;\ms \\
 \dis \U_t^* = y_0 - \int_0^t f(\cd,\pa_x\dbU\sigma(\cd, I^\dbU_1)+I^\dbU_2, I^\dbU_1) (s,X_s^*,\U_s^*)ds\ms \\
 \dis\qq + \int_{0}^t [\pa_x\dbU\sigma(\cd, I^\dbU_1)+I^\dbU_2](s,X_s^*,\U_s^*)dB_s.
 \ea\right.
 \eea
 
 \vspace{-9mm}
\qed
\end{rem}

\begin{rem}
\label{rem-verification}
(i) Under the setting of above (ii), the $\a^*$ is an optimal argument (at least locally) for the scalarized optimization problem: $\sup_{\a\in \cA_0} \bn(0,x_0, y_0) \cd Y^{0,x_0,\a}_0$. We refer to Subsection \ref{sect-moving} below for more detailed analysis along this line.

(ii) When $\si$ is nondegenerate, by \eqref{X*2} we have
\beaa
&&\dis \U_t^* = y_0 +\int_{0}^t \Big[\big[\pa_x\dbU\sigma(\cd, I^\dbU_1)+I^\dbU_2\big] \si^{-1}(\cd,I^\dbU_1)\Big](s,X_s^*,\U_s^*)dX^*_s\\
&&\dis - \int_0^t \Big[f(\cd,\pa_x\dbU\sigma(\cd, I^\dbU_1)+I^\dbU_2, I^\dbU_1) + \big[\pa_x\dbU\sigma(\cd, I^\dbU_1)+I^\dbU_2\big] \si^{-1}b(\cd, I^\dbU_1)\Big] (s,X_s^*,\U_s^*)ds
 \eeaa
Then we may write $\U_t^*$ as a function of $X^*_{[0,t]}$, thus as a closed loop control $\a^*_t = \a^*(t, X^*_{[0,t]})$ is path dependent. Such path dependence appears often in multivariate setting. However, we note that $(X^*, \U^*)$ is jointly Markovian, so by adding the state $\U^*$, the optimal control $\a^*$ becomes Markovian, or more precisely state dependent. Therefore, the above verification theorem does help to construct Markovian optimal controls in this sense.  
\qed
\end{rem}

\no{\bf Proof of Theorem \ref{thm-unique}.} We proceed in three steps. Denote $T_\d:= T-\d$ for $\d>0$ small.

{\bf Step 1.} We first show that $\dbV(0, x_0) \subset \dbU(0,x_0)$. By the same arguments, we can also show that $\dbV(t, x) \subset \dbU(t,x)$ for all $(t,x)\in [0, T)\times \dbR^d$.

Fix $\d>0$ small and $\a\in \cA_0$.  Denote $(X^\a, Y^\a, Z^\a):= (X^{0,x_0,\a}, Y^{0,x_0,\a}, Z^{0,x_0,\a})$. Since $\dbV(T,x) = \{g(x)\}= \dbU(T, x)$, by Assumption \ref{assum-coefficients} and the continuity of $\dbU$, there exists $\phi_\d\in \dbL^2(\cF_{T_\d})$ such that $\phi_\d\in \dbU_b(T_\d, X^\a_{T_\d})$, a.s. and
\bea
\label{phid}
  \dbE\big[|Y^\a_{T_\d} - \phi_\d|^2\big] \le C\dbE\big[|Y^\a_{T_\d}-g(X^\a_T)|^2 +|g(X^\a_T)-g(X^\a_{T_\d})|^2+ |\phi_\d-g(X^\a_{T_\d})|^2\big]\to 0,\hspace{-2em}
\eea 
as $\d\to 0$. Recall \eqref{YT0} and set $(Y^{\a,\d}, Z^{\a,\d}) := (Y^{T_\d, \phi_\d; 0, x_0, \a}, Z^{T_\d, \phi_\d; 0, x_0, \a})$. Then by the standard BSDE estimates we have
\bea
\label{Yad}
\lim_{\d\to 0} |Y^{\a,\d}_0 - Y^\a_0|=0.
\eea

 As in \eqref{dr}, by standard It\^{o}'s formula, we have
\beaa
&&\dis d \br_{\dbU}(t, X^\a_t, Y^{\a,\d}_t) = \L(t, X^\a_t, Y^{\a,\d}_t, Z^{\a,\d}_t, \a_t) dt + \tilde Z^{\a,\d}_t  dB_t,\q\mbox{where}\\
&&\dis \L:= \td_t \br_\dbU + \td_x \br_\dbU \cd b - \td_y \br_\dbU\cd f + {1\over 2} \tr\big(\si^\top\td_{xx}\br_\dbU \si + 2 z^\top \td_{xy} \br_{\dbU}  \si+ z^\top \td_{yy} \br_\dbU z\big);\\
&&\dis \tilde Z^{\a,\d}:= \td_x \br_{\dbU} \si + (\td_y \br_\dbU)^\top Z^{\a,\d}. 
\eeaa
Denote
\beaa
\pi^{\a,\d}_t := \pi_{\dbU}(t, X^\a_t, Y^{\a,\d}_t),\q \zeta^{\a,\d}_t := Z^{\a,\d}_t - \bn_\dbU\bn_\dbU^\top Z^{\a,\d}_t.
\eeaa
 Then, by \eqref{LM} and since $\dbU$ is a classical solution of \eqref{HJB}, we have 
\beaa
\L(t, X^\a_t, \pi^{\a,\d}_t, \pa_x\dbU \si+\zeta^{\a,\d}_t, \a_t) = -\bn_\dbU \!\cd\! \big[\pa_t \dbU + h_\dbU\big(t, X^\a_t, \pi^{\a,\d}_t, \pa_x\dbU, \pa_{xx} \dbU, \a_t, \zeta^{\a,\d}_t\big)\big]\ge 0.
\eeaa
Note that, for appropriate processes $\tilde b, \tilde \si$,
\beaa
\L(t, X^\a_t, Y^{\a,\d}_t, Z^{\a,\d}_t, \a_t) - \L(t, X^\a_t, \pi^{\a,\d}_t, \pa_x\dbU\si+\zeta^{\a,\d}_t, \a_t) = -\big[\tilde b_t \br_\dbU + \tilde \si_t \tilde Z^{\a,\d}_t\big].
\eeaa
Here, due to the regularity of $\dbU\in C^{1,2}([0, T_\d]\times \dbR^d; \cD^m_2)$, as in \eqref{LMest} there exists a constant $C_\d>0$, which may depend on $\d$, such that for $0\le t\le T_\d$,
\beaa
|\tilde b_t|\le C_\d[1+|Z^{\a,\d}_t|^2],\q |\tilde \si_t|\le C_\d[1+|Z^{\a,\d}_t|].
\eeaa
 Then, for the $\tilde \G$ in \eqref{tildeG} we have
 \beaa
 d\Big(\tilde \G_t \br_{\dbU}(t, X^\a_t, Y^{\a,\d}_t)\Big) = \tilde \G_t \L(t, X^\a_t, \pi^{\a,\d}_t, \pa_x\dbU\si+\zeta^{\a,\d}_t, \a_t) dt + \tilde \G_t (\tilde Z_t^{\a,\d} - \br_\dbU \tilde\si_t)dB_t.
 \eeaa
Since $\br (T_\d, X^\a_{T_\d}, Y^{\a,\d}_{T_\d})=0$, a.s. then, 
\beaa
\br_\dbU(0, x_0, Y^{\a,\d}_0) = -\dbE\Big[\int_0^{T_\d} \tilde \G_t     \L(t, X^\a_t, \pi^{\a,\d}_t, \pa_x\dbU\si+\zeta^{\a,\d}_t, \a_t) \Big]\le 0.
\eeaa
That is, $Y^{\a,\d}_0 \in \dbU(0,x_0)$. Send $\d\to 0$, by \eqref{Yad} and the closedness of $\dbU(0,x_0)$, we have $Y^{\a}_0 \in \dbU(0,x_0)$. Moreover, since $\a\in \cA_0$ is arbitrary, we obtain $\dbV(0,x_0) \subset \dbU(0, x_0)$.

\ss
{\bf Step 2.} We next prove (ii) and show that in this case $\dbU(0,x_0) \subset \dbV(0,x_0)$. Indeed, consider an arbitrary $y_0\in \dbU_b(0,x_0)$. First by the It\^{o} formula Theorem \ref{thm-Ito} we see that $\U_t^* \in \dbV_b(t,X_t^*)$, $0\le t\le T$, a.s. In particular, this implies $\U^*_T = g(X^*_T)$. Note that, by the optimality of $I^\dbU_1, I^\dbU_2$, 
\beaa
h_\dbU(\cd, \pa_x \dbU, \pa_{xx}\dbU, I^\dbU_1, I^\dbU_2)(s,X_s^*,\U_s^*) = H_\dbU(\cd, \pa_x \dbU, \pa_{xx}\dbU, I^\dbU_1, I^\dbU_2)(s,X_s^*,\U_s^*).
\eeaa
Since $\dbU$ satisfies the PDE \eqref{HJB} and by \eqref{H}, at $(s, X^*_s, \U^*_s)\in \dbG_\dbU$ we have
 \bea
 \label{IU3=0}
\qq \bn_\dbU \cd \tilde I^\dbU_3=0;\q
\pa_t \dbU+h^0_\dbU(\cd, \pa_x \dbU, \pa_{xx}\dbU, I^\dbU_1, I^\dbU_2)+I^\dbU_3 =-f(\cd, \pa_x\dbU\sigma(\cd, I^\dbU_1)+I^\dbU_2, I^\dbU_1).
 \eea
This implies that $Y^{0,x_0, \a^*}_t = \U^*_t$. In particular, $y_0 = \U^*_0 = Y^{0,x_0, \a^*}_0 \in \dbV(0,x_0)$. Thus $\dbU_b(0,x_0)\subset \dbV(0, x_0)$, which implies that  $\dbU(0,x_0)\subset \dbV(0, x_0)$.

\ss
  {\bf Step 3.}  We now prove $\dbU(0,x_0) \subset \dbV(0,x_0)$ in the general case, without assuming the additional conditions in (ii).  As in Step 2 of Theorem \ref{thm-existence}, this step is more involved than Step 2 above, because it requires to construct approximately "optimal controls". Fix $(0,x_0, y_0)\in \dbG_\dbU$ and $ \d>0$. 
    Since $\dbU\in C^{1,2}([0, T_\d]\times\dbR^d;\cD^m_2)$, we assume $\br_\dbU$ is smooth in $O^{T_\d}_{\e_0}(\dbG_\dbU)$ for some $\e_0>0$. In the rest of this proof, let $C_\d$ be a generic constant which may depend on $\d$, more precisely on the $c_{T_\d}$ in \eqref{HVbound} and the regularity of $\dbU$ on $[0, T_\d]\times \dbR^d$. 
  
 Since $\dbU$ satisfies \eqref{HJB}, by \eqref{HVbound} there exist $\bar a_0\in A$ and $\bar\zeta^0\in (\dbT_\dbU(0,x_0, y_0))^d$ such that 
\bea
\label{zeta0}
|\bar\zeta^0|\le C_\d\q \mbox{and}  \q 0\le - \bn_\dbU\cd \big[\pa_t \dbU+ h_\dbU(\cd, \pa_x\dbU, \pa_{xx}\dbU, \bar a_0, \bar \zeta^0)\big](0,x_0,y_0)<\d. 
 \eea
  Set $\t_0:=0$, $\a^1_t \equiv \bar a_0$, $0\le t\le T_\d$, and define
   \beaa
   X_t^1 = x_0 + \int_0^t b(s,X_s^1, \a^1_s) ds + \int_0^t \sigma(s,X_s^1,\a^1_s)dB_s,\q 0\le t\le T_\d. 
   \eeaa 
 Recall Remark \ref{rem-rextension} and introduce random fields $(\xi^1, \zeta^1): [0, T_\d]\times \O \times \dbR^m \to (\dbR^m, \dbR^{m\times d})$: 
      \beaa
   \label{zeta1}
   \left.\ba{c}
\dis     \zeta^1_t(y) := \bar\zeta^0 - \bn_\dbU\bn_\dbU^\top\bar\zeta^0 (t,X_t^1,y),\q \xi^1_t(y):= \tilde\xi^1_t(y) - [\bn_\dbU \cd \tilde \xi^1_t] \bn_\dbU(t, X^1_t, y),\ms\\
\dis \mbox{where}\q     \tilde\xi^1_t(y) :=  -\big[\pa_t \dbU+ h_\dbU(\cd, \pa_x\dbU, \pa_{xx}\dbU, \a^1, \zeta^1)\big](t, X^1_t, y).
\ea\right.
  \eeaa
  Then  $\xi^1_t(y) \!\in\! \dbT_\dbU(t, X^1_t, y)$, $\zeta^1_t(y) \!\in\! (\dbT_\dbU(t, X^1_t, y))^d$, $\forall y\in \dbU_b(t, X^1_t)$, and $\xi^1, \zeta^1$ are uniformly Lipschitz continuous in $y$, with a Lipschitz constant depending on $\d$. Consider the SDE:
   \bea
    \label{Upsilon1} 
    \left.\ba{c}
    \dis \U_t^1 = y_0 + \int_0^t \Big[\pa_t \dbU+ h^0_\dbU(\cd, \pa_x\dbU, \pa_{xx}\dbU, \a^1_s, \zeta^1_s) + \xi^1\Big](s,X_s^1,\U_s^1)ds\\[1.9ex]
  \dis  + \int_0^t \Big[\pa_x\dbU(s, X^1_s, \U^1_s)\sigma(s, X^1_s, \a^1_s) + \zeta^1_s(\U_s^1)\Big] dB_s.
  \ea\right.
   \eea 
   By the It\^{o} formula Theorem \ref{thm-Ito} we have $\U^1_t \in \dbU_b(t, X^1_t)$, $0\le t\le T_\d$. Note that \eqref{zeta0} implies $\bn_\dbU(0, X^1_0,\U^1_0)\cd \tilde \xi^1_0(\U^1_0) \le \d$, and by our construction, $\a^1, \zeta^1$ and hence $\tilde \xi^1$ are continuous in $t$. We then set
   \beaa
   &\dis \t_1 := \inf\big\{t>\t_0:  \bn_\dbU(t, X^1_t,\U^1_t)\cd \tilde \xi^1_t(\U^1_t) \ge  2\d \Big\} \wedge T_\d.
   \eeaa
   
   Next, on $\{\t_1 < T_\d\}$, by measurable selection theorem, there exist $\cF_{\t_1}$-measurable random variables $\bar\a^1_{\t_1}\in A$ and $\bar\zeta^1_{\t_1} \in (\dbT_\dbV(\t_1, X^1_{\t_1}, \U^1_{\t_1}))^d$ such that 
   \beaa
 |\bar\zeta^1_{\t_1}|\le C_\d\q\mbox{ and}\q   0\le - \bn_\dbU\cd \big[\pa_t \dbU+ h_\dbU(\cd, \pa_x\dbU, \pa_{xx}\dbU, \bar \a^1_{\t_1}, \bar \zeta^1_{\t_1})\big](\t_1, X^1_{\t_1}, \U^1_{\t_1}) <\d. 
 \eeaa
Set  $\a^2_t \equiv \bar \a^1_{\t_1}$, $\t_1\le t\le T_\d$, and define
   \beaa
   X_t^2 = X^1_{\t_1} + \int_{\t_1}^t b(s,X_s^2, \a^2_s) ds + \int_{\t_1}^t \sigma(s,X_s^2,\a^2_s)dB_s,\q \t_1\le t\le T_\d. 
   \eeaa
   Similarly introduce, for $\t_1\le t\le T_\d$, 
   \beaa
   \label{zeta2}
   \left.\ba{c}
\dis     \zeta^2_t(y) := (\mathbf{I}_m - \bn_\dbU\bn_\dbU^\top)\bar\zeta^1_{\t_1}(t,X_t^2,y),\q \xi^2_t(y):= (\mathbf{I}_m-\bn_\dbU\bn_\dbU^\top)\tilde\xi^2_t(y)\bn_\dbU(t, X^2_t, y),\ms\\
\dis \mbox{where}\q     \tilde\xi^2_t(y) :=  -\big[\pa_t \dbU+ h_\dbU(\cd, \pa_x\dbU, \pa_{xx}\dbU, \a^2, \zeta^2)\big](t, X^2_t, y),
\ea\right.
  \eeaa
 and consider the SDE:
   \beaa
    \label{Upsilon2} 
    \left.\ba{c}
    \dis \U_t^2 =  \U^1_{\t_1} + \int_{\t_1}^t\Big[\pa_t \dbU+ h^0_\dbU(\cd, \pa_x\dbU, \pa_{xx}\dbU, \a^2_s, \zeta^2_s) + \xi^2\Big](s,X_s^2,\U_s^2)ds\\
  \dis  + \int_{\t_1}^t \Big[\pa_x\dbU(s, X^2_s, \U^2_s)\sigma(s, X^2_s, \a^2_s) + \zeta^2_s(\U_s^2)\Big] dB_s.
  \ea\right.
   \eeaa
Then $\U^2_t \in \dbU_b(t, X^2_t)$, $\t_1\le t\le T_\d$, and we may set
   \beaa
   \t_2 := \inf\big\{t>\t_1:  \bn_\dbU(t, X^2_t,\U^2_t) \cd \tilde\xi^2_t(\U^2_t)  \ge 2\d \Big\} \wedge T_\d.
   \eeaa

Repeat the arguments, we obtain a sequence $(\t_n, \a^n, \zeta^n, \tilde \xi^n, \xi^n,  X^n, \U^n)$, $n\ge 0$, satisfying the desired properties. We first show that $\t_n = T_\d$ for $n$ large enough, a.s. Indeed,  on $E_\d:= \cap_{n\ge 1}\{\t_n<T_\d\}$, we have, 
\beaa
\bn_\dbU(\t_n, X_{\t_n},\U_{\t_n}) \cd \tilde \xi^n_{\t_n}(\U_{\t_n})\le \d,\q  \bn_\dbU(\t_{n+1}, X_{\t_{n+1}},\U_{\t_{n+1}}) \cd  \tilde\xi^n_{\t_{n+1}}(\U_{\t_{n+1}})= 2\d,\q \forall n. 
\eeaa   
Then, for any $n$,
\beaa
\d \dbP(E_\d) \le \dbE\Big[\Big| \bn_\dbU(\t_{n+1}, X_{\t_{n+1}},\U_{\t_{n+1}}) \cd  \tilde\xi^n_{\t_{n+1}}(\U_{\t_{n+1}}) - \bn_\dbU(\t_n, X_{\t_n},\U_{\t_n}) \cd  \tilde \xi^n_{\t_n}(\U_{\t_n})\Big|\Big].
\eeaa
Send $n\to \infty$, by the desired regularity and in particular $|\zeta|\le C_\d$, we obtain $\dbP(E_\d) = 0$.

We now define
\beaa
(\a_t, \zeta_t, X_t, \U_t, \xi_t) := (\a^n_t,  \zeta^n_t, X^n_t, \U^n_t, \xi^n_t),\q t\in [\t_n, \t_{n+1}), n=0, 1, \cds.
\eeaa
Note that  $X_{T_\d}:= \lim_{t\uparrow T_\d} X_t$ and $\U_{T_\d}:=\lim_{t\uparrow T_\d} \U_t$ are well defined.
Define
\beaa
Z_t := \pa_x\dbU(t, X_t, \U_t)\sigma(t, X_t, \a_t) + \zeta_t(\U_t),\q \eta_t := [\bn_\dbU \cd \tilde \xi_t] \bn_\dbU(t, X_t, \U_t),\q 0\le t<T_\d.
\eeaa
Then, $|\eta|\le 2\d$, and by \eqref{H} and \eqref{Upsilon1}  we have
\beaa
 \U_t = y_0- \int_0^t \Big[f(s, X_s, \U_s, Z_s, \a_s) + \eta_s\Big]ds + \int_0^t Z_s dB_s,\q 0\le t\le T_\d.
\eeaa
Equivalently, we may rewrite it backwardly:
\beaa
\label{UBSDE}
 \U_t = \U_{T_\d} +\int_t^{T_\d} \Big[f(s, X_s, \U_s, Z_s, \a_s) + \eta_s\Big]ds - \int_t^{T_\d} Z_s dB_s,\q 0\le t\le T_\d.
\eeaa
Compare this with \eqref{YT0}, by standard BSDE estimates we have
\bea
\label{YTd}
\big|y_0 - Y^{T_\d, \U_{T_\d}; 0, x_0, \a}_0\big|^2= \big|\U_0 - Y^{T_\d, \U_{T_\d}; 0, x_0, \a}_0\big|^2\le C\dbE\Big[\int_0^{T_\d} |\eta_s|^2ds\Big] \le C\d^2.
\eea

Finally, fix an arbitrary $a_*\in A$, and extend $\a$ with $\a_t\equiv a_*$, $t\in [T_\d, T]$. Since $\dbV(T,x) = \{g(x)\}= \dbU(T, x)$, by Assumption \ref{assum-coefficients} and the continuity of $\dbU$, similarly to \eqref{phid} we have
\beaa
  &&\dis\dbE\big[|Y^\a_{T_\d} - \U_{T_\d}|^2\big] \\
  &&\dis \le C\dbE\Big[|Y^\a_{T_\d}-g(X^\a_T)|^2 + |g(X^\a_T)-g(X^\a_{T_\d})|^2 +|\U_{T_\d}-g(X^\a_{T_\d})|^2\Big]\le \rho(\d),
  \eeaa 
for some modulus of continuity function $\rho$, independent of $\a$. Then, by standard BSDE estimates again,
\beaa
\big|Y^{T_\d, \U_{T_\d}; 0, x_0, \a}_0- Y^\a_0\big|^2 = \big|Y^{T_\d, \U_{T_\d}; 0, x_0, \a}_0- Y^{T_\d, Y^\a_{T_\d}; 0, x_0, \a}_0\big|^2 \le  \dbE\big[|Y^\a_{T_\d} - \U_{T_\d}|^2\big] \le \rho(\d).
\eeaa
Combine this with \eqref{YTd}, we have
\beaa
|y_0 - Y^\a_0|\le C\d + \sqrt{\rho(\d)}.
\eeaa
Since $Y^\a_0\in \dbV(0,x_0)$ and $\d>0$ is arbitrary, we obtain $y_0\in \dbV(0,x_0)$.
\qed

We conclude this section with  a simple example where $\dbV$ is indeed a classical solution. 
\begin{eg}
\label{eg-classical} Set $d=1$, $m=2$, $A= \{a\in \dbR^2: |a|\le 1\}$, and 
\beaa
b=0, \q \si =1,\q f = f^0(t,x) + a,
\eeaa
where $f^0$ and $g$ are smooth and bounded. Then it is straightforward to check that
\beaa
\dbV(t, x) = \Big\{ y\in \dbR^2: |y-w(t,x)| \le T-t\Big\},
\eeaa
where $w = (w_1, w_2)^\top$ is the classical solution to the following heat equations:
\beaa
\td_t w_i + {1\over 2} \td_{xx} w_i + f^0_i =0,\q w_i(T,x) = g_i(x),\q i=1,2.
\eeaa
We shall prove in Appendix that $\dbV\in C^{1,2}_0([0, T)\times \dbR;\cD^2_2)$, and  the conditions in Theorem \ref{thm-unique} (ii) hold true. Then it follows from Theorems \ref{thm-existence} and \ref{thm-unique} that $\dbV$ is the unique classical solution of the HJB equation \eqref{HJBr}.
\qed
\end{eg}

\section{An application: the moving scalarization}
\label{sect-moving}
\setcounter{equation}{0}

Recall Remark \ref{rem-MV}, in particular \eqref{MVV0} and \eqref{MVdbV=V0} for the mean variance optimization problem. This problem is time inconsistent in the following sense. Consider the general setting \eqref{FBSDE} and \eqref{Vtx}. Given $(0, x_0)$ and  $\f \in C(\dbR^m; \dbR)$, let $\a^*_{[0, T]}$ be an optimal control for the problem 
\bea
\label{V0phi}
V_0 := \sup_{\a\in \cA_0} \f(Y^{0,x_0, \a}_0) = \sup_{y\in \dbV(0,x_0)} \f(y).
\eea
 If we follow $\a^*$ on $[0, t]$ and denote $X^*_t := X^{0,x_0, \a^*}_t$. Then $\a^*_{[t, T]}$ is not optimal for the optimization problem at $t$\footnote{Here we are using the notations heuristically. Rigorously we shall either consider $\esup$ in the left side or consider  $X^*$ and $\a_{[t, T]}$ in a pathwise manner.}:
\beaa
\label{Vtphi}
\sup_{\a_{[t, T]}} \f\big(Y^{t, X^*_t,\a_{[t, T]}}_t\big) = \sup_{y\in \dbV(t, X^*_t)} \f(y).
\eeaa
It was proposed in \cite{KMZ} to find a so called dynamic utility function $\Phi(t, \bx_{[0, t]}, y)$ such that $\Phi(0, x_0, y) = \f(y)$ and $\a^*_{[t, T]}$ remains optimal for the alternative optimization problem
\bea
\label{VtPhi}
\sup_{\a_{[t, T]}} \Phi\Big(t, X^*_{[0, t]}, Y^{t, X^*_t,\a_{[t, T]}}_t\Big) = \sup_{y\in \dbV(t, X^*_t)} \Phi\big(t, X^*_{[0, t]}, y\big).
\eea
In Subsection \ref{sect-MV} below we will find such an $\Phi$ for the mean variance problem explicitly. In the next subsection we first consider the case that $\f$ is linear.

\subsection{The linear scalarization}
\label{sect-linear}
When $\f$ is linear: $\f(y) = \l_0 \cd y$ for some $\l_0 \in \dbR^m$, we require $\Phi$ to be linear as well: $\Phi(t, \bx_{[0, t]}, y) = \L(t, \bx_{[0,t]}) \cd y$. This $\L$ is exactly the moving scalarization proposed in \cite{FR2022}. That is, we want to find $\L$ such that $\L(0,x_0) = \l_0$ and $\a^*_{[t, T]}$ is optimal for the problem:
\bea
\label{VtPhilinear}
 \sup_{\a_{[t, T]}} \L(t, X^*_{[0, t]})\cd Y^{t, X^*_t, \a_{[t, T]}}_t = \sup_{y\in \dbV(t, X^*_t)} \L(t, X^*_{[0, t]})\cd y.
\eea

Our set-valued HJB equation provides a solution to this interesting problem, provided that  \eqref{HJB} is wellposed in the sense of Theorem \ref{thm-unique} (ii) and $\dbV(t,x)$ is strictly convex. Consider a slightly more general setting by letting $\l: \dbR^d\to \dbR^m$ be such that $\l(x_0)=\l_0$. Assume without loss of generality that $|\l(x)|=1$ for all $x\in \dbR^d$. Since $\dbV(0,x)$ is compact and strict convex, we may find a unique optimal argument $y_\l(x) \in \dbV_b(0, x)$ for the problem: $V(0,x) := \sup_{y\in \dbV(0,x)} \l(x)\cd y$. Recalling $\dbU=\dbV$, we construct $X^*, \U^*, \a^*$ as in Theorem \ref{thm-unique} (ii) with initial data $(0, x, y_\l(x))\in \dbG_\dbV$. Assume further that $\si\in \dbR^d$ is nondegenerate, then as in Remark \ref{rem-verification} (ii) $\U^*$ is $\dbF^{X^*}$-progressively measurable and hence there exists $\L$ such that
\bea
\label{MovingScalarization}
\L(t, X^*_{[0,t]}) = \bn_\dbV(t, X^*_t, \U^*_t).
\eea 
We argue that this $\L$ is a desired moving scalarization. 

First,  since $\l(x) \cd y_\l(x) = \sup_{y\in \dbV(0,x)}\l(x)\cd y$ and $|\l(x)|=1$, we see that 
\beaa
\l(x) = \bn_\dbV(0,x, y_\l(x)) = \L(0, x).
\eeaa
 Next, from the construction in Theorem \ref{thm-unique} (ii), it is clear that $\U^*_t = Y^{t, X^*_t, \a^*_{[t, T]}}_t$. Then, since $\dbV(t, X^*_t)$ is convex, by \eqref{MovingScalarization} we see that 
\beaa
&&\dis \L(t, X^*_{[0, t]}) \cd Y^{t, X^*_t, \a^*_{[t, T]}}_t = \L(t, X^*_{[0, t]})  \cd \U^*_t \\
&&\dis =   \sup_{y\in \dbV(t, X^*_t)} \L(t, X^*_{[0, t]}) \cd y =\sup_{\a_{[t, T]}} \L(t, X^*_{[0, t]}) \cd Y^{t, X^*_t, \a_{[t, T]}}_t.
\eeaa
This exactly means $\a^*_{[t, T]}$ is an optimal control for the dynamic optimization problem \eqref{VtPhilinear}.

We remark that the mapping $\L: [0, T]\times C([0, T]; \dbR^d) \to \dbR^m$, which is path dependent in an adapted way, is time consistent in the following sense.  Consider the problem at time $0$ with initial condition $(x, \l)$. Let $X^*$ and $\L$ be as above, but denoted as $X^{0,x,\l, *}$ and $\L^{0,\l}$ to indicate their dependence on the initial conditions. Now fix $0<t<T$, consider the problem on $[t, T]$ with initial condition $X^{0,x,\l, *}_{[0, t]}$ and $\L^{0,\l}(t, \cd)$, we can easily see that the moving scalarization we find following the same procedure coincides with the original $\L$ found at time $0$: 
\beaa
\label{moving-consistency}
X^{t, X^{0,x,\l, *}_{[0, t]}, \L^{0,\l}(t, \cd)}_s = X^{0,x,\l, *}_s,\q \L^{t, \L^{0,\l}(t,\cd)}(s,\cd) = \L^{0,\l}(s,\cd),\q t\le s\le T.
\eeaa

\begin{rem}
\label{rem-nonconvex}
When $\dbV(t, X^*_t)$ is nonconvex, as in Example \ref{eg-nonconvex}, the $\L$ in \eqref{MovingScalarization} can be viewed as a local asymptotic moving scalarization in the following sense:
\beaa
\label{localMS}
\left.\ba{c}
\dis \L(t, X^*_{[0, t]}) \cd \U^*_t \ge \L(t, X^*_{[0, t]}) \cd y  - o(|y- \U^*_t|),\q \forall y\in \dbV(t, X^*_t);\q\mbox{or equivalently,}\ms\\
\dis\!\! \L(t, X^*_{[0, t]}) \!\cd\! Y^{t, X^*_t, \a^*_{[t, T]}}_t\! \ge \L(t, X^*_{[0, t]}) \!\cd\! Y^{t, X^*_t, \a_{[t, T]}}_t  \!\!\!- o\big(\big|Y^{t, X^*_t, \a_{[t, T]}}_t- Y^{t, X^*_t, \a^*_{[t, T]}}_t\big|\big),\forall \a_{[t, T]}.
\ea\right.
\eeaa
\qed
\end{rem}

\begin{rem}
\label{rem-nonlinear}
(i) When the $\f$ in \eqref{V0phi} is nonlinear, since $\dbV(0,x_0)$ is compact, one may still find an optimal argument $y_0\in \dbV(0, x_0)$  for the problem in the right side of \eqref{V0phi}. We emphasize that it is possible that $y_0\in \dbV_o(0, x_0)$ and such $y_0$ may not be unique. Fix an arbitrary $\a^0\in \cA_0$ and $Z^0$, for example $\a^0 \equiv a_0\in A$ and $Z^0\equiv 0$. Denote $X^0:= X^{0,x_0, \a^0}$ and
\beaa
\label{Y0}
\left.\ba{c}
\dis Y^0_t = y_0 - \int_0^t f(s, X^0_s, Y^0_s, Z^0_s, \a^0_s) ds + \int_0^t Z^0_s dB_s,\\
\dis \t_0 := \inf\big\{t\ge 0: (t, X^0_t, Y^0_t) \in \dbG_\dbV\big\}.
\ea\right.
\eeaa
It is clear that $\t_0 \le T$ and $(\t, X^0_\t, Y^0_\t)\in \dbG_\dbV$.  Applying Theorem \ref{thm-unique} (ii) on $(\t, X^0_\t, Y^0_\t)$ (assuming all the conditions are satisfied) and following the measurable selection theorem we may construct $\a^*$ on $[\t_0, T]$ with initial condition $(\t, X^0_\t, Y^0_\t)$. Then one can easily see that $\a^0\oplus_{\t_0} \a^*$ is an optimal argument for the left side of \eqref{V0phi}. That is, Theorem \ref{thm-unique} (ii) can help us to construct an optimal control for \eqref{V0phi} even when $\f$ is nonlinear. However, in this case it is not clear how to construct naturally a  (nonlinear) moving scalarization $\Phi$ as in \eqref{VtPhi}. In particular, when $\f$ has certain structure, for example the linear quadratic structure for the mean variance problem in Remark \ref{rem-MV}, we may naturally expect $\Phi$ to have the same structure, which will add the difficulty for constructing a desired $\Phi$.

(ii) For some nonlinear $\f$, it is possible to linearize it through certain transformation. Indeed, let $\psi$ be a diffeomorphism\footnote{We refer to \cite[Theorem A]{Gordon} for a characterization of diffeomorphisms.}  on $\dbR^m$ and set $\tilde \dbV(t,x) := \{\psi(y): y\in \dbV(t,x)\}$. Then 
\beaa
 \sup_{y\in \dbV(t,x)} \f(y) = \sup_{\tilde y\in \tilde \dbV(t,x)} \tilde \f(\tilde y),\q \mbox{where}\q \tilde \f(\tilde y) := \f(\psi^{-1}(\tilde y)).
 \eeaa
  If one can choose $\psi$ such that $\tilde \f$ is linear, then one can apply the analysis in this subsection to find a linear moving scalarization $\tilde \Phi$ for $\tilde \dbV$, which leads to a desired nonlinear moving scalarization for the original $\dbV$: $\Phi(t, X^*_{[0,t]}, y) := \tilde \Phi(t, X^*_{[0,t]}, \psi(y))$. We remark that $X^*$ stands for $X^{\a^*}$ for some optimal control $\a^*$, so it remains the same after the transformation. However, in this case $\bn_\dbV(t, X^*_t, \U^*_t)$ does not lead to a desired moving scalarization. 
   \qed 
\end{rem}

\subsection{The mean variance problem}
\label{sect-MV}
In this subsection we find a desired moving scalarization for the mean variance problem in Remark \ref{rem-MV}, by employing the idea in Remark \ref{rem-nonlinear} (ii). We first remark that in this case $\dbV$ is not bounded. However, since $\dbV$ is explicit as in \eqref{MVV}, we may still apply the results in Theorem \ref{thm-unique} (ii).

\begin{thm}
\label{thm-MV}
Consider the optimization problem \eqref{MVV0} and introduce:
\bea
\label{MVLamda}
\L(t, \bx_{[0,t]}) :=   {\l e^{T-t}\over  e^T -  \l (\bx_t-\bx_0)}\q\forall\bx\in C([0, T], \dbR)~\mbox{s.t.}~ \sup_{0\le t\le T}[\bx_t -\bx_0] < {1\over \l}e^T.
\eea
Then the following dynamic mean variance problem is time consistent:
\bea
\label{MVdynamic}
V_t := \esup_\a\Big\{ \dbE[X^\a_T|\cF_t] -  \tfrac{1}{2}\L(t, X^*_{[0,t]})  \text{\emph{Var}}(X^\a_T|\cF_t)\Big\}.
\eea
Here $X^*$ is the optimal trajectory for \eqref{MVV0} and it satisfies $\sup_{0\le t\le T}[X^*_t -x_0] < {1\over \l}e^T$, a.s. Moreover, the optimal control and the optimal value are\footnote{The optimal control $\a^*$ is the same as the static optimal control in \cite[Theorem 3.3 (A)]{PP}, with the correspondence $\a^*_t = X^*_t u^s_*(t, X^*_t)$. However, our $V_t$ is neither equal to the static optimal value nor to the dynamic optimal value in \cite[Theorem 3.3]{PP}, except that at $t=0$ it is equal to the static optimal value there.}:
\bea
\label{MVoptimal2}
\left.\ba{c}
\dis\a^*_t =  - X^*_t + x_0 +\tfrac{1}{ \l} e^T;\ms\\
\dis  V_t =  \tfrac{1}{2}(1+e^{t-T})X_t^* + \tfrac{1}{2}(1-e^{t-T})x_0  + \tfrac{e^T}{2\l}(1-e^{t-T}). 
\ea\right.
\eea
\end{thm}

\proof In light of Remark \ref{rem-nonlinear} (ii), we introduce an obvious diffeomorphism 
\bea
\label{MVpsi}
\psi(y_1, y_2) := (y_1, \tilde y_2) := (y_1, y_2-|y_1|^2).
\eea
Then, by \eqref{MVV} we have
\begin{equation}
  \begin{aligned}
    \label{MVtildeV}
      \tilde \dbV(t,x) :=& \Big\{\psi(y): y\in \dbV(t,x)\Big\}
      = \Big\{ ( y_1, \tilde y_2):
      y_1\in \dbR, \tilde y_2 \ge \phi_1(t)(y_1 - x)^2\Big\}, \\
      \tilde \dbV_b(t,x) =&
      \Big\{ ( y_1, \tilde y_2): y_1\in \dbR, \tilde y_2 =
      \phi_1(t)(y_1 - x)^2\Big\},~\mbox{where}~ \phi_1(t) := {1\over e^{T-t}-1}.
  \end{aligned}
\end{equation}
Note that $\tilde \dbV$ is convex, so the concern in Remark \ref{rem-nonconvex} is irrelevant and we are  finding a true moving scalarization.  
We shall denote $\tilde y = (y_1, \tilde y_2)$, and for $\tilde y\in \tilde \dbV_b(t,x)$, clearly it suffices to specify $y_1$. Moreover, recall \eqref{MVdbV} and denote $\tilde Y := \psi(Y) = (Y^1, \tilde Y^2)$, by the standard It\^{o} formula we have
\beaa 
Y^{t,x,\a, 1}_s = X^{t,x,\a}_T - \int_s^T Z^{t,x,\a, 1}_r dB_r,~ \tilde Y^{t,x,\a, 2}_s = \int_s^T |Z^{t,x,\a, 1}_r|^2 dr- \int_s^T \tilde Z^{t,x,\a, 2}_r dB_r.
 \eeaa
 That is, in light of \eqref{FBSDE} and denoting $\tilde z = (z_1, \tilde z_2)$,
 \bea
 \label{MVtildef}
 \tilde f(t,x, \tilde y, \tilde z, a) = (0, |z_1|^2)^\top.
 \eea

Given \eqref{MVtildeV}, one can easily compute that
\bea
\label{MVn}
\left.\ba{c}
\dis \tilde \bn(t,x, \tilde y) := \bn_{\tilde \dbV}(t,x, \tilde y) =  {1\over \phi_3} \begin{bmatrix}
 \phi_2 \\
    -1
  \end{bmatrix}(t,x,\tilde y),\q (t,x,\tilde y)\in \dbG_{\tilde \dbV},\\
\dis  \mbox{where}\q  \phi_2(t,x, \tilde y):= 2\phi_1(t)(y_1-x),\q \phi_3:= \sqrt{1+|\phi_2|^2}.
\ea\right.
\eea
Next, fix $(t, x, \tilde y)\in \dbG_{\tilde \dbV}$ and set $\U(x') := \big(y_1, \phi_1(t) (y_1 - x')^2\big)^\top$, $x'\in \dbR$. Clearly $ \U(x') \in \tilde\dbV_b(x')$ for all $x'$ and ${d\over dx'}\U(x')\Big|_{x'=x} = (0, -\phi_2(t,x,\tilde y))^\top$. Then by \eqref{patxV} and \eqref{paxVU} we have
\beaa
 \pa_x \tilde \dbV(t,x,\tilde y) = \Big((0, -\phi_2)^\top \cd \tilde \bn\Big)\tilde \bn(t,x,\tilde y) =  {\phi_2\over \phi_3^2} \begin{bmatrix}
 \phi_2 \\
    -1
  \end{bmatrix}(t,x,\tilde y).
 \eeaa
 The right side of above and \eqref{MVn} provide natural extensions of $\tilde\bn$ and $\pa_x\tilde \dbV$ on $[0, T]\times \dbR \times \dbR^2$. Then by \eqref{payfhat} and \eqref{paxfhat} we may compute straightforwardly that, at $(t,x,\tilde y)\in \dbG_{\tilde \dbV}$,
\beaa
 \pa_{xx} \tilde \dbV = {2\phi_1\over \phi_3^6} \begin{bmatrix}
 -2\phi_2 \\
    1-|\phi_2|^2
  \end{bmatrix},\q
  \pa_x\tilde\bn = -\frac{2\phi_1}{\phi_3^5}
  \begin{bmatrix}
    1\\
    \phi_2
  \end{bmatrix},\q  
  \pa_{\tilde y}\tilde \bn   = \frac{2\phi_1}{\phi_3^5}
  \begin{bmatrix}\vspace{1em}
    1\\
    \phi_2
  \end{bmatrix}
  \begin{bmatrix}
    1 &\q \phi_2
  \end{bmatrix}.
\eeaa
Moreover, as the tangent space is one dimensional, it is clear that
\beaa
 &\dis \zeta \in \dbT_{\tilde\dbV}(t,x,\tilde y) \iff
  \exists \zeta_0\in\dbR ~\text{such that}~
  \zeta = \zeta_0
  \begin{bmatrix}
    1 \\
    \phi_2
  \end{bmatrix}.
\eeaa
Then, recalling \eqref{H} and \eqref{MVtildef}, we may compute straightforwardly that
\beaa
&&\dis \tilde \bn \cd h_{\tilde \dbV}(t, x, \tilde y, \pa_x\tilde \dbV, \pa_{xx}\tilde \dbV, a,\zeta) \\
&&\dis = \tilde \bn \cd\Big[a \pa_x\tilde \dbV + {a^2\over 2} \pa_{xx}\tilde \dbV +  (0, |a (\pa_x\tilde \dbV)_1 + \zeta_1|^2)^\top\Big] -\Big[a\zeta^\top\pa_x\tilde \bn + {1\over 2}\zeta^\top\pa_{\tilde y}\tilde \bn\zeta\Big]\\
&&\dis = {\phi_2\over \phi_3^2}a - {\phi_1\over \phi_3^5} a^2 - {1\over \phi_3}( {\phi_2^2\over \phi_3^2}a+\zeta_0)^2+ \frac{2\phi_1}{\phi_3^3} a\zeta_0 -  \frac{\phi_1}{\phi_3} \zeta_0^2\\
&&\dis= - {1\over \phi_3}\Big[(1+\phi_1)\zeta_0^2 + {\phi_1+\phi_2^4\over \phi_3^4}a^2 - {2(\phi_1-\phi_2^2)\over \phi_3^2} a\zeta_0 - \phi_2 a\Big].
\eeaa
This is quadratic in $(a, \zeta_0)$, and one may obtain immediately the optimal arguments:
\bea
\label{MVI12}
\left.\ba{c}
\dis I^{\tilde \dbV}_1 = a^* = {(1+\phi_1)\phi_2\over 2\phi_1} = (1+\phi_1(t))(y_1-x),\\
\dis \zeta_0^*= {\phi_2(\phi_1-\phi_2^2)\over 2\phi_1\phi_3^2},\q \mbox{and thus}\q I^{\tilde \dbV}_2 = \zeta^* = {\phi_2(\phi_1-\phi_2^2)\over 2\phi_1\phi_3^2}  \begin{bmatrix}
    1 \\
    \phi_2
  \end{bmatrix}.
\ea\right.
\eea

We next derive \eqref{X*} for $\tilde\dbV$, with the solution denoted as $(X^*, \tilde \U^*) = (X^*, \U^{*,1}, \tilde \U^{*,2})$. Since by Theorem \ref{thm-unique} we have $\tilde\U^*_t\in \tilde\dbV_b(t, X^*_t)$, it suffices to specify the equations for $(X^*, \U^{*,1})$. Note that, with $(\cd)_1$ denoting the first component, 
\beaa
\Big(\pa_x \tilde \dbV \si(\cd, I^{\tilde \dbV}_1) + I^{\tilde \dbV}_2\Big)_1(t,x,\tilde y) = {\phi_2^2\over \phi_3^2} \times {(1+\phi_1)\phi_2\over 2\phi_1} + {\phi_2(\phi_1-\phi_2^2)\over 2\phi_1\phi_3^2} = \phi_1(t) (y_1-x).
\eeaa
Note further that $\tilde f_1=0$. Then, by recalling \eqref{X*2} in Remark \ref{rem-X*} and \eqref{MVoptimal}, we have
\bea
\label{MVX*}
\left.\ba{lll}
\qq \dis X_t^* = x_0 + \int_{0}^t (1+\phi_1(s))(\U_s^{*,1}-X^*_s)ds + \int_0^t (1+\phi_1(s))(\U_s^{*,1}-X^*_s)dB_s;\\\qq
 \dis \U_t^{*,1} = x_0 + {1\over \l}[e^T-1] + \int_{0}^t \phi_1(s)(\U_s^{*,1}-X^*_s)dB_s.
 \ea\right.
 \eea
 Thus:
 \bea
\label{MVU-X}
\qq \U_t^{*,1} - X^*_t =  {1\over \l}[e^T-1] - \int_{0}^t (1+\phi_1(s))(\U_s^{*,1}-X^*_s)ds - \int_0^t (\U_s^{*,1}-X^*_s)dB_s.
\eea
Clearly we can solve this and hence \eqref{MVX*} explicitly. More relevantly for the moving scalarization, as in Remark \ref{rem-verification} (ii) we may rewrite \eqref{MVU-X} as
\beaa
\U_t^{*,1} - X^*_t = {1\over \l}[e^T-1] - \int_{0}^t {1\over 1+\phi_1(s)} d X^*_s -  \int_{0}^t \phi_1(s)(\U_s^{*,1}-X^*_s)ds.
\eeaa
Then, denoting $\G_t := e^{\int_0^t \phi_1(s)ds} = {e^T-1\over e^T-e^t}$, 
\beaa
\G_t(\U_t^{*,1}-X^*_t) &=&  {1\over \l}[e^T-1] - \int_{0}^t {\G_s\over 1+\phi_1(s)}  d X^*_s={1\over \l}[e^T-1] - \int_{0}^t (1-e^{-T}) d X^*_s \\
&=&  {1\over \l}[e^T-1] - (1-e^{-T}) X^*_t + (1-e^{-T}) x_0.
\eeaa
Thus
\bea
\label{MVX*sol}
\U_t^{*,1}-X^*_t = {1\over \l}(e^T-e^t) - (1-e^{t-T})(X^*_t -x_0).
\eea

By abusing the notation $\L$ with the previous subsection, our goal in this subsection is to find a moving scalarization $\L_t:=\L(t, X^*_{[0,t]})$ such that the following dynamic problem 
\bea
\label{MVtildeVconsistent}
 \sup_{\tilde y \in \tilde \dbV(t, X^*_t)} \Big(y_1 - {\L_t\over 2} \tilde y_2 \Big)~\mbox{is time consistent}.
\eea
From the analysis in the previous subsection, this implies that $(1, -{\L_t\over 2})^\top$ is parallel to $\tilde\bn (t, X^*_t, \U^*_t)$. By \eqref{MVn}, this implies that $-{\L_t\over 2} = {-1\over \phi_2(t, X^*_t, \U^*_t)}$. Thus, by \eqref{MVX*sol},
\beaa
\L(t, X^*_{[0,t]}) =  \L_t = {2\over \phi_2(t, X^*_t, \U^*_t)} = {e^{T-t}-1\over  \U_t^{*,1}-X^*_t} = {\l e^{T-t}\over  e^T -  \l (X^*_t-X^*_0)}.
\eeaa
This proves \eqref{MVLamda}. We remark that, by \eqref{MVU-X} clearly $\U_t^{*,1}-X^*_t >0$, then it follows from \eqref{MVX*sol} that $\sup_{0\le t\le T}[X^*_t -x_0] < {1\over \l}e^T$, a.s. 

For the original $\dbV$, by \eqref{MVtildeVconsistent} the following dynamic problem is time consistent:
\beaa
 \sup_{ y \in  \dbV(t, X^*_t)}\Phi(t, X^*_{[0, t]}, y),\q\mbox{where}\q  \Phi(t, X^*_{[0, t]}, y):= y_1 +{ \L(t, X^*_{[0,t]}) \over 2} |y_1|^2 -{ \L(t, X^*_{[0,t]}) \over 2} y_2.
 \eeaa
 This is clearly equivalent to the time consistency of  the dynamic problem \eqref{MVdynamic}.

 Finally, plugging \eqref{MVX*sol} into the first line of \eqref{MVI12}, we obtain the expression of $\a^*$ in \eqref{MVoptimal2} immediately. Moreover,  by \eqref{MVtildeV}, \eqref{MVLamda},
 and \eqref{MVX*sol} we have
 \begin{equation*}
   \begin{aligned}
     &V_t = \U_t^{*,1} - \tfrac{\L_t}{2} \tilde\U_t^{*,2} =\U_t^{*,1} -{e^{T-t}-1\over 2( \U_t^{*,1}-X^*_t)}\times \phi_1(t)(\U_t^{*,1} - X_t^*)^2
     \\&    =
     \tfrac{1}{2} (\U_t^{*,1} - X_t^*) + X_t^*=\tfrac{1}{2}(1+e^{t-T})X_t^* + \tfrac{1}{2}(1-e^{t-T})x_0 +\tfrac{e^T}{2\l}(1-e^{t-T}).
   \end{aligned}
 \end{equation*}
 This proves \eqref{MVoptimal2}, and completes the proof of the theorem.
 \qed

\section{Further discussions}
\label{sect-further}
\setcounter{equation}{0}
\subsection{The case with nondegenerate terminal}
\label{sect-nondegenerate}
As pointed out in Remark \ref{rem-target} (ii), given a general $\dbG: \dbR^d \to \cD^m_0$, we may define $\dbV$ by \eqref{VtxG}. This is equivalent to
\bea
\label{VtxG2}
 \dbV(t,x) := {\rm cl} \big\{Y^{T, \phi; t,x, \a}_t:  \a\in \cA_t, \phi \in \dbL^2(\cF^t_T) ~\mbox{s.t.} ~ \phi\in \dbG(X^{t,x,\a}_T),~\mbox{a.s.}\big\}.
  \eea
Then we have

\begin{thm}
\label{thm-nondegenerate}
Let Assumption \ref{assum-coefficients} (i), (ii) hold, and $\dbG$ is bounded and continuous.  Assume the $\dbV$ defined by \eqref{VtxG} or \eqref{VtxG2} is in $C^{1,2}_0([0, T)\times \dbR^d; \cD^m_2)$. Then $\dbV$ is the unique classical solution of the HJB equation \eqref{HJB} with terminal condition $\dbV(T,x) = \dbG(x)$.
\end{thm}
The proof is essentially the same as in the previous sections, we thus omit it. In particular, when $\dbG(x)\in \cD^m_2$ and $\dbV\in C^{1,2}([0, T]\times \dbR^d; \cD^m_2)$, the proof is actually slightly easier.

\subsection{Comparison with Soner-Touzi \cite{ST2003}}
\label{sect-ST}
In the contexts of stochastic target problem, \cite{ST2003} derived a geometric equation to characterize the reachable set of the problem. This work is very closely related to our problem. In this subsection we provide some detailed analyses on the connection and the differences between the two works. We shall introduce their approach, but in our contexts and using our notations, and all the discussions are heuristic.

We first note that, the stochastic target problem \eqref{targetV} (or the more general one \eqref{VtxG}) can be rewritten equivalently as:
\beaa
\label{hatVtarget}
\wh\dbV(t) := \Big\{(x, y)\in \dbR^{d+m}: \exists (\a, Z)~\mbox{such that}~ Y^{t,x,y, \a, Z}_T = g(X^{t,x,\a}_T), \mbox{a.s.}\Big\}.
\eeaa
Here $(X^{t,x,\a},  Y^{t,x,y, \a, Z})$ becomes a $d+m$-dimensional controlled state process with control $(\a, Z)$. It is clear that $\wh \dbV$ and our $\dbV$ are equivalent in the following sense:
\bea
\label{hatVV}
\left.\ba{c}
\dis \wh\dbV(t) = \big\{(x, y): x\in \dbR^d, y\in \dbV(t,x)\big\},\q \wh\dbV_b(t) = \big\{(x, y): x\in \dbR^d, y\in \dbV_b(t,x)\big\};\ms\\
\dis\mbox{and}\q \dbV(t,x) = \big\{y: (x, y)\in \wh \dbV(t)\big\},\q \dbV_b(t, x) = \big\{y: (x,y)\in \wh\dbV_b(t)\big\}.
\ea\right.
\eea
Then $\dbG_{\wh\dbV} = \dbG_\dbV$.  Naturally we may define, for some $\e>0$,
\beaa
\label{hatnr}
\left.\ba{c}
\dis \bn_{\wh\dbV}(t,x,y) := \bn_{\wh \dbV(t)}(x, y)\in \dbR^{d+m}, \q (t, x, y)\in \dbG_{\wh\dbV};\\
\dis  \br_{\wh \dbV}(t,x,y) := \br_{\wh\dbV(t)}(x, y)\in \dbR,\q  (t,x,y) \in O_\e(\dbG_{\wh\dbV}).
\ea\right.
\eeaa
The work \cite{ST2003} characterized the square of the distance function\footnote{The reason to consider the squared function is that, in the degenerate case, $\br_{\wh\dbV}$ is typically not smooth while $|\br_{\wh\dbV}|^2$ is. In the nondegenerate case as in this paper, actually one may study $\br_{\wh\dbV}$ directly.}: $\eta(t,x,y) := {1\over 2} |\br_{\wh\dbV}(t,x,y)|^2$ by the following PDE: denoting $\hat y := (x, y)$ and noting the time change in \cite{ST2003},
\bea
\label{HJBrhat}
\left.\ba{lll}
\dis \td_{\hat y} \td_t \eta(t,\hat y) + \td_{\hat y}\Big[F(t, \hat y, \td_{\hat y} \eta(t,\hat y), \td_{\hat y\hat y}\eta(t,\hat y))\Big]  =0,\q (t,\hat y)\in \dbG_{\wh\dbV},\q\mbox{where}\ms\\
\dis  F(t, \hat y, \td_{\hat y} \eta, \td_{\hat y\hat y} \eta) := \inf_{(a,z)\in \cN(t, \hat y, \td_{\hat y}\eta)}\Big[b(t,x,a)\cd \td_x \eta(t, \hat y) - f(t,\hat y, z, a) \cd \td_y \eta(t, \hat y)\ms\\
\dis  \q + {1\over 2}\tr\big(\si^\top(t,x,a)\td_{xx}\eta(t, \hat y) \si(t,x,a) + 2 z^\top \td_{xy}\eta(t, \hat y) \si(t,x,z) + z^\top \td_{yy}\eta(t, \hat y) z\big)\Big],\ms\\
\dis  \cN(t, \hat y, \td_{\hat y}\eta):= \Big\{(a, z): [\si^\top(t,x,a), z^\top]^\top [\si^\top(t,x,a), z^\top] \td_{\hat y} \eta(t, \hat y)=0\Big\}.
\ea\right.
\eea 
 We remark that \eqref{HJBrhat} holds only on $\dbG_{\hat\dbV}$, and thus it is also not a standard PDE.

We first note that $\hat \dbV$ is a function of $t$ only, so it does not invoke the set-valued It\^{o} formula. While this may seem to be technically easier, the set-valued It\^{o} formula has independent interest and is one of the main contributions of this paper. For example, it provides  microstructure of the flow on the boundary surface, as we see in Theorem \ref{thm-unique} (ii). More importantly, the roles of $x$ and $y$ are symmetric in $\hat \dbV(t)$. In many applications, however, the input $x$ and the output $y$ play different roles: $x$ is the observed state, while the value $y$ (more precisely the set of $y$) is our objective, so it is more natural to study $\dbV(t,x)$.

Technically, the major difference is that, as we see in Example \ref{eg-hatV} below,
\beaa
\label{hatrnotr}
\br_{\wh\dbV}(t,x,y) \neq \br_\dbV(t,x,y).
\eeaa
In general, recalling \eqref{hatVV} we have
\beaa
\label{hatr<r}
\left.\ba{lll}
\dis |\br_{\wh\dbV}(t,x,y)|^2 = \inf_{(\tilde x, \tilde y)\in \wh\dbV_b(t)}\big[|x-\tilde x|^2 + |y-\tilde y|^2\big] =\inf_{\tilde x\in \dbR^d} \inf_{\tilde y\in \dbV_b(t, \tilde x)}\big[|x-\tilde x|^2 + |y-\tilde y|^2\big] \ms\\
\dis\qq\qq\q ~= \inf_{\tilde x\in \dbR^d}\big[|x-\tilde x|^2 +  |\br_{\dbV}(t,\tilde x, y)|^2\big] \le |\br_{\dbV}(t,x,y)|^2
\ea\right.
\eeaa

\begin{eg}
\label{eg-hatV}
Set $d=m=1$ and consider time invariant set values:
\beaa
\dbV(x) = [x-1, x+1]\subset \dbR,\q \wh\dbV := \{(x, y): x\in \dbR, y\in [x-1, x+1]\} \subset \dbR^2.
\eeaa
Clearly $\dbV_b(x) = \{x-1, x+1\}$. One can easily verify that, 
\beaa
\dis \br_\dbV(x, y) = y - (x+1),\q \mbox{for}~ y\approx x+1;\q\mbox{and}\q \br_\dbV(x, y) =  (x-1) - y,\q \mbox{for}~ y\approx x-1;\\
\dis \br_{\wh\dbV}(x, y) = {y - (x+1)\over \sqrt{2}},\q \mbox{for}~ y\approx x+1;\q\mbox{and}\q \br_{\wh\dbV}(x, y) =  {(x-1) - y\over \sqrt{2}},\q \mbox{for}~ y\approx x-1.
\eeaa
We also observe directly from above that, although $\br_\dbV = \br_{\wh\dbV} =0$ on $\dbG_\dbV=\dbG_{\wh\dbV}$, their derivatives are in general not equal.
\qed
\end{eg}

Consequently, although both \eqref{HJBr} and \eqref{HJBrhat} characterize the same set $\dbG_\dbV = \dbG_{\wh\dbV}$, the two equations are different. This is partially explained by the above observation that $\br_{\dbV}$ and $\br_{\wh\dbV}$ have different derivatives on $\dbG_\dbV$. At below we provide more detailed calculation for the set-valued heat equation in Example \ref{eg-paV} (ii), but with $d=m=1$. 

\begin{eg}
\label{eg-heat}
Set $d=m=1$, $b = 0$, $\si=1$, $f=0$, and the terminal $\dbG(x) = [-\psi(x), \psi(x)]$, where $\psi: \dbR\to (0, \infty)$ is smooth. Then, similar to Example \ref{eg-paV} (ii), we have
\beaa
\dbV(t,x) = [-u(t,x), u(t,x)],\q \dbV_b(t,x) = \{-u(t,x), u(t,x)\},
\eeaa
where $u$ is the unique classical solution of the heat equation
\bea
\label{heat1d}
\pa_t u + {1\over 2} \pa_{xx} u =0,\q u(T,x) = \psi(x).
\eea
We shall prove in Appendix that $\wh \br := \br_{\wh\dbV}$ satisfies the following equation:
\bea
\label{heathat}
\td_t \wh \br  + {1\over 2}\Big[ \td_{xx} \wh \br  - 2 \td_{xy}\wh \br {\td_x \wh \br \over \td_y\wh \br } +  \td_{yy}\wh \br \big|{\td_x \wh \br \over \td_y \wh \br }\big|^2\Big]=0,\q\mbox{on}~ \dbG_{\wh\dbV}.
\eea

\vspace{-3mm}
\qed
\end{eg}

In this scalar case, by Remark \ref{rem-HJB} (i) we see that the set-valued HJB equation \eqref{HJBr} reduces back to the standard PDE for $\ol v(t,x) = u(t,x)$ and $\ul v(t,x) = - u(t,x)$, both of which identify with the heat equation \eqref{heat1d}. So \eqref{HJBr} is indeed a natural extension of the HJB equation to the multivariate case. The equation \eqref{heathat}, however, is quite different from \eqref{heat1d}. So in this sense, it is more natural to study \eqref{HJBr} than to study \eqref{HJBrhat}.

Another advantage of \eqref{HJBr} is that, as we saw in Section \ref{sect-moving}, the normal vector $\bn_\dbV(t, X^*_t, \U^*_t)$ provides naturally a moving scalarization for the time inconsistent multivariate optimization problem. The vector $\bn_{\wh\dbV}$ (at certain optimal paths) does not serve for this purpose. In fact, $\bn_{\wh\dbV}\in \dbR^{d+m}$, while a moving scalarization $\L$ is by nature $m$-dimensional.

Finally, we remark that \cite[Theorem 2.1]{ST2003} showed that $\wh\dbV$ is the unique classical solution of \eqref{HJBrhat}, provided that optimal controls exist, in the same spirit of Theorem \ref{thm-unique} (ii). We instead proved the existence and uniqueness under weak conditions in Theorems \ref{thm-existence} and \ref{thm-unique} (i).

\subsection{Comparison with Ararat-Ma-Wu \cite{AMW}}
\label{sect-AMW}
Mainly motivated by dynamic set-valued risk measures for multi-asset  financial models,  \cite{AMW} studied the following set-valued BSDE:
\bea
\label{setBSDE}
\dbY_t = \dbE\Big[\dbG(B_T) + \int_t^T \dbF(s, B_s, \dbY_s)ds\big|\cF_t\Big].  
\eea
Here, denoting by $\cD^m_{cc}$ the space of convex compact sets $\dbD\in \cD^m_0$,  the terminal $\dbG: \dbR^d\to \cD^m_{cc}$, and the driver $\dbF: [0, T]\times \dbR^d \times \cD^m_{cc} \to \cD^m_{cc}$ (abusing the notation $\dbF$ here).  We note that \cite{AMW} actually allows $\dbG$ and $\dbF$ to depend on the paths of $B$. By relying on the sophisticated set-valued stochastic analysis, especially the Hukuhara difference, \cite{AMW} established the wellposedness of the above set-valued BSDE. The general case that $\dbF$ depends on $\dbZ$, and the martingale representation with the term $\dbZ_t dB_t$ seem to be a quite remote goal. 

Formally, the set-valued BSDE \eqref{setBSDE} is associated with our set-valued HJB equation \eqref{HJB} in the case $x_0=0$, $b=0$, $\si=1$, $f = f(t,x,y,a)$. Then $X = B$, and we may naturally define 
\bea
\label{dbY}
\dbY_t := \dbV(t, B_t),\q \dbF(t,x, \dbD) := \big\{f(t,x,y,a): y\in \dbD, a\in A\big\}.
\eea
In the linear case: $f=f(t,x,a)$ and thus $\dbF(t,x) = \big\{f(t,x,a):  a\in A\big\}$ is independent of $\dbD$, the random set-valued process  $\dbY_t:= \dbV(t, B_t)$ indeed satisfies \eqref{setBSDE} in the sense of \cite{AMW}. 

However, when $f$ depends on $y$, the $\dbY, \dbF$ in \eqref{dbY} do not satisfy \eqref{setBSDE}. That is, \eqref{setBSDE} is not the stochastic counterpart of \eqref{HJB}. The reason is the same as in Remark \ref{rem-fundamental} (ii). In \eqref{FBSDE}, the $Y$ in the left side and that in the right side of the BSDE  are required to be the same process. In \eqref{setBSDE}, however, one allows to consider different selectors for the $\dbY$ in the left side and that in the right side of the equation. Consequently, the solutions to \eqref{HJB} and to \eqref{setBSDE} are typically not equal. We shall remark that, the applications mentioned in Introduction typically fall into our framework, although technically many of them are not covered by the current form of our HJB equation \eqref{HJB}.

\begin{appendix}
\section{Some technical proofs}
\label{sect-Appendix}
\setcounter{equation}{0}
\no{\bf Proof of Proposition \ref{prop-geodesic}.}  Again we denote $\br, \bn$ for notational simplicity. We prove it only for $x>x_0$. Fix $x_1>x_0$. Without loss of generality, we assume $\th$ is absolutely continuous in $x\in [x_0,x_1]$ with appropriate derivative function $\th'$, as otherwise the length of $\th$ would be $\infty$. Thus we have
\beaa
\th(x) = y_0 - \int_{x_0}^x \th'(\tilde x) d\tilde x,
\eeaa
Note that $\br (x, \th(x))=0$. Then, for Lebesgue-a.e. $x$,
\beaa
0 &=& {d\over dx} \br (x, \th(x)) = \td_x \br (x, \th(x)) - \td_y \br (x, \th(x)) \cd \th'(x) \\
&=& \td_x \br (x, \th(x)) - \bn(x, \th(x)) \cd \th'(x),
\eeaa
and thus
\beaa
\zeta(x) :=  \th'(x) - \td_x \br\ \bn (x, \th(x)) \in \dbT_\dbV(x, \th(x)).
\eeaa
Therefore,
\beaa
&&\dis L_\th(x_0, x_1) = \int_{x_0}^{x_1} \sqrt{1+|\th'(x)|^2} ~dx = \int_{x_0}^{x_1} \sqrt{1+\big|\td_x \br\ \bn (x, \th(x)) + \zeta(x) \big|^2} ~dx\\
&&= \int_{x_0}^{x_1} \sqrt{1+\big|\td_x \br (x, \th(x))\big|^2 + |\zeta(x)|^2} ~dx \ge  \int_{x_0}^{x_1} \sqrt{1+\big|\td_x \br (x, \th(x))\big|^2} ~dx
\eeaa
This implies that
\beaa
&&\dis \limsup_{x_1\downarrow x_0} {1\over x_1-x_0}\Big[L_\U(x_0, x_1) - L_\th(x_0, x_1)\Big]\\
&&\dis \le \limsup_{x_1\downarrow x_0} {1\over x_1-x_0}\Big[\int_{x_0}^{x_1} \sqrt{1+\big|\td_x \br (x, \U(x))\big|^2} ~dx - \int_{x_0}^{x_1} \sqrt{1+\big|\td_x \br (x, \th(x))\big|^2} ~dx\Big]\\
&&\dis \le \limsup_{x_1\downarrow x_0} {1\over x_1-x_0}\Big[\int_{x_0}^{x_1} \big|\sqrt{1+\big|\td_x \br (x, \U(x))\big|^2} - \sqrt{1+\big|\td_x \br (x_0, y_0)\big|^2}\big| ~dx \\
&&\dis\qq\qq\qq + \int_{x_0}^{x_1} \big|\sqrt{1+\big|\td_x \br (x, \th(x))\big|^2} - \sqrt{1+\big|\td_x \br (x_0, y_0)\big|^2}\big|~dx\Big]=0.
\eeaa

\vspace{-9mm}
\qed

\bs
\no{\bf Proof of Lemma \ref{lem-paxxV}.} Recall \eqref{nD} and consider the natural extension $\wh \bn = \td_y \br $.  By \eqref{patxV}, \eqref{product}, and \eqref{paxfhat} we have, for $i,j=1,\cds, d$,  and $(t,x,y)\in \dbG_\dbV$,
\bea
\label{paxxV1}
\dis \pa_{x_ix_j} \dbV(t,x, y) &=& -\pa_{x_i}\big(\td_{x_j} \br \bn\big)(t,x,y) = -\Big[ \pa_{x_i}(\td_{x_j} \br ) \bn + \td_{x_j} \br \pa_{x_i} \bn\Big](t,x, y)\nonumber\\
&=&- \Big[ \td_{x_ix_j} \br + \td_{x_i}r \td_{x_j y} \br \cd \bn \Big] \bn(t,x, y) - \td_{x_j} \br \pa_{x_i} \bn(t,x, y).
\eea
Recall \eqref{nD}, at  $(t,x,y)\in O_\e(\dbG_\dbV)$ we have
\beaa
\td_{x_j y} \br \cd \td_y \br = {1\over 2} \td_{x_j}\big(|\td_y \br |^2\big) = 0.
\eeaa
In particular, $\td_{x_j y} \br \cd \bn(t,x,y)=0$ for $(t,x,y)\in \dbG_\dbV$. Plugging this into \eqref{paxxV1} we obtain \eqref{paxxV} immediately.

Moreover, again considering the extension $\wh \bn^i = \td_{y_i} \br $, by \eqref{paxfhat}  and \eqref{payfhat} we have
\beaa
\pa_x n^i = \td_{x y_i} \br - \td_x \br (\td_{y_iy}\br \cd\bn),\qq \pa_y \bn^i = (\mathbf{I}_m - \bn\bn^\top)\td_{y_i y}\br.
\eeaa
Similarly, by \eqref{nD} we have $\td_{y_iy}\br \cd\bn =0$, which implies  \eqref{pan} immediately. 
\qed

\ms
\no{\bf Proof of Proposition \ref{prop-convex}.} Under Assumption \ref{assum-coefficients}, clearly $\dbV(t,x)$ is bounded. Then the compactness follows from its closeness. To show the convexity, we recall \eqref{targetV} and introduce
\bea
\label{W}
~~~~~~J(t,x,y,\a, Z) := \dbE\Big[\big|Y^{t,x, y, \a, Z}_T - g(X^{t,x,\a}_T)\big|^2\Big],~ W(t,x,y) := \inf_{\a, Z} J(t,x,y,\a, Z).
\eea
It is clear that $W\ge 0$. Following \cite{KMZ} we have the duality: 
\bea
\label{duality}
\dbV(t,x) = \{y: W(t,x,y)=0\}.
\eea
 Then it suffices to show that $W$ is convex in $y$.\footnote{The convexity of $W$ is a sufficient but not a necessary condition of  the convexity of $\dbV$, since the latter involves only the zero level set of $W$.} Indeed, for any $\ul y, \ol y\in \dbV(t,x)$ and $y := \th \ul y + (1-\th)\ol y$, $\th \in (0,1)$, by \eqref{duality} we have $W(t,x,\ul y) = W(t,x,\ol y)=0$, and then by the convexity of $W$ we have
 \beaa
 0\le W(t,x, y) \le \th W(t,x,\ul y) + (1-\th) W(t,x,\ol y) =0.
 \eeaa
 Thus by \eqref{duality} again we have $y\in \dbV(t,x)$. That is, $\dbV(t,x)$ is convex.
 
 We now prove the convexity of $W$ in two steps. Throughout this proof, $\rho$ denotes a generic modulus of continuity function which may vary from line to line.

\ss
{\it Step 1.} In this step we establish the DPP and the desired regularity of $W$. First, note that $J(t,\cd, \a,Z)$ is uniformly continuous in $(x,y)$, and thus $W$ is uniformly continuous in $(x,y)$. By this regularity, it follows from standard arguments to  have the DPP for deterministic time: 
\bea
\label{convex-DPP}
W(t,x, y) = \inf_{\a,Z} \dbE\Big[ W(t+\d, X^{t,x,\a}_{t+\d},Y^{t,x, y, \a, Z}_{t+\d})\Big].
\eea
By setting $\a\equiv a_0\in A$ and $Z\equiv 0$ in \eqref{convex-DPP}, we have, 
\bea
\label{DPPle}
W(t,x,y) \le \dbE\Big[ W(t+\d, X^{t,x,a_0}_{t+\d},Y^{t,x, y, a_0, 0}_{t+\d})\Big].
\eea
This implies immediately that $W(t,x,y) \le \liminf_{\d\downarrow 0} W(t+\d, x, y)$. Moreover, for any $\a, Z$, by the tower properties of SDEs and conditional expectations, it is clear that
\beaa
J(t,x,y,\a, Z) \ge \dbE\Big[ W(t+\d, X^{t,x,\a}_{t+\d},Y^{t,x, y, \a, Z}_{t+\d})\Big].
\eeaa
Send $\d\to 0$ we have $J(t,x,y,\a, Z) \ge   \limsup_{\d\downarrow 0} W(t+\d, x, y)$. Then $W(t,x,y) \ge  \limsup_{\d\downarrow 0} W(t+\d, x, y)$. Thus $W$ is right continuous in $t$:
\bea
\label{Wrightt}
 \lim_{\d\downarrow 0} W(t+\d, x, y) = W(t,x,y).
\eea

We note that, by \eqref{convex-DPP},  similarly to \eqref{DPPle} we can show that $\limsup_{\d\downarrow 0} W(t-\d, x, y) \le W(t,x,y)$. So $W$ is upper semicontinuous in $t$.  
However, in general it is not clear whether or not $W$ is left continuous in $t$.

\ms

{\it Step 2.} We now prove the convexity of W and we use two distinct
  approaches: the first one based on the viscosity solution theory and the
  second one on probabilistic arguments. Both approaches have
  independent interest. 

{\it Method 1.}
By DPP \eqref{convex-DPP} and the upper semicontinuity of $W$ in Step 1, it is clear that $W$ is a viscosity subsolution to the following standard HJB equation with unbounded control $z$:
\bea
\label{WHJB}
\left.\ba{c}
\dis \pa_t W(t,x,y) + \inf_{a,z} \Big[\tr\Big(\tfrac{1}{2} \pa_{xx} W \si\si^\top(t,x,a) + \tfrac{1}{2} \pa_{yy} W  zz^\top + \pa_{xy} W\si(t,x,a) z^\top\Big) \ss\\
\dis + \pa_x W \cd b(t,x,a) - \pa_y W\cd  f(t,x, y, z, a)\Big]\ge 0.
\ea\right.
\eea
Here, by abusing notations we use $\pa_t, \pa_x$ etc to denote standard partial derivatives. For any $(t,x,y)$ and any smooth test function $\f$ at $(t,x,y)$, set $a:= a_0$ for some fixed $a_0\in A$ and $z := n [e_1, \cds, e_1]$, where $e_1:= (1, 0, \cds, 0)^\top$. Then, by sending $n\to \infty$, it follows from the Lipschitz continuity of $f$ in $z$ that 
\beaa
\pa_{y_1y_1} \f(t,x,y) \ge 0.
\eeaa
Now fix $(t,x, y)$ and introduce the function $U(\tilde y_1) := W(t,x, (\tilde y_1, y_2, \cds, y_m)^\top)$, $\tilde y_1 \in \dbR$. Then $U$ is uniformly continuous, and it follows from   \cite[Proposition 5.10]{Touzi-Lecture} that $U$ is a viscosity subsolution of $U''(\tilde y_1)\ge 0$. Applying  \cite[Lemma 5.25]{Touzi-Lecture}, we obtain that $U$ is convex in $\tilde y_1$, which implies that $W$ is convex in $y_1$.

Moreover, for any orthogonal matrix $Q\in \dbR^{m\times m}$, introduce $\tilde W(t,x,y) := W(t,x,Qy)$. It is clear that $\tilde W$ satisfies a PDE similar to \eqref{WHJB}, with a quadratic term $\tfrac{1}{2} \tr\big(\pa_{yy} \tilde W (Qz)(Qz)^\top\big)$. 
Then, by choosing $z = n Q^{-1}[e_1,\cds, e_1]$, one can show as above that $\tilde W$ is convex in $y_1$. 
Now, for any $\ul y, \ol y\in\dbR^m$ and $y = \th \ul y + (1-\th) \ol y$, $\th\in[0,1]$, set $Q$ such that $Q^{-1}(\ol y - \ul y) = |\ol y - \ul y|e_1$. Then, $Q^{-1} y = Q^{-1}\ol y - \th Q^{-1} (\ol y-\ul y) = Q^{-1}\ol y - \th |\ol y-\ul y| e_1$, namely $\th$ affects only the first component. Therefore, by the claimed convexity of $\tilde W$, we have  
\beaa
W(t,x,y) &=& \tilde W(t,x, Q^{-1}y) \leq \th W(t,x,Q^{-1}\ul y) + (1-\th)W(t,x,Q^{-1}\ol y) \\
&=& \th W(t,x,\ul y) + (1-\th)W(t,x,\ol y)
\eeaa
This is the desired convexity.

\ms {\it Method 2.}
Fix $\ul y \neq \ol y\in \dbR^m$ and $y = \th \ul y + (1-\th) \ol y$ for some $\th\in (0,1)$. For any $\d>0$, $\a$, and $n$, set $\eta := {\ol y-\ul y\over |\ol y-\ul y|}$, $X := X^{0,x_0, \a}$, $Y^{n}:= Y^{0,x_0, y, \a, n[\eta,\cds,\eta]}$, and 
\bea
\label{convex-taun}
\left.\ba{c}
\dis\t^n_\d := \inf\big\{t>0: \eta\cd Y^n_t = \eta\cd \ul y ~\mbox{or}~ \eta\cd \ol y\big\}\wedge \d,\ss\\
\dis Z^{n,\d}_t := n[\eta,\cds, \eta] \1_{[0, \t^n_\d]}(t),\q Y^{n,\d} := Y^{0,x_0, y, \a, Z^{n,\d}}.
\ea\right.
\eea
By DPP \eqref{convex-DPP} we have, for some appropriate modulus of continuity function $\rho$,
\bea
\label{convex-Ei}
&&W(0,x_0, y) \le  \dbE\Big[ W(\d, X_{\d},Y^{n,\d}_{\d})\Big] \le \dbE\Big[ W(\d, x_0,Y^{n,\d}_{\d})\Big] + \rho(\d)\nonumber\\
&&\qq\qq\q =\dbE\Big[ W(\d, x_0,Y^{n,\d}_{\d}) \big[ \1_{E^{n,\d}_1} + \1_{E^{n,\d}_2} + \1_{\{\t^n_\d =\d\}}\big]\Big] + \rho(\d),\\
&&\mbox{where}\q E^{n,\d}_1:= \{\t^n_\d <\d, \eta\cd Y^n_{\t^n_\d}=\eta\cd \ul y\},~~ E^{n,\d}_2:= \{\t^n_\d <\d, \eta\cd Y^n_{\t^n_\d}=\eta\cd \ol y\}.\nonumber
\eea
Here the second inequality used the regularity of $W$ in $x$ and the following obvious estimate:
\beaa
\dbE[|X_\d-x|^2] \le C\d.
\eeaa 

We next prove a crucial estimate: 
\bea
\label{convexest1}
|Y^{n,\d}_t| \le C(1+n\t^n_\d), \q \mbox{and}\q |f(t, X_t, Y^{n,\d}_t, Z^{n,\d}_t,\a_s)|\le Cn,\q t\le \t^n_\d.
\eea
Indeed, by \eqref{convex-taun} $\eta\cd Y^{n,\d}_t$ is bounded. Denote $\tilde Y^{n,\d}:= (\eta\cd Y^{n,\d})\eta = \eta \eta^\top  Y^{n,\d}$ and $\D Y^{n,\d} := Y^{n,\d}-\tilde Y^{n,\d}$. Note that $\eta \eta^\top Z^{n,\d} = Z^{n,\d}$. Then, for $0\le t\le \t^n_\d$,
\beaa
&&\tilde Y^{n,\d}_t = \eta\eta^\top y+ \int_0^t \eta\eta^\top f ds + \int_0^t Z^{n,\d}_s dB_s;\\
&&\D Y^{n,\d}_t = (\mathbf{I}_m-\eta \eta^\top)y +  \int_0^t (\mathbf{I}_m-\eta\eta^\top) f(s, X_s, \tilde Y^{n,\d}_s + \D Y^{n,\d}_s, Z^{n,\d}_s,\a_s) ds. 
\eeaa
Since $|\tilde Y^{n,\d}|\le C$ and $|Z^{n,\d}|\le n$, it follows from the Lipschitz property of $f$ that $|\D Y^{n,\d}_t| \le C(1+nt)$. This, together with the boundedness of $\tilde Y^{n,\d}$ again, implies $|Y^{n,\d}_t| \le C(1+nt)\le C(1+n\t^n_\d)$. Moreover, by Assumption \ref{assum-coefficients} and $|Z^{n,\d}|\le n$, we obtain \eqref{convexest1}.

We now estimate the terms in \eqref{convex-Ei}. We first consider
$E^{n,\d}_1$. In this case, note that
$\eta\cd (Y^n_{\t^n_\d}-\ul y)=0$. Moreover,
$(\mathbf{I}_m-\eta\eta^\top)(Y_0^n - \ul y) = (\mathbf{I}_m-\eta\eta^\top)(y - \ul y) =
(1-\th)|\ol y - \ul y|(\mathbf{I}_m-\eta\eta^\top)\eta = 0$ and
$(\mathbf{I}_m - \eta\eta^\top)Z_t^{n,\d}= 0, 0\leq t \leq \t_\d^n$. Then, by  \eqref{convexest1} we have
\bea
\label{convexest2}
|Y^n_{\t^n_\d}-\ul y| = \big|(\mathbf{I}_m-\eta\eta^\top)(Y^n_{\t^n_\d}-\ul y) \big|= \Big| \int_0^{\t^n_\d} (\mathbf{I}_m-\eta\eta^\top) \cd f ds\Big| \le Cn \t^n_\d.
\eea
Note further that $Z^{n,\d}_t =0$, $\t^n_\d\le t<\d$, then by Assumption \ref{assum-coefficients}, \eqref{targetV}  and \eqref{convexest1} we have 
\beaa
&&\dis |Y^{n,\d}_t|\le C(|Y^{n,\d}_{\t^n_\d}| + t-\t^n_\d) \le C(1 + \d + n\t^n_\d),\q \t^n_\d \le t\le \d;\\
&&\dis  \big|Y^{n,\d}_\d - Y^{n,\d}_{\t^n_\d}\big| \le \int_{\t^n_\d}^\d |f(s, X_s, Y^{n,\d}_s, 0, \a_s)| ds \le C\d(1+n\t^n_\d).
\eeaa
This, together with \eqref{convexest2}, implies
\bea
\label{convexest3}
\big|Y^{n,\d}_\d - \ul y\big| \le \big|Y^{n,\d}_\d - Y^{n,\d}_{\t^n_\d}\big| + |Y^n_{\t^n_\d}-\ul y| \le Cn \t^n_\d + C\d,\q \mbox{a.s. on}~ E^{n,\d}_1.
\eea
Similarly, we have
\bea
\label{convexest4}
\big|Y^{n,\d}_\d - \ol y\big| \le Cn \t^n_\d + C\d,\q \mbox{a.s. on}~ E^{n,\d}_2.
\eea
Plug \eqref{convexest3}, \eqref{convexest4}, and the first estimate of \eqref{convexest1} into \eqref{convex-Ei}, we obtain 
\bea
\label{convex-Ei2}
W(0,x_0, y_0) &\le& W(\d, x_0,\ul y) \dbP(E^{n,\d}_1) +  W(\d, x_0,\ol y) \dbP(E^{n,\d}_2) \\
&&+  C\dbE[n\t_\d^n] + C\dbP(n\t_\d^n > \d) + C(1+n\d) \dbP(\t^n_\d =\d) +   \rho(\d).\nonumber
\eea
Here we used the fact $\rho(n\t^n_\d) \le \rho(\d) + C(\1_{\{n\t^n_\d >\d\}} + n\t^n_\d)$, by assuming its linear growth.

We next estimate $\t^n_\d$. Note that
\bea
\label{convex-barB}
 &&\int_0^t Z_s^n dB_s =  n\sqrt{d} \bar B_t \eta,\q\mbox{where}\q \bar B_t := {1\over \sqrt{d}}\sum_{i=1}^d B^i_t;\nonumber\\
 &&\eta\cd Y^n_t = \eta \cd y + \int_0^t \eta \cd f ds +  n\sqrt{d} \bar B_t =  \eta \cd y +   n\sqrt{d} \bar B^n_t,\\
 &&\dis \mbox{where}\q \bar B^n_t:=\bar B_t - \int_0^t \th^n_sds,\q \th^n_t := -{1\over n\sqrt{d}} \eta\cdot f(t, X_t, Y^n_t, Z^n_t, \a_s).\nonumber
  \eea
Note that $Y^n_t = Y^{n,\d}_t$ for $t\le \t^n_\d$, by \eqref{convexest1} $\th^n_t$ is bounded for $t\le \t^n_\d$. Note further that $\bar B$ is a Brownian motion, introduce 
\beaa
d\dbP^n := M^n_t d\dbP,\q M^n_t := e^{\int_0^t \th^n_s d\bar B_s - {1\over 2} \int_0^t |\th^n_s|^2ds}.
\eeaa
Then $\bar B^n$ is a $\dbP^n$-Brownian motion. By the middle line of \eqref{convex-barB} and by rescaling the time for the Brownian motion it is straightforward to check that 
\bea
\label{convexest5}
\dbE^{\dbP^n}[\t^n_\d] \le {C\over n^2}.
\eea
\sloppy Notice again that $\th^n$ is bounded for $t\le \t_\d^n$, then by standard estimates we have $\dbE^{\dbP^n}\big[(M^n_{\t^n_\d})^{-4}\big]\le C$. 
 This, together with \eqref{convexest5}, implies
 \begin{equation}
   \label{convexest6}
   \dbE[\t^n_\d] \le \dbE[(\t^n_\d)^{3\over 4}] =
   \dbE^{\dbP^n}\big[(M^n_{\t^n_\d})^{-1}(\t^n_\d)^{3\over 4}\big] \le
   \Big(\dbE^{\dbP^n}\big[(M^n_{\t^n_\d})^{-4}\big]\Big)^{1\over 4}
   \Big(\dbE^{\dbP^n}[\t^n_\d]\Big)^{3\over 4}  \le {C\over n^{3\over 2}}.
   \hspace{-1em}
 \end{equation}
In particular, this implies
\bea
\label{convexest7}
n\dbE[\t^n_\d] \le {C\over \sqrt{n}},\q \dbP(\t^n_\d=\d) \le {1\over \d} \dbE[\t^n_\d] \le {C\over \d n^{3\over 2}},\q \text{and}\q \dbP(n\t_\d^n > \d) \leq \frac{C}{\d n^{\frac{1}{2}}.\hspace{-1em}}
\eea
Then, by fixing $\d$ and sending $n\to \infty$ in \eqref{convex-Ei2}, we have
\begin{equation}  
  \label{convex-Ei3}
  W(0,x_0, y_0) \le W(\d, x_0,\ul y) \limsup_{n\to \infty}\dbP(E^{n,\d}_1) +
  W(\d, x_0,\ol y) \limsup_{n\to \infty}\dbP(E^{n,\d}_2) + \rho(\d).
  \hspace{-1em}
\end{equation}

Moreover, note that
\beaa
&&\dis \th \eta\cd \ul y + (1-\th) \eta\cd \ol y = \eta\cd y = \eta\cd \dbE\Big[Y^n_{\t^n_\d} + \int_0^{\t^n_\d}  f ds\Big] \\
&&\dis = \eta\cd \ul y \dbP(E^{n,\d}_1) + \eta\cd \ol y \dbP(E^{n,\d}_2) +\eta\cd \dbE\Big[Y^n_{\t^n_\d}\1_{\{\t^n_\d=\d\}} + \int_0^{\t^n_\d} f ds\Big].
\eeaa
By the same analysis we have
\beaa
 \lim_{n\to \infty} \eta\cd \dbE\Big[Y^n_{\t^n_\d}\1_{\{\t^n_\d=\d\}} + \int_0^{\t^n_\d} f ds\Big]=0,\q \lim_{n\to\infty}\Big[\dbP(E^{n,\d}_1)+\dbP(E^{n,\d}_2)\Big] =1.
 \eeaa
  Then we must we have
\beaa
\lim_{n\to\infty}| \dbP(E^{n,\d}_1) - \th| = \lim_{n\to\infty}|\dbP(E^{n,\d}_2) - (1-\th)| =0.
 \eeaa
Plug this into \eqref{convex-Ei3}, we obtain
\beaa
 W(0,x_0, y_0) \le  \th W(\d, x_0, \ul y)   + (1-\th)W(\d, x_0, \ol y)  +C\d+\rho(\d)
 \eeaa
 Now send $\d\to 0$, by the right continuity $W$ in $t$ we obtain the desired convexity.
\qed

\ms
\no{\bf Proof of Example \ref{eg-nonconvex}.} We first prove \eqref{V-quadratic}. Denote
\beaa
\tilde f_1(a, y) := a_1,\q \tilde f_2(a, y) := a_2,\q \tilde Y^\a_t = \int_t^T \tilde f(\a_s, \tilde Y^\a_s)ds - \int_t^T \tilde Z_s^\a dB_s.
\eeaa
Then one can easily check that
\bea
\label{tildeVb}
\left.\ba{c}
\dis \tilde V(t) := \big\{\tilde Y^\a_t: \a\in \cA_t\big\} = \big\{\tilde y\in \dbR^2: |\tilde y|\le T-t\big\},\\
\dis\tilde \dbV_b(t):= \big\{\tilde y\in \dbR^2: |\tilde y|=T-t\big\} = \big\{(T-t)\big(\cos \th, \sin\th\big)^\top: \th \in [0, 2\pi)\big\}.
\ea\right.
\eea
Consider a function  $\psi: \dbR^2\to \dbR^2$ and set
\begin{equation*}
  \begin{aligned}
    &Y^\a_t := \psi(\tilde Y^\a_t),\ \
    Z_t^\a := \big[\td\psi_1(\tilde Y_t^\a),  \td\psi_2(\tilde Y_t^\a)\big]^\top \tilde Z_t^\a,
    \\&\mbox{where}\q \psi_1(\tilde y) := \tilde y_1,\q
    \psi_2(\tilde y) := [1+|\tilde y_1|^2] \tilde y_2.   
  \end{aligned}
\end{equation*}
Then one may check straightforwardly that
\begin{equation*}
  \begin{aligned}
    d Y^{\a, 1}_t &= d \tilde Y^{\a,1}_t
    = -\a^1_t dt + \tilde Z_s^{\a,1}dB_t
    = -f_1(\a_t, Y^\a_t, Z_t^\a)dt + Z_s^{\a,1}dB_t;
    \\d Y_t^{\a,2} &= 2 \tilde Y^{\a,1}_t\tilde Y^{\a,2}_t  d \tilde Y^{\a,1}_t
    + [1+ |\tilde Y^{\a,1}_t|^2] d \tilde Y^{\a,2}_t
    + \big[\tilde Y_t^{\a,2}|\tilde Z_t^{\a,1}|^2
    + 2 \tilde Y_t^{\a,1} \tilde Z_t^{\a,1} \tilde Z_t^{\a,2}\big]dt
    \\&= -\Big[2 \tilde Y^{\a,1}_t\tilde Y^{\a,2}_t  \a^1_t
    +  [1+ |\tilde Y^{\a,1}_t|^2] \a^2_t -
    \tilde Y_t^{\a,2}|\tilde Z_t^{\a,1}|^2
    - 2 \tilde Y_t^{\a,1} \tilde Z_t^{\a,1} \tilde Z_t^{\a,2}\Big]dt
    \\&\q + \Big[2\tilde Y_t^{\a,1}\tilde Y_t^{\a,2} \tilde Z_t^{\a,1}
    + [1+|\tilde Y_t^{\a,1}|^2]\tilde Z_t^{\a,2}\Big]dB_t
    \\&= -\Big[{2 Y^{\a,1}_tY^{\a,2}_t\over 1+|Y^{\a,1}_t|^2}  \a^1_t
    +  [1+ |Y^{\a,1}_t|^2] \a^2_t
    \\&\hspace{2.5em}
    - \frac{Y_t^{\a,2}|Z_t^{\a,1}|^2+ 2Y_t^{\a,1}Z_t^{\a,1}Z_t^{\a,2}}{1+|Y_t^{\a,1}|^2}
    + \frac{4|Y_t^{\a,1}|^2Y_t^{\a,2}|Z_t^{\a,1}|^2}{(1+|Y_t^{\a,1}|^2)^2}
    \Big]dt + Z_t^{\a,2} dB_t
    \\&= -f_2(\a_t,Y_t^\a,Z_t^\a) dt + Z_t^{\a,2}dB_t
  \end{aligned}
\end{equation*}
That is, $Y^\a$ satisfies \eqref{FBSDE} for given $f$. Therefore, $\dbV(t) = \big\{ \psi(\tilde y): \tilde y\in \tilde \dbV(t)\big\}$. Note further that $\psi_1(\tilde y_1)=\tilde y_1$, and $\psi_2$ is strictly increasing in $\tilde y_2$. It is clear that
\beaa
\dbV_b(t) = \big\{\psi(\tilde y): \tilde y\in \tilde \dbV_b(t)\big\}.
\eeaa
Plug \eqref{tildeVb} into it, we obtain \eqref{V-quadratic} immediately.

We next analyze the convexity of $\dbV(t)$. Assume for simplicity that $t=0$. Note that
\beaa
&\dis \tilde \dbV_b(0) = \Big\{ \big(y_1, \sqrt{T^2-|y_1|^2}\big), ~ \big(y_1, -\sqrt{T^2-|y_1|^2}\big): |y_1|\le T\Big\};\\
&\dis \dbV_b(0) = \Big\{ \big(y_1, \f(y_1)\big), \big(y_1, -\f(y_1)\big): |y_1|\le T\Big\},~\mbox{where}~ \f(y_1) := [1+|y_1|^2]\sqrt{T^2-|y_1|^2}.
\eeaa
One may compute straightforwardly that: 
\beaa
\f''(y_1) = { 6 |y_1|^4 - 9 T^2|y_1|^2 + 2T^4 -T^2 \over (T^2-|y_1|^2)^{3\over 2}},\q |y_1|<T.
\eeaa
Note that
\beaa
\sup_{|y_1|<T} \big[6 |y_1|^4 - 9 T^2|y_1|^2\big] = 0.
\eeaa
So when $T\le {1\over \sqrt{2}}$ and thus $2T^4 - T^2\le 0$, we have $\f''(y_1) \le 0$ for $|y_1|<T$, and  in this case $\dbV(0)$ is indeed convex. However, when $T> {1\over \sqrt{2}}$, we find that $\f''(y_1) <0$ for $|y_1|\approx T$, but $\f"(0) = {2T^4-T^2\over T^3}> 0$, then $\dbV(0)$ is nonconvex.
\qed

\ms

\no{\bf Proof of Proposition \ref{prop-reg}.}  Given the convexity of $\dbV(t,x)$ in Proposition \ref{prop-convex},  Part (ii) is a direct consequence of Part (i) and Lemma \ref{lem-dD} below. So we shall only prove (i) here. We remark that this regularity does not rely on the convexity of $\dbV(t,x)$. Throughout the proof, $\rho$ denotes a generic modulus of continuity function, which may vary from line to line.  We proceed in four steps. 

{\bf Step 1.} We first prove the regularity in $x$. For any $y = Y^{t,x,\a}_t\in \dbV(t,x)$, denote $\tilde y := Y^{t,\tilde x,\a}_t\in \dbV(t, \tilde x)$. By standard SDE and BSDE estimates we have
\beaa
d(y, \dbV(t,\tilde x)) \le |y-\tilde y| = |Y^{t,x,\a}_t - Y^{t,\tilde x,\a}_t| \le \rho(|x-\tilde x|).
\eeaa
Since $d$ is Lipschitz continuous in $y$ with Lipschitz constant $y$, this clearly implies $\sup_{y\in \dbV(t,x)} d(y, \dbV(t,\tilde x)) \le \rho(|x-\tilde x|)$. Similarly we have $\sup_{\tilde y\in \dbV(t,\tilde x)}d(\tilde y, \dbV(t,x)) \le \rho(|x-\tilde x|)$. Then we obtain $d(\dbV(t,x), \dbV(t, \tilde x)) \le \rho(|x-\tilde x|)$. 

{\bf Step 2.} In this step we estimate $\sup_{\tilde y\in\dbV(t+\d,x)}d (\tilde y, \dbV(t,x))$.   For $\tilde y = Y^{t+\d,x,\tilde \a}_{t+\d}\in \dbV(t+\d,x)$ with $\tilde \a\in \cA_{t+\d}$, fix some $a_0\in A$ and consider $\a\in \cA_t$ defined by: $\a_s(\o) := a_0 \1_{[t, t+\d)}(s) + \tilde \a_{s}(\o^{t+\d}) \1_{[t+\d, T]}(s)$ where $\o_s^{t+\d} := \o_s-\o_{t+\d} \in \O_{t+\d}$. Then $y := Y^{t,x,\a}_t \in \dbV(t,x)$. Note that
\beaa
\Big|Y^{t,x,\a}_{t+\d} - \tilde y\Big|= \Big|Y^{t+\d, X^{t,x,a_0}_{t+\d}, \tilde \a}_{t+\d} - Y^{t+\d,x,\tilde \a}_{t+\d}\Big| \le \rho\big(X^{t,x,a_0}_{t+\d} - x\big).
\eeaa
Recall \eqref{FBSDE} for $Y^{t,x, \a}$ and introduce another BSDE on $[t, t+\d]$:
\beaa
&&Y^{t,x,\a}_s = Y^{t,x,\a}_{t+\d} + \int_s^{t+\d} f(r, X^{t,x,a_0}_r, Y^{t,x,\a}_r, Z^{t,x,\a}_r, a_0) dr - \int_s^{t+\d} Z^{t,x,\a}_r dB_r;\\
&&\tilde Y_s = \tilde y+ \int_s^{t+\d} f(r, X^{t,x,a_0}_r, \tilde Y_r, \tilde Z_r, a_0) dr - \int_s^{t+\d} \tilde Z_r dB_r.
\eeaa
By standard BSDE estimates we get, for possibly different modulus of continuity function $\rho$, 
\beaa
|y - \tilde Y_t|^2= |Y^{t,x,\a}_t- \tilde Y_t|^2 \le \dbE\Big[\big|Y^{t,x,\a}_{t+\d} - \tilde y\big|^2\Big] \le \dbE\Big[\rho^2\big(X^{t,x,a_0}_{t+\d} - x\big)\Big] \le \rho(\d),
\eeaa
Moreover, since the FBSDE system for $(X^{t,x,a_0}, \tilde Y, \tilde Z)$ is Markovian and $\tilde y\in \dbV(t+\d, x)$ is bounded, under Assumption \ref{assum-coefficients} one can easily show that $(\tilde Y, \tilde Z)$ is bounded.  This implies that
$
\big|\tilde Y_t - \tilde y\big| \le C\sqrt{\d}.
$
Then $|y- \tilde y|\le C\sqrt{\d} + \rho(\d)\le \rho(\d)$, and thus we can easily get:
\bea
\label{Vcontt1}
\sup_{\tilde y\in \dbV(t+\d, x)} d(\tilde y, \dbV(t,x)) \le  \rho(\d).
\eea

{\bf Step 3.} To see the opposite direction, we first show that $d(g(x), \dbV(T-\d, x))\le \rho(\d)$. Given arbitrary $y = Y^{T-\d,x,\a}_{T-\d} \in \dbV(T-\d, x)$ with $\a\in \cA_{T-\d}$. Denote $\f^\a:= \f^{T-\d,x,\a}$ for $\f:= X, Y, Z$. Under Assumption \ref{assum-coefficients}, by standard SDE and BSDE estimates we have
\beaa
\dbE\big[|X^{\a}_T - x|^2\big] \le C\d,\q \dbE\Big[\sup_{T-\d\le s\le T} |Y^{\a}_s|^2 + \int_{T-\d}^T |Z^\a_s|^2 ds\Big] \le C.
\eeaa
Then, for a generic modulus of continuity function $\rho$, 
\beaa
&&\dis |y- g(x)| = \Big|\dbE\big[g(X^\a_T) + \int_{T-\d}^T f(t, X^\a_s, Y^\a_s, Z^\a_s, \a_s)ds - g(x)\big]\Big|\\
&&\dis \le \dbE\Big[\rho(|X^\a_T-x|) + C\int_{T-\d}^T [1+|Y^\a_s| + |Z^\a_s|]ds \Big]\\
&&\dis \le \rho(\d) + C\Big(\dbE\Big[\d \int_{T-\d}^T [1+|Y^\a_s|^2 + |Z^\a_s|^2]ds \Big]\Big)^{1\over 2} \le \rho(\d) + C\sqrt{\d}\le \rho(\d).
\eeaa
This implies that $d(g(x), \dbV(T-\d, x))\le \rho(\d)$.

{\bf Step 4.} We now show that $\sup_{y\in\dbV(t,x)}d (y, \dbV(t+\d,x))\le \rho(\d)$. This, together with the estimate in Step 2, implies the desired regularity in $t$. Fix $y= Y^{t,x,\a}_t \in \dbV(t,x)$ for some $\a\in \cA_t$. As in Step 3 we denote $\f^\a:= \f^{t,x,\a}$ for $\f=X, Y, Z$. Abusing the notation $(\tilde Y, \tilde Z)$ with Step 2, we consider the following BSDE on $[t, T-\d]$:
\bea
\label{Vreg-tildeY}
\tilde Y_s = g(X^\a_{T-\d}) + \int_t^{T-\d} f(r, X^\a_r, \tilde Y_r, \tilde Z_r, \a_r) dr - \int_t^{T-\d}\tilde Z_r dB_r.
\eea
Note that $(Y^\a, Z^\a)$ satisfies the BSDE on $[t, T-\d]$ with terminal condition $Y^\a_{T-\d}$. One can easily see that, as in \eqref{DPP-claim} with deterministic $\t$, $Y^\a_{T-\d}\in \dbV(T-\d, X^\a_{T-\d})$. Then, by Step 3, $\big|Y^\a_{T-\d} -  g(X^\a_{T-\d})\big| \le \rho(\d)$. Then, by comparing the two BSDEs, we obtain immediately that
\beaa
|\tilde Y_t - y| = |\tilde Y_t - Y^\a_t| \le \rho(\d).
\eeaa

We next introduce another shifted FBSDE, still on $[t, T-\d]$:
\bea
\label{Vreg-barY}
\left.\ba{lll}
\dis\bar X_s := x + \int_t^s b(r+\d, \bar X_r, \a_r) dr + \int_t^s \si(r+\d, \bar X_r, \a_r) dB_r;\\
\dis \bar Y_s = g(X^\a_{T-\d}) + \int_s^{T-\d} f(r+\d, X^\a_r, \bar Y_r, \bar Z_r, \a_r) dr - \int_r^{T-\d}\bar Z_r dB_r.
\ea\right.
\eea
By the continuity of $b,\si, f$ in $t$, it is clear that $|\bar Y_t - \tilde Y_t|\le \rho(\d)$.

Finally, note that $\a\in \cA_t$ means that $\a_s = \f(s, B^t_{[t,s]})$, $s\in [t, T]$, for some appropriate (deterministic) mapping $\f$. Introduce $\bar \a_s := \f(s-\d, B^{t+\d}_{[t+\d, s]})$, $s\in [t+\d, T]$. It is clear that $(\bar \a_s, B^{t+\d}_s)_{s\in [t+\d, T]}$ has the same distribution as $(\a_s, B^t_s)_{s\in [t, T-\d]}$. Then, compare \eqref{Vreg-barY} with the FBSDE for $(X^{\bar \a}, Y^{\bar \a}, Z^{\bar \a})$, we see that $(X_s^{\bar \a}, Y_s^{\bar \a}, Z_s^{\bar \a})_{s\in [t+\d, T]}$ has the same distribution as $(\bar X_s, \bar Y_s, \bar Z_s)_{s\in [t, T-\d]}$. In particular, this means that $Y^{\bar \a}_{t+\d} = \bar Y_t$. Therefore,
\beaa
\big|y- Y^{\bar \a}_{t+\d}\big| = \big|y- \bar Y_t\big|  \le \big|y-\tilde Y_t\big| + \big|\tilde Y_t - \bar Y_t\big|\le \rho(\d).
\eeaa
Since $\bar \a\in \cA_{t+\d}$, this implies that $d(y, \dbV(t+\d,x)) \le \rho(\d)$, and the proof is complete.
\qed

\begin{lem}
\label{lem-dD}
Assume $\dbD$ and $\tilde \dbD$ are compact and convex. Then 
\bea
\label{Dest}
d(\dbD_b, \tilde \dbD_b) = d(\dbD, \tilde \dbD),\q\mbox{and}\q \sup_{y\in \dbR^m} \big|\br_{\dbD}(y) - \br_{\tilde \dbD}(y)\big| \le 3d(\dbD, \tilde \dbD).
\eea
\end{lem}
\proof (i) The result $d(\dbD_b, \tilde \dbD_b) = d(\dbD, \tilde \dbD)$ is not new, see e.g. \cite{W2007}. However, for the
  convenience of the readers, we present here an elementary proof in two
  steps.

{\bf Step 1.} In this step we show that $\sup_{y\in \dbD} d(y, \tilde \dbD) \le d(\dbD_b, \tilde \dbD_b)$. Then similarly we have $\sup_{\tilde y\in \dbD} d(\tilde y, \dbD)\le d(\dbD_b, \tilde \dbD_b)$ and thus $d(\dbD, \tilde \dbD)\le d(\dbD_b, \tilde \dbD_b)$. Note that the inequality is trivial when $\dbD\subset\tilde\dbD$, thus at below we assume $\dbD\backslash \tilde \dbD\neq\emptyset$, and then,
\beaa
\sup_{y\in \dbD} d(y, \tilde \dbD) = \sup_{y\in \dbD\backslash \tilde \dbD} d(y, \tilde \dbD) = \sup_{y\in \dbD\backslash \tilde \dbD} d(y, \tilde \dbD_b). 
\eeaa
For each $y\in \dbD\backslash \tilde \dbD$, by compactness and convexity of $\tilde\dbD$,  there exists unique $\tilde y\in \tilde \dbD_b$ such that $d(y, \tilde \dbD_b) = |y-\tilde y|$. Note that, by convexity of $\tilde \dbD$, in this case hyperplane tangent to $y-\tilde y$ at $\tilde y$ is a supporting hyperplane for $\tilde \dbD$. Since $y\in \dbD$, there exists $\bar y\in \dbD_b$ such that $\bar y - \tilde y = c (y - \tilde y)$ for some $c\ge1$. By the convexity of $\tilde \dbD$, we see that $|\bar y - \tilde y| = d(\bar y, \tilde \dbD_b)$. Then, for any $y\in \dbD\backslash \tilde \dbD$,
\beaa
d(y, \tilde \dbD) = |y-\tilde y| \le |\bar y - \tilde y| = d(\bar y, \tilde \dbD_b) \le d(\dbD_b, \tilde \dbD_b).
\eeaa
This implies the desired estimate: $\sup_{y\in \dbD} d(y, \tilde \dbD) \le d(\dbD_b, \tilde \dbD_b)$.

{\bf Step 2.} In this step we prove the opposite direction: $d(y, \tilde \dbD_b) \le d(\dbD, \tilde \dbD)$ for all $y\in \dbD_b$. If $y \notin \tilde \dbD$, then clearly $d(y, \tilde \dbD_b) = d(y, \tilde \dbD) \le d(\dbD, \tilde \dbD)$. So we assume $y \in \tilde \dbD$. Since $y\in \dbD_b$, let $\eta$ denote a unit outward vector corresponding to a supporting hyperplane of $\dbD$ at $y$. Since $y\in \tilde \dbD$, there exists $\tilde y\in \tilde \dbD_b$ such that $\tilde y - y = c \eta$ for some $c>0$. By the definition of $\eta$, this implies that $d(\tilde y, \dbD)= |\tilde y - y|$. On the other hand, since $\tilde y\in \tilde \dbD_b$, by definition $d(y, \tilde \dbD_b)\le |y-\tilde y|$. So
\beaa
d(y, \tilde \dbD_b) \le |\tilde y - y| = d(\tilde y, \dbD) \le d(\dbD, \tilde \dbD),\q\forall y\in \dbD_b.
\eeaa
Similarly we have $d(\tilde y, \dbD_b) \le d(\dbD, \tilde \dbD)$ for all $\tilde y\in \tilde \dbD_b$. Thus $d(\dbD_b, \tilde \dbD_b) \le d(\dbD, \tilde \dbD)$.

(ii) We now estimate $\br_\dbD - \br_{\tilde \dbD}$. Since $\dbD$ and $\tilde \dbD$ are symmetric, it suffices to prove the following estimate and we proceed in four cases:
\bea
\label{rVcont0}
\br_{\tilde \dbD}(y) -  \br_{\dbD}(y) \le 3d(\dbD, \tilde \dbD).
\eea

{\it Case 1.} $y \in \tilde\dbD \backslash \dbD$. Then $\br_{\tilde \dbD}(y) \le 0< \br_{\dbD}(y)$, and \eqref{rVcont} is obvious.

{\it Case 2.}  $y \notin \dbD\cup \tilde \dbD$. Then $\br_{\tilde \dbD}(y) -  \br_{\dbD}(y)  = d(y, \tilde\dbD_b) - d(y, \dbD_b)$. By the compactness there exists $y'\in \dbD_b$ such that $\br_{\dbD}(y)  = d(y, \dbD_b) = |y-y'|$. Then, 
\beaa
&&\br_{\tilde\dbD}(y)= d(y,\tilde\dbV_b)\le |y-y'| + d(y', \tilde\dbD_b)\le \br_{\dbD}(y) + d(\dbD_b, \tilde \dbD_b)=\br_{\dbD}(y) + d(\dbD, \tilde \dbD).
\eeaa

{\it Case 3.} $y \in \dbD\cap \tilde \dbD$. Then $\br_{\tilde \dbD}(y) -  \br_{\dbD}(y)  = d(y, \dbD_b) - d(y, \tilde\dbD_b)$
, and \eqref{rVcont} can be estimated similarly as in Case 2.

{\it Case 4.} $y \in \dbD \backslash \tilde\dbD$. Then $\br_{\tilde \dbD}(y) -  \br_{\dbD}(y)  = d(y, \dbD_b) + d(y, \tilde\dbD_b)$.
First, 
\beaa
d(y, \tilde\dbD_b) = d(y, \tilde\dbD) \le d(\dbD, \tilde \dbD),
\eeaa
where the first equality is due to the fact that $y\notin  \tilde\dbD$ and the inequality is due to the fact that $y\in \dbD$. Next, it is obvious that
\beaa
d(y, \dbD_b) \le d(y, \tilde \dbD_b) + d( \tilde \dbD_b, \dbD_b)  \le 2d(\dbD, \tilde \dbD).
\eeaa
Put together with obtain \eqref{rVcont}.
\qed

We note that the first equality in \eqref{Dest} fails in general when $\dbD$ or
$\tilde \dbD$ is not convex. At below we present a simple example where $\lim_{\e\to 0}d(\dbD, \tilde\dbD^\e) = 0 < 1 = \lim_{\e\to  0}d(\dbD_b,\tilde\dbD_b^\e)$.
\begin{eg}
  For $\e>0$, set
  $\dbD := \{y\in \dbR^2: |y|\le 1\}$ and
  \begin{equation*}
    \tilde \dbD^\e := \{ r(\cos\th,\sin\th)\in \dbR^2\ :\
    0\leq r \leq 1,\ \ \e\leq \th\leq 2\pi-\e\} \subset \dbD.
  \end{equation*}
 One can easily see that
  $d(\dbD, \tilde\dbD^\e) = d((1,0), \tilde\dbD_b^\e) \to 0$, and  $d(\dbD_b, \tilde \dbD_b^\e) = d((0,0), \dbD_b) = 1$.
  \qed
\end{eg}

\ss

 \no{\bf Proof of Lemma \ref{lem-boundary}.}  Denote 
\beaa
\t_\d := \inf\{t\ge 0:  |\br (t, X^\a_t, Y^\a_t)|  > \d\}\wedge T_0,
\eeaa
 and consider the linear BSDE with solution pair $(\k, \b)$:
\beaa
\label{kbBSDE}
\k_t =  \1_{\{\t_\d < T_0\}} - \int_t^{T_0}\beta_sdB_s,\q 0\le t\le T_0,
\eeaa
where $\k\in \dbR$, $\b\in \dbR^{1\times d}$. It is clear that $|\k|\le 1$, $\beta_t = 0$ for $\t_\d \leq t \leq T_0$, and $\int_0^\cd \b_s dB_s$ is an BMO martingale. Thus, 
\beaa
\label{bBMO}
\dbE\Big[\exp\big(c_0 \int_0^{T_0} |\b_t|^2dt\big)\Big] \le C_0 < \infty,\q\mbox{for some $c_0, C_0>0$}.
\eeaa
For $n\ge 1$, we truncate $\b$ by $n$ and denote it as $\b^n$. Define
\beaa
\k^n_t := \k_0 + \int_0^t \b^n_s dB_s.
\eeaa
Then it is obvious that, for any $p\ge 1$,
\begin{equation}
  \label{knconv}
  \dbE\big[\sup_{0\le t\le T_0}|\k^n_t|^p\big]\le C_p<\infty;
  ~\mbox{and}~
  c^n_p:= \Big(\dbE\big[\sup_{0\le t\le T_0} |\k^n_t - \k_t|^p\big]\Big)^{1\over p}
  \to 0,~\mbox{as}~n\to\infty.  
\end{equation}

Introduce two random fields $\nu^n(t,\o, y)$ and  $\rho^n(t,\o, y)$:
\beaa
\dis \nu^n_{ij}(t, y) &:=& -\sum_{k=1}^m \int_0^1 \td_{z_{kj}}f^i(t, X^\a_t, Y^\a_t, Z^\a_t + \th y\b^n_t, \a_t) y_k d\th;\\
\dis \rho^n_{i}(t, y) &:=& - \int_0^1\Big[ \td_yf^i \cd y + \tr((\pa_z f^i)^\top \nu^n(t,y))\Big]\\
&&\dis\big(t, X^\a_t, Y^\a_t + \th \k^n_t y, Z^\a_t +  y\b^n_t + \th \k^n_t \nu^n(t,y), \a_t\big)  d\th;
\eeaa
where $ \nu^n = [\nu^n_{ij}]_{1\le i\le m, 1\le j\le d}\in \dbR^{m\times d}$, $\rho^n=[\rho^n_i]_{1\le i\le m}\in \dbR^m$. One can easily verify that
\beaa
\left.\ba{lll}
\dis \nu^n(t, y)\b^n_t :=f(t, X^\a_t, Y^\a_t, Z^\a_t, \a_t) -  f(t, X^\a_t, Y^\a_t, Z^\a_t + y\b^n_t, \a_t);\ss\\
\dis \k^n_t \rho^n(t, y) :=  f(t, X^\a_t, Y^\a_t, Z^\a_t + y\b^n_t,\a_t) \ss\\
\dis\qq\qq\q -  f(t, X^\a_t, Y^\a_t + \k^n_t y, Z^\a_t + y\b^n_t  + \k^n_t \nu(t, y),\a_t).
\ea\right.
\eeaa
Moreover, by Assumption \ref{assum-coefficients} (ii) we have
\bea
\label{nurhoest}
\left.\ba{c}
\dis |\nu^n(t, y)| \le C|y|,\q  |\rho^n(t, y)|  \le C|y|;\\
\dis |\nu^n(t, y)-\nu^n(t, \tilde y)| \le C_n\big[1+|y| \big] |y-\tilde y|,\\
\dis  |\rho^n(t, y)-\rho^n(t, \tilde y)| \le C_n\big[1+|y|+|\k^n_t| +  |\k^n_t||y|^2\big] |y-\tilde y|.
\ea\right.
\eea
Next, consider the following SDE:
\beaa
    \label{etaSDE} 
    \eta^n_t := \l \bn_\dbV(0,x_0,\pi(0,x_0,Y_0^{\a})) + \int_0^t\rho^n(s, \eta^n_s)ds + \int_0^t \nu^n(s, \eta^n_s) dB_s, 
    \eeaa
 where $\l>0$ is a small number which will be determined later. By the standard stopping arguments for stochastic Lipschitz continuous coefficients, and by the uniform linear growth in the first line of \eqref{nurhoest}, the above SDE is wellposed, and for any $p\ge 1$, 
 \bea
 \label{etaest}
 \dbE\Big[\sup_{0\le t\le T_0} |\eta^n_t|^p\Big] \le C_p |\l|^p,
 \eea
 where $C_p$ does not depend on $n$.
 
 Denote
\beaa 
\tilde Y^n_t := Y_t^{\a} + \k^n_t\eta^n_t,\q \tilde Z^n_t :=  Z_t + \k^n_t \nu^n_t +  \eta^n_t\b^n_t.
\eeaa
Then by the standard It\^{o} formula we have
\beaa
    \label{BSDEtildeY} 
    \tilde Y^n_t = Y_{\t_\d}^{\a} + \k^n_{\t_\d}\eta^n_{\t_\d} + \int_t^{\t_\d} f(s,X_s^{\a},\tilde Y^n_s,\tilde Z^n_s, \a_s)ds - \int_t^{\t_\d}\tilde Z^n_s dB_s. 
    \eeaa
    Moreover,  introduce  the BSDE 
    \beaa 
    \hat Y^n_t = Y^\a_{\t_\d} + \k^n_{\t_\d}\eta^n_{\t_\d}\ind{|\k^n_{\t_\d}| \le 2, |\eta^n_{\t_\d}|<{\d\over 2}} + \int_t^{\t_\d}f(s,X_s^{\a},\hat Y^n_s,\hat Z^n_s,\a_s)ds - \int_t^{\t_\d}\hat Z^n_s dB_s.
     \eeaa 
     By \eqref{knconv}, \eqref{etaest}, and noting that $|\k_{\t_\d}|\le 1$, it follows from standard BSDE estimates that
     \beaa 
     |\hat Y^n_0 - \tilde Y^n_0|^2 &\leq& C\dbE\Big[\big|\k^n_{\t_\d}\eta^n_{\t_\d}\big|^2 \big[\ind{|\k^n_{\t_\d}-\k_{\t_\d}|>1} + \ind{|\eta^n_{\t_\d}| \geq {\d\over 2}}\big]\Big]\\
     &\le& C\dbE\Big[\big|\k^n_{\t_\d}\eta^n_{\t_\d}\big|^2 |\k^n_{\t_\d}-\k_{\t_\d}| + {1\over \d^2}|\k^n_{\t_\d}|^2|\eta^n_{\t_\d}|^4\Big] \leq Cc^n_{2}+ C{|\l|^4\over \d^2}.
      \eeaa 
      Thus
      \bea
      \label{hatYntilde}
       |\hat Y^n_0 - \tilde Y^n_0| \le  C\big[\sqrt{c^n_{2}}+{|\l|^2\over \d}\big].
      \eea

Note that 
\beaa
\hat Y^n_{\t_\d} = Y^\a_{\t_\d} + \k^n_{\t_\d}\eta^n_{\t_\d}\ind{|\k^n_{\t_\d}| \le 2, |\eta^n_{\t_\d}|<{\d\over 2}}.
\eeaa
On $\{{|\k^n_{\t_\d}| \le 2, |\eta^n_{\t_\d}|<{\d\over 2}}\}^c$, we have $\hat Y^n_{\t_\d} = Y^\a_{\t_\d}\in \dbV(\t_\d, X^\a_{\t_\d})$. On $\{{|\k^n_{\t_\d}| \le 2, |\eta^n_{\t_\d}|<{\d\over 2}}\}$, noting again that $Y^\a_{\t_\d}\in \dbV(\t_\d, X^\a_{\t_\d})$, we have $\br_\dbV(\t_\d, X^\a_{\t_\d}, Y^\a_{\t_\d}) = -\d$ and $|\k^n_{\t_\d}\eta^n_{\t_\d}|\le \d$, then $\hat Y^n_{\t_\d}  = Y^\a_{\t_\d} +\k^n_{\t_\d} \eta^n_{\t_\d}  \in \dbV(\t_\d, X^\a_{\t_\d})$. So in both cases $\hat Y^n_{\t_\d} \in \dbV(\t_\d, X^\a_{\t_\d})$. Then by DPP \eqref{DPP} we have $\hat Y^n_0 \in \dbV(0, x_0)$. Thus, by \eqref{hatYntilde},
\beaa
\br_\dbV(0, x_0, \tilde Y^n_0) \le |\tilde Y^n_0 - \hat Y^n_0|  \le C\big[\sqrt{c^n_2}+{|\l|^2\over \d}\big].
\eeaa
On the other hand,  note that $\tilde Y^n_0 = Y^\a_0 + \k_0\l\bn_\dbV(0,x_0,\pi(0,x_0,Y_0^{\a}))$, for $\k_0\l$ small we have 
\beaa
\br_\dbV(0, x_0, \tilde Y^n_0) = \br_\dbV(0, x_0, Y^\a_0) + \k_0\l.
\eeaa
Thus
\beaa
\k_0\l = \br_\dbV(0, x_0, \tilde Y^n_0) - \br_\dbV(0, x_0, Y^\a_0) \le  C\big[\sqrt{c^n_2}+{|\l|^2\over \d}\big] + \e.
\eeaa
Send $n\to \infty$ and set $\l := \sqrt{\e\d}$, we obtain \eqref{boundaryest}:
\beaa
\dbP(\t_\d < T_0) = \k_0 \le C{\l\over \d} + {\e\over \l} = C\sqrt{\e\over \d}.
\eeaa

Finally, if $Y^\a_0\in \dbV_b(0, x_0)$, then $\e=0$. We see that
$\dbP(\t_\d < T_0) = 0$ for all $\d>0$ and all
$T_0<T$. This implies immediately that $Y^\a_t \in \dbV_b(t, X^\a_t)$,
$0\le t<T$, a.s. Moreover, note that $Y^\a_T = g(X^\a_T)$ and
$\dbV(T,x)=\{g(x)\}$, we have $Y^\a_T \in \dbV_b(T, X^\a_T)$ as well.
\qed

\ms
\no{\bf Proof of Example \ref{eg-classical}.} As usual we drop the subscript $_\dbV$ in $\br$ and $\bn$. 

(i) We first show that $\dbV\in C^{1,2}_0([0, T)\times \dbR;\cD^2_2)$. Fix $\d>0$ and denote $T_\d:=T-\d$. By Example \ref{eg-paV} we have, with $u(t,x) = T-t\ge \d$ there for $t\in [0, T_\d]$,
\beaa
\br (t,x,y) =  |y-w(t,x)| - (T-t).
\eeaa
Then it is clear that $\dbV\in C^{1,2}([0, T_\d]\times \dbR;\cD^2_2)$. By \eqref{circle}, for $|y-w(t,x)|=T-t$, we have
\beaa
 \left.\ba{c}
\dis \bn = {y-w\over T-t};\q \pa_t \dbV = \big[ \td_t w \cd \bn-1\big] \bn;\q \pa_{x} \dbV = \big[ \td_{x} w \cd \bn\big] \bn;\ms\\
\dis \pa_{x}\bn = {1\over T-t}\big[- \td_{x}w + [\bn\cd\td_{x}w] \bn\big],\q \pa_{y} \bn = {1\over T-t}\big[I_{2\times 2}- \bn\bn^\top\big];\ms\\
\dis \pa_{xx} \dbV = -{1\over T-t}\Big[\big[ |\td_{x }w|^2- \td_{xx}w \cd \bn + |\td_{x} w\cd\bn|^2 \big] \bn +[\td_{x}w\cd \bn]\big[\td_{x}w - (\td_{x}w\cd\bn)\bn\big]\Big].
\ea\right.
\eeaa
In particular, $c_{T_0} = {1\over T}$ in \eqref{HVbound}, and thus  $\dbV\in C^{1,2}_0([0, T)\times \dbR;\cD^2_2)$.

(ii) We next verify the conditions in Theorem \ref{thm-unique} (ii). For any $a\in A$ and $\zeta\in \dbT_\dbV(t,x,y)$, by \eqref{H} we have: at $(t,x,y)\in \dbG_\dbV$,
\beaa
&&\dis h^0_\dbV(t, x,y, \pa_x\dbV, \pa_{xx}\dbV, a,\zeta)  =\frac{1}{2}\pa_{xx}\dbV- \Big[\zeta\cd \pa_x\bn + \frac{1}{2}\zeta^\top \pa_y\bn\zeta\Big] \bn\\
&&\dis = -{1\over 2(T-t)}\Big[\big[ |\td_{x }w|^2- \td_{xx}w \cd \bn + |\td_{x} w\cd\bn|^2 \big] \bn +[\td_{x}w\cd \bn]\big[\td_{x}w - (\td_{x}w\cd\bn)\bn\big]\Big]\\
&&\dis\q - {1\over T-t}\Big[- \zeta\cd \td_x w + {1\over 2}|\zeta|^2\Big] \bn.
\eeaa
Thus
\beaa
&&\dis \bn\cd h_\dbV(t, x,y, \pa_x\dbV, \pa_{xx}\dbV, a,\zeta) = \bn\cd h^0_\dbV(t, x,y, \pa_x\dbV, \pa_{xx}\dbV, a,\zeta) + \bn \cd [f^0(t,x) +a]\\
&&\dis = -{1\over 2(T-t)}\Big[ |\td_{x }w|^2- \td_{xx}w \cd \bn + |\td_{x} w\cd\bn|^2  - 2 \zeta\cd \td_x w +|\zeta|^2\Big] +\bn\cd\big[f^0+a\big]\\
&&\dis = -{1\over 2(T-t)}\Big[ |\zeta-\td_{x }w|^2- \td_{xx}w \cd \bn + |\td_{x} w\cd\bn|^2 \Big] +\bn\cd\big[f^0+a\big].
\eeaa
Recall $|a|\le 1$, then clearly  the optimal arguments are:
\beaa
a^* = I^\dbV_1(t,x,y) := \bn(t,x,y),\q \zeta^* = I^\dbV_2(t,x,y) := \td_x w - [\td_x w\cd\bn]\bn.
\eeaa
Together with \eqref{xi*}, this implies further that
\beaa
\dis &&\tilde I^\dbV_3 := -\Big[\pa_t \dbV+h^0_\dbV(\cd, \pa_x \dbV, \pa_{xx}\dbV, I^\dbV_1, I^\dbV_2) + f^0 + I^\dbV_1\Big]\\
\dis &&=  -\big[ \td_t w \cd \bn-1\big] \bn + {1\over T-t}\big[- \zeta\cd \td_x w + {1\over 2}|\zeta|^2\big] \bn -[ f^0 + \bn]\\
\dis &&\q  + {1\over 2(T-t)}\Big[\big[ |\td_{x }w|^2- \td_{xx}w \cd \bn + |\td_{x} w\cd\bn|^2 \big] \bn +[\td_{x}w\cd \bn]\big[\td_{x}w - (\td_{x}w\cd\bn)\bn\big]\Big]; \\
\dis && I^\dbV_3 := -f^0 + [\bn\cd f^0] \bn +  {1\over 2(T-t)}[\td_{x}w\cd \bn]\big[\td_{x}w - (\td_{x}w\cd\bn)\bn\big].
\eeaa
Plug $I^\dbV_1, I^\dbV_2, I^\dbV_3$ into \eqref{X*}, clearly the resulted SDE is wellposed.
\qed

\ms

\no{\bf Proof of Example \ref{eg-heat}.} We now compute the equation \eqref{HJBrhat} in this case. First, 
\beaa
\cN(t,x,y, \td_x \eta, \td_y\eta) &=& \big\{(a, z): [1, z]^\top[\td_x \eta + z\td_y \eta] = 0\big\} \\
&=& \big\{(a, z): z = -{\td_x \eta \over \td_y \eta}\big\} = \big\{\big(a, -{\td_x \wh \br (t,x,y) \over \td_y \wh \br (t,x,y)}\big)\big\}.
\eeaa
Then, recalling $\eta = {1\over 2}|\wh \br |^2$, 
\beaa
F &=& {1\over 2}\Big[\td_{xx}\eta - 2\td_{xy}\eta {\td_x \wh \br \over \td_y\wh \br } +  \td_{yy}\eta\big|{\td_x \wh \br \over \td_y \wh \br }\big|^2\Big](t,x,y)\\
&=&{1\over 2}\Big[\wh \br \td_{xx} \wh \br + |\td_x\wh \br |^2 - 2[\wh \br \td_{xy}\wh \br + \td_x \wh \br \td_y \wh \br ]{\td_x \wh \br \over \td_y\wh \br } + [\wh \br \td_{yy}\wh \br + |\td_y\wh \br |^2]\big|{\td_x \wh \br \over \td_y \wh \br }\big|^2\Big](t,x,y)\\
&=&{1\over 2} \wh \br \Big[ \td_{xx} \wh \br  - 2 \td_{xy}\wh \br {\td_x \wh \br \over \td_y\wh \br } +  \td_{yy}\wh \br \big|{\td_x \wh \br \over \td_y \wh \br }\big|^2\Big](t,x,y). 
\eeaa
Note further that $\wh \br = 0$ on $\dbG_\dbV$. Then, for $(t,x,y)\in \dbG_\dbV$, we have
\beaa
\td_x F &=& {1\over 2}\td_x\wh \br \Big[ \td_{xx} \wh \br  - 2 \td_{xy}\wh \br {\td_x \wh \br \over \td_y\wh \br } +  \td_{yy}\wh \br \big|{\td_x \wh \br \over \td_y \wh \br }\big|^2\Big](t,x,y);\\
\td_y F &=& {1\over 2}\td_y\wh \br \Big[ \td_{xx} \wh \br  - 2 \td_{xy}\wh \br {\td_x \wh \br \over \td_y\wh \br } +  \td_{yy}\wh \br \big|{\td_x \wh \br \over \td_y \wh \br }\big|^2\Big](t,x,y).
\eeaa
On the other hand, note that $\td_t \eta = \wh \br \td_t \wh \br $. Then, again at $(t,x,y)\in \dbG_\dbV$,
\beaa
\td_x \td_t \eta = \td_x \wh \br \td_t \wh \br ,\q \td_y \td_t \eta = \td_y \wh \br \td_t \wh \br .
\eeaa
Plug these into \eqref{HJBrhat}, we obtain \eqref{heathat} immediately.
\qed
\end{appendix}

\bibliographystyle{plain}
\bibliography{SetPDE}

\end{document}